\documentclass{amsart}
\usepackage{amsfonts,amssymb,amscd,amsmath,enumerate,verbatim}

\usepackage[all]{xy}
\SelectTips{eu}{} 
\xyoption{curve}

\newcommand{\lra}{\longrightarrow}

\newcommand{\xra}{\xrightarrow}

\newcommand{\up}[1]{{{}^{#1}\!}}
\newcommand{\card}{\operatorname{card}}

\newcommand{\wh}{\widehat}
\newcommand{\wt}{\widetilde}
\newcommand{\twedge}{\textstyle{\bigwedge}}
\newcommand{\les}{{\scriptscriptstyle\leqslant}}
\newcommand{\ges}{{\scriptscriptstyle\geqslant}}

\newcommand{\dcat}{{\mathcal D}}
\newcommand{\dtensor}[3]{{#1}\otimes_{#2}^{\bf{L}}{#3}}
\newcommand{\koszul}[2]{{\operatorname{K}[#1;#2]}}

\newcommand{\rkoszul}[2]{{\operatorname{K}^{#1}[#2]}}
\newcommand{\Rhom}[3]{{\bf{R}\!\Hom}_{#1}(#2,#3)}
\newcommand{\RHom}[3]{{\bf{R}\!\Hom}_{#1}\big(#2,#3\big)}
\newcommand{\shift}{{\sf\Sigma}}

\newcommand\dd{\partial}
\newcommand\HH{\operatorname{H}}

\newcommand\Tor{\operatorname{Tor}}
\newcommand\Ext{\operatorname{Ext}}
\newcommand\Hom{\operatorname{Hom}}

\newcommand\Ker{\operatorname{Ker}}
\newcommand\Coker{\operatorname{Coker}}
\newcommand\Image{\operatorname{Im}}
\newcommand\Soc{\operatorname{Soc}}
\newcommand\Spec{\operatorname{Spec}}
\newcommand\Supp{\operatorname{Supp}}

\newcommand{\astar}{a^*}
\newcommand{\grass}[2]{\operatorname{gr}_{#1}(#2)}
\newcommand{\lch}[3]{\HH^{#1}_{#2}(#3)}

\newcommand\codim{\operatorname{codim}}
\newcommand\depth{\operatorname{depth}}
\newcommand\edim{\operatorname{edim}}
\newcommand\gd{\operatorname{G-dim}}
\newcommand\height{\operatorname{height}}
\newcommand{\lol}{{\ell\ell}}
\newcommand{\ord}{\operatorname{ord}}
\newcommand\rad{\operatorname{rad}}
\newcommand\rank{\operatorname{rank}}
\newcommand\serre{\operatorname{spread}}
\newcommand\type{\operatorname{type}}

\newcommand{\hann}[2]{\operatorname{Ann}_{\dcat(#1)}(#2)}
\newcommand{\hell}[2]{\lol_{\dcat(#1)}(#2)}
\newcommand{\Ann}{\operatorname{Ann}}

\newcommand\idmap{\operatorname{id}}

\newcommand{\rP}[2]{P^{#1}_{#2}(t)}
\newcommand{\rI}[2]{I_{#1}^{#2}(t)}
\newcommand{\lP}[2]{F^{#1}_{#2}(t)}

\newcommand{\betti}[3]{\beta_{#1}^{#2}(#3)}
\newcommand{\bass}[3]{\mu^{#1}_{#2}(#3)}

\newcommand\cxy{\operatorname{cx}}
\newcommand\curv{\operatorname{curv}}
\newcommand\fd{\operatorname{flat\,dim}}
\newcommand\id{\operatorname{inj\,dim}}
\newcommand\injcxy{\operatorname{inj\,cx}}
\newcommand\injcurv{\operatorname{inj\,curv}}
\newcommand\pd{\operatorname{proj\,dim}}

\newcommand\fa{{\mathfrak a}}
\newcommand\fb{{\mathfrak b}}
\newcommand\fl{{\mathfrak l}}
\newcommand\fm{{\mathfrak m}}
\newcommand\fn{{\mathfrak n}}

\newcommand\fp{{\mathfrak p}}
\newcommand\fq{{\mathfrak q}}
\newcommand\fr{{\mathfrak r}}

\newcommand\bsf{{\boldsymbol f}}

\newcommand\bsu{{\boldsymbol u}}
\newcommand\bsv{{\boldsymbol v}}
\newcommand\bsw{{\boldsymbol w}}
\newcommand\bsx{{\boldsymbol x}}
\newcommand\bsy{{\boldsymbol y}}
\newcommand\bsz{{\boldsymbol z}}

\newcommand\BN{{\mathbb N}}
\newcommand\BQ{{\mathbb Q}}
\newcommand\BR{{\mathbb R}}
\newcommand\BZ{{\mathbb Z}}

\newcommand{\vf}{{\varphi}}

\newcommand{\var}{{\hskip1pt\vert\hskip1pt}}

\theoremstyle{plain}

\newtheorem{theorem}{Theorem}[section]

\newtheorem{proposition}[theorem]{Proposition}
\newtheorem{lemma}[theorem]{Lemma}

\newtheorem{corollary}[theorem]{Corollary}

\newtheorem{subtheorem}{Theorem}[subsection]
\newtheorem{subproposition}[subtheorem]{Proposition}
\newtheorem{sublemma}[subtheorem]{Lemma}
\newtheorem{subcorollary}[subtheorem]{Corollary}

{\Alph{theorem}}

\theoremstyle{definition}

\newtheorem{subexample}[subtheorem]{Example}
\newtheorem{subconstruction}[subtheorem]{Construction}

\newtheorem{example}[theorem]{Example}

\newtheorem{subremark}[subtheorem]{Remark}
\newtheorem{remark}[theorem]{Remark}

\theoremstyle{remark}

\newtheorem{chunk}[theorem]{}

\newtheorem{subquestion}[subtheorem]{Question}
\newtheorem{subchunk}[subtheorem]{}

\newtheorem*{Claim}{Claim}

\newtheorem*{Remark}{Remark}

\numberwithin{equation}{subtheorem}

\begin{document}

\title{Homology over local homomorphisms}
\date{\today}

\author[L.~L.~Avramov]{Luchezar L.~Avramov}
\address{Department of Mathematics,
University of Nebraska, Lincoln, NE 68588, U.S.A.}
\email{avramov@math.unl.edu}

\author[S.~Iyengar]{Srikanth Iyengar}
\address{Department of Mathematics,
University of Nebraska, Lincoln, NE 68588, U.S.A.}
\email{iyengar@math.unl.edu}

\author[C.~Miller]{Claudia Miller}
\address{Department of Mathematics, Syracuse University,
Syracuse, NY 13244, U.S.A.}
\email{cmille06@syr.edu}

\thanks
{The authors were partly supported by grants from the National Science
Foundation.\endgraf
Part of the work on this paper was done during the Spring semester of 2003
at MSRI, where the authors participated in the Program in Commutative
Algebra, and during visits of C.M.\ (August 2002--December 2002) and
S.I.\ (August 2003--May 2004) to the University of Nebraska, Lincoln.
The hospitality of both institutions is gratefully acknowledged.}

\begin{abstract}
The notions of Betti numbers and of Bass numbers of a finite module $N$
over a local ring $R$ are extended to modules that are only assumed to
be finite over $S$, for some local homomorphism $\vf\colon R\to S$.
Various techniques are developed to study the new invariants and to
establish their basic properties.  In several cases they are computed
in closed form.  Applications go in several directions.  One is to
identify new classes of finite $R$-modules whose classical Betti
numbers or Bass numbers have extremal growth.  Another is to transfer
ring theoretical properties between $R$ and $S$ in situations where
$S$ may have infinite flat dimension over $R$.  A third is to obtain
criteria for a ring equipped with a `contracting' endomorphism---such
as the Frobenius endomorphism---to be regular or complete intersection;
these results represent broad generalizations of Kunz's characterization
of regularity in prime characteristic.
 \end{abstract}


\maketitle


\section*{Introduction}

The existence of a homomorphism $\vf\colon R\to S$ of commutative noetherian rings does
not imply a relationship between ring theoretical properties of $R$ and $S$, such as
regularity, normality, Cohen-Macaulayness, etc.  It is therefore remarkable that certain
homological conditions on the $R$\emph{-module} $S$ force stringent relations between the
\emph{ring} structures of $R$ and $S$.  A classical chapter of commutative algebra,
started by Grothendieck, deals with the case when $S$ is \emph{flat} over $R$.  Parts of
this theory have been extended to a situation where $S$ is only assumed to have finite
flat dimension over $R$.

An initial motivation for this investigation was to find conditions on the $R$-module $S$
that allow a transfer of properties between the rings even in cases of infinite flat
dimension.  It became rapidly apparent that such a program requires new invariants.  Our
first objective is to introduce homological measures for finite $S$-modules, which reflect
their structure as $R$-modules. In the special case when $\vf$ is the identity map of $R$,
they reduce to classical invariants of finite $R$-modules.  In general, they have
properties that adequately extend those of their counterparts in the finite case.  Our
main goal is to demonstrate the usefulness of the new concepts for studying homomorphisms
of commutative noetherian rings.  The central case is when the homomorphism $\vf$ is
local, which means that the rings $R$ and $S$ are local and $\vf$ maps
the unique maximal ideal of $R$ into that of $S$.

We start by constructing sequences of invariants modeled on the sequences of integers, the
Betti numbers and the Bass numbers, classically attached to a finite $R$-module $M$. One
way to introduce them is as ranks of modules in minimal resolutions. Another is as ranks
of the vector spaces $\Tor^R_n(k,M)$ or $\Ext^n_R(k,M)$, where $k$ is the residue field of
$R$.  However, for a finite $S$-module $N$ neither approach provides finite numbers in
general.

We define Betti numbers $\betti n{\vf}N$ and Bass numbers $\bass n{\vf}N$ of $N$ over
$\vf$, in a way that ensures their finiteness when $N$ is finite over $S$.  To this end we
use the fact that $\Tor^R_n(k,N)$ or $\Ext^n_R(k,N)$ have natural structures of finite
$S$-modules, and that this holds even when $N$ is a homologically finite complex of
$S$-modules.  This is the contents of Section \ref{Sequences}.  The necessary machinery is
assembled in the first three sections of the paper. It is put to different use in Section
\ref{Dimensions} where it is applied in conjunction with the ``Bass conjecture'' to prove
that if a finite $S$-module has injective dimension over $R$, then $R$ is Cohen-Macaulay.

Under special conditions we compute in closed form the
entire sequence of Betti numbers or of Bass numbers.  Results are often
best stated in terms of the corresponding generating function, called the
Poincar\'e series or the Bass series of $N$ over $\vf$, respectively.
Section \ref{Illustrations} contains instances of such computations,
intended both as illustration and for use later in the paper.

We start Section \ref{Separation} by establishing upper bounds for the Poincar\'e series
and the Bass series of $N$ over $\vf$, in terms of expressions where the contributions of
$R$, $S$, and $\vf$ appear as separate factors.  When the bound for Poincar\'e series is
reached the module $N$ is said to be separated over $\vf$; it is said to be injectively
separated over $\vf$ when the Bass series reaches its bound. We study such modules in
significant detail.  Using numerical invariants of the $S$-module $N$ obtained from a
Koszul complex and analyzed in Section \ref{Koszul invariants}, we prove that the Betti
numbers of separated modules share many properties with the classical Betti numbers of $k$
over $R$.  It came as a surprise (to us) that separated modules occur with high frequency.
For example, when the ring $S$ is regular every $S$-module is separated over $\vf$.

As we do not assume the homological dimensions of $N$ over $R$ to be finite, the sequences
of Betti numbers and of Bass numbers of $N$ over $\vf$ may contain a lot of inessential
information.  Some of our main results show that the asymptotic behavior of these
sequences captures important aspects of the structure of $N$.  Comparisons of Betti
sequences to polynomial functions and to exponential functions lead us in Section
\ref{Asymptotes} to the notions of complexity $\cxy_\vf N$ and curvature $\curv_\vf N$,
respectively; injective invariants are similarly derived from Bass numbers.

The next four sections are devoted to the study of various aspects of these new
invariants.  In Section \ref{Comparisons} we analyze their dependence on $S$, and prove
that this ring can be replaced by any other ring $S'$ over which $N$ is finite module,
provided its action is compatible with that of $S$.  In particular, it follows that when
the $R$-module $N$ is finite, then its complexities and curvatures over $\vf$ are equal to
those over $R$, although the corresponding Betti numbers or Bass numbers may differ
substantially.  In Section \ref{Composition} we investigate changes of (injective)
complexities and curvatures under compositions of homomorphisms.  In Section
\ref{Localization} 
we prove that they do not go up under localization.  In Section
\ref{Extremality} we give conditions that ensure that the asymptotic invariants of $N$
over $\vf$ have the maximal possible values, namely, those of the corresponding invariants
of $k$ over $R$.

In Section \ref{Endomorphisms} we focus on the case when $\vf$ is a
contracting endomorphism of $R$, by which we mean that for every non-unit
$x\in R$ the sequence $(\vf^i(x))$ converges to $0$ in the natural
topology of $R$.  The motivating example is the Frobenius endomorphism,
but interesting contractions exist in all characteristics.

The final Section \ref{Local homomorphisms} contains applications of the
methods developed in the paper to the study of local homomorphisms. We
obtain results on the descent of certain ring theoretical properties
from $S$ to $R$.  We also prove that when $\vf$ is a contraction its
homological properties determine whether $R$ is regular or complete
intersection.  Thus, we obtain vast generalizations and a completely
new proof of Kunz's famous characterization of regularity in prime
characteristic.  Some of our results concerning the Frobenius endomorphism
are announced in Miller's survey \cite{Mi} of the homological properties
of that map.

Even when dealing with modules, in some constructions and in many
proofs we use complexes.  We have therefore chosen to develop the entire
theory in terms of the derived categories of $R$ and $S$, stating both
definitions and results for homologically finite complexes of modules
over $S$, rather than just for finite $S$-modules.

\section*{Notation}

Throughout this paper $(R,\fm,k)$ denotes a \emph{local ring}:
this means  (for us) that $R$ is a commutative noetherian ring with unique maximal
ideal $\fm$ and residue field $k=R/\fm$.  We also fix a \emph{local
homomorphism}
 \[
\vf\colon (R,\fm,k)\lra (S,\fn,l)
 \]
that is, a homomorphism of noetherian local rings such that
$\vf(\fm)\subseteq\fn$.

As usual, $\edim S$ stands for the \emph{embedding dimension} of $S$,
defined as the minimal number of generators of its maximal ideal $\fn$. In
addition, we set
 \[
\edim\vf=\edim(S/\fm S)
 \]
A \emph{set of generators of $\fn$ modulo $\fm S$} is a finite subset
$\bsx$ of $\fn$ whose image in $S/\fm S$ generates the ideal
$\fn/\fm S$.  Such a set is \emph{minimal} if no proper subset of
$\bsx$ generates $\fn$ modulo $\fm S$; by Nakayama's Lemma, this
happens if and only if $\card(\bsx)=\edim\vf$.

Throughout the paper, $N$ denotes a complex of $S$-modules
 \[
\cdots\lra N_{n+1}\xra{\dd^N_{n+1}} N_{n}\xra{\ \dd^N_{n}\ }
N_{n-1}\lra\cdots
 \]
which is \emph{homologically finite}, that is, the $S$-module $\HH_n(N)$
is finite for each $n$ and vanishes for almost all $n\in\BZ$.  

Complexes of $S$-modules are viewed as complexes of $R$-modules via $\vf$.
Modules are identified with complexes concentrated in degree $0$.

We let $\wh S$ denote the $\fn$-adic completion of $S$ and set $\wh
N=N\otimes_S\wh S$.  This complex is homologically finite over $\wh S$:
the flatness of $\wh S$ over $S$ yields $\HH_n(\wh N)\cong
\HH_n(N)\otimes_S\wh S$.

\section{Complexes}
\label{Complexes}
This article deals mainly with homologically finite complexes. However, it is often
convenient, and sometimes necessary, to operate in the full derived category of complexes.

\subsection{Derived categories}
\label{derived categories}
For each complex $M$, we set 
 \[
\sup M=\sup\{n\in\BZ\mid M_n\ne0\}
\quad\text{and}\quad
\inf M=\inf\{n\in\BZ\mid M_n\ne0\}
 \]
 We say that $M$ is \emph{bounded} when both numbers are finite and that it is
 \emph{homologically bounded} when $\HH(M)$ is bounded.

Let $\dcat(R)$ denote the derived category of complexes of $R$-modules, obtained from the
homotopy category of complexes of $R$-modules by localizing at the class of homology
isomorphisms. The procedure for constructing this localization is the same as for derived
categories of bounded complexes, see Verdier \cite{Ve}, once appropriate ``resolutions''
are provided, such as the $K$-resolutions of Spaltenstein \cite{Spa}. 

The symbol $\simeq$ denotes isomorphism in a derived category, and $\shift$ the shift
functor.  We identify the category of $R$-modules with the full subcategory of $\dcat(R)$
consisting of the complexes with homology concentrated in degree zero, and let $\dcat^{\rm
  f}(R)$ denote the full subcategory of homologically finite complexes.

The derived functors of tensor products and of homomorphisms are denoted $\dtensor
{(-}R{-)}$ and $\Rhom R--$, respectively. These may be defined as follows: for complexes
of $R$-modules $X$ and $Y$, let $P$ be a $K$-projective resolution of $X$, and set
\[
\dtensor XRY\simeq P\otimes_RY
\qquad\text{and}\qquad
\Rhom RXY\simeq\Hom_R(P,Y)
 \]
 In particular, when $Y$ is a complex of $S$-modules, so are $P\otimes_RY$ and
 $\Hom_R(P,Y)$.  In this way, $\dtensor -RN$ and $\Rhom R-N$ define functors from
 $\dcat(R)$ to $\dcat(S)$.

For each integer $n$ we set
 \[
\Tor^R_n(X,Y) = \HH_n(\dtensor XRY)
\qquad\text{and}\qquad
\Ext_R^n(X,Y
) = \HH_{-n}(\Rhom RXY)
\]
When $Y$ is a complex of $S$-modules, $\Tor^R_n(X,Y)$ and $\Ext_R^n(X,Y)$ inherit
$S$-module structures from $P\otimes_RY$ and $\Hom_R(P,Y)$, respectively.

\smallskip

The rest of the section is a collection of basic tools used frequently in the paper.

\subsection{Endomorphisms}
Let $X$ be a complex of $S$-modules. Let $\theta\colon X\to X$ be a morphism in $\dcat(S)$
and let $X\xra{\theta} X\to C\to $ be a triangle in $\dcat(S)$.

The homology long exact sequence yields

\begin{subchunk}
\label{koszul:les}
For each integer $n$, there exists an exact sequence of $S$-modules
 \[
0\lra\Coker\HH_n(\theta)\lra \HH_n(C)
\lra\Ker\HH_{n-1}(\theta)\lra0
  \]
 In particular, since $S$ is noetherian, if $\HH_n(X)$ is finite for each $n\in\BZ$
 (respectively, if $X$ is homologically finite), then $C$ has the corresponding property.
\end{subchunk}

The next statement is a slight extension of \cite[(1.3)]{FI}.

\begin{subchunk}
\label{koszul:sup}
The following inequalities hold:
 \[
\sup \HH(C) \le\sup \HH(X) +1
\qquad\text{and}\qquad
\inf \HH(X)\le\inf \HH(C)
 \]
In addition, if $\Image\HH(\theta)\subseteq\fn\HH(X)$, then
 \[
\sup \HH(X)\le\sup\HH(C)
\qquad\text{and}\qquad
\inf \HH(X)=\inf\HH(C)
 \]

Indeed, the (in)equalities follow from the exact sequences above;
under the additional hypothesis, Nakayama's Lemma has to be invoked
as well.
 \end{subchunk}

The symbol $\ell_S$ denotes length over $S$.

\begin{sublemma}
\label{rod:sandwich}
If each $\ell_S\HH_n(X)$ is finite, and $\HH(\theta)^v=0$ for some $v\in\BN$
then for each $n\in\BZ$ there are inequalities
 \[
v^{-1}\big(\ell_S\HH_n(X) + \ell_S\HH_{n-1}(X)\big)
\leq \ell_S \HH_n(C) \leq
\ell_S\HH_n(X) + \ell_S\HH_{n-1}(X)
 \]
\end{sublemma}

\begin{proof}
Set $\alpha_n=\HH_n(\theta)$.
A length count on the short exact sequence \eqref{koszul:les} yields
 \[
\ell_S\HH_n(C) =\ell_S\Coker(\alpha_n)+\ell_S\Ker(\alpha_{n-1})
 \]
Since $\alpha$ is nilpotent of degree $v$, each $S$-module
$\HH_n(X)$ has a filtration
 \[
0=\Image(\alpha_n^v)\subseteq \Image(\alpha_n^{v-1})\subseteq
\cdots\subseteq \Image(\alpha_n)\subseteq \Image(\alpha_n^0)=
\HH_n(X)
 \]
by $S$-submodules.  As $\alpha_n$ induces for each $i$ an epimorphism
 \[
\frac{\Image(\alpha_n^{i-1})}{\Image(\alpha_n^{i})}\lra
\frac{\Image(\alpha_n^{i})}{\Image(\alpha_n^{i+1})}
 \]
we get the inequality on the left hand side in the following formula
 \[
\frac 1v \ell_S\HH_n(X)\leq\ell_S\Coker(\alpha_n)=
\ell_S\Ker(\alpha_{n})\leq\ell_S\HH_n(X)
 \]
The inequality on the right is obvious.  The equality comes from the
exact sequence
 \[
0\lra\Ker(\alpha_{n})\lra \HH_n(X)\xra{\ \alpha_n\ }\HH_n(X)\lra
\Coker(\alpha_n)\lra 0
 \]
To get the desired inequalities put together the numerical
relations above.
 \end{proof}

\subsection{Actions}
Here is another way to describe the action of $S$ on homology; see \eqref{derived
categories}. We let $\lambda^X_s\colon X\to X$ denote multiplication
with $s\in S$ on a complex $X$.

\begin{subchunk}
\label{homotopies}
Recall that $\Ext^0_S(N,N)$ is the set of homotopy classes of $S$-linear
morphisms $G\to G$, where $G$ is a projective resolution of $N$ over
$S$.  Compositions of morphisms commute with the formation of homotopy
classes, so $\Ext^0_S(N,N)$ is a ring.  The map assigning to each $s\in S$
the morphism $\lambda_s\colon G\to G$ defines a homomorphism of rings
$\eta_S\colon S\to \Ext^0_S(N,N)$ whose image lies in the center of
$\Ext^0_S(N,N)$, and
\begin{equation*}
\xymatrixrowsep{3pc}
\xymatrixcolsep{3pc}
\xymatrix{ S
\ar@{->}[r]^-{\eta_S} &\Ext^0_S(N,N)
\ar@{->}[d]^{\Ext^0_{\vf}(N,N)}
\\
R
\ar@{->}[u]^{\vf}
\ar@{->}[r]^-{\eta_R}
&\Ext^0_R(N,N)
}
\end{equation*}
is a commutative diagram of homomorphisms of rings.  If $\sigma$ is
any morphism in the homotopy class $\Ext^0_{\vf}(N,N)\eta_S(s)$, then
the maps $\HH(\dtensor MR\sigma)$ and $\HH\Rhom RM\sigma$ coincide with
multiplication by $s$ on $\Tor^R(M,N)$ and $\Ext_R(M,N)$, respectively.
 \end{subchunk}

\begin{sublemma}
\label{newfinite}
If $\inf\HH(M)>-\infty$ and each $R$-module $\HH_n(M)$ is finite,
then all $S$-modules $\Tor^R_n(M,N)$ and $\Ext_R^n(M,N)$ are finite, and
are trivial for all $n\ll0$.
 \end{sublemma}

\begin{proof}
  Since $N$ is homologically finite over the noetherian ring $S$, it is isomorphic in the
  derived category of $S$-modules to a finite complex of finite $S$-modules.  Changing
  notation if necessary, we may assume that $N$ itself has these properties.  On the other
  hand, since $R$ is noetherian, one may choose a $K$-projective resolution $P$ such that
  for each $n$, the $R$-module $P_n$ is finite and projective and $P_n=0$ for
  $n<\inf\HH(M)$. The complexes $P\otimes_RN$ and $\Hom_R(P,N)$ then consist of finite
  $S$-modules, so the desired assertion follows.
 \end{proof}

\subsection{Annihilators}
\label{Annihilators}
The annihilator of an $S$-module $H$ is denoted $\Ann_S(H)$.

The \emph{homotopy annihilator} of $N$ over $S$, introduced by Apassov
\cite{Ap}, is the ideal $\hann SN=\Ker(\eta_S)$.  With $\rad(I)$
denoting the radical of $I$, the first inclusion below follows from
\eqref{homotopies} and the second is \cite[Theorem, \S2]{Ap}.

\begin{subchunk}
\label{hann}
There are inclusions $\hann SN\subseteq\Ann_S(\HH(N))
\subseteq\rad\big(\!\hann SN\big)$.
 \end{subchunk}

The first assertion below follows from \eqref{homotopies}, the
second from \eqref{hann}.

\begin{subchunk}
\label{hann2}
The ideal $\hann RM S + \hann SN$ is contained in the homotopy
annihilators of the complexes of $S$-modules $\dtensor MRN$ and $\Rhom
RMN$, and hence annihilates $\Tor^R_n(M,N)$ and $\Ext_R^n(M,N)$
for each $n\in\BZ$.
 \end{subchunk}

\begin{subchunk}
\label{primary}
If $\HH(N)/\fm\HH(N)$ is artinian over $S$ (in particular, if
$\HH(N)$ is finite over $R$, or if the ring $S/\fm S$ is artinian),
then the ideal $\fm S+\hann SN$ is $\fn$-primary.

Indeed, it follows from \eqref{hann} that the radical of $\fm S+\hann SN$
equals the radical of $\fm S + \Ann_S(\HH(N))$, which is $\fn$ in view
of the hypothesis.
 \end{subchunk}

\subsection{Koszul complexes}
Let $\bsx$ be a finite set of elements in $S$.
The \emph{Koszul complex} $\koszul\bsx S$ is the DG (= differential
graded) algebra with underlying graded algebra the exterior algebra
on a basis $\{e_x\colon |e_x|=1\}_{x\in\bsx}$ and differential given by
$\dd(e_x)=x$ for each $x\in \bsx$. Let $\koszul{\bsx}N$ be the
DG module $\koszul{\bsx}S\otimes_SN$ over $\koszul\bsx S$.

\begin{subchunk}
\label{free koszul}
In the derived category of $S$-modules there are isomorphisms
\begin{gather*}
\dtensor MR{\koszul\bsx N}
\simeq \dtensor MR{\dtensor {\koszul\bsx S}S{N}}
\simeq \koszul\bsx{\dtensor MRN}\\
\Rhom RM{\koszul\bsx N} \simeq \koszul\bsx{\Rhom RMN}
\end{gather*}
because $\koszul\bsx S$ is a finite free complex over $S$.
\end{subchunk}

\begin{subchunk}
\label{koszul dual}
There is an isomorphism of complexes of $S$-modules
 \[
\Hom_S(\koszul{\bsx}S,N) \cong \shift^{-\card(\bsx)} \koszul{\bsx}N
 \]

Indeed, the canonical morphism of complexes of $S$-modules
 \[
\Hom_S(\koszul{\bsx}S,S)\otimes_SN\to \Hom_S(\koszul{\bsx}S,N)
 \]
is bijective, since $\koszul{\bsx}S$ is a finite free complex over $S$,
and the self-duality of the exterior algebra yields an isomorphism
$\Hom_S(\koszul{\bsx}S,S)\cong\shift^{-\card(\bsx)}\koszul{\bsx}S$.
 \end{subchunk}

\begin{subchunk}
\label{hann koszul}
The homotopy annihilator of $\koszul{\bsx}N$ contains $\bsx S+\hann SN$.

\smallskip
Indeed, the Leibniz rule for the DG module $K=\koszul{\bsx}N$ shows that
left multiplication with $e_x\in\koszul{\bsx}S$ on $K$ is a homotopy between
$\lambda^K_x$ and $0$, so $\hann S{K}$ contains $\bsx$.
If follows from the definitions that it also contains $\hann SN$.
 \end{subchunk}

It is known that Koszul complexes can be described as iterated mapping cones:

\begin{subchunk}
\label{cone koszul}
For each $x\in\bsx$, the Koszul complex $\koszul{\bsx}N$ is isomorphic to
the mapping cone of the morphism $\lambda^{\koszul{\bsx'}N}_{x}$,
where $\bsx'=\bsx\smallsetminus\{x\}$.
 \end{subchunk}

\begin{sublemma}
\label{finite koszul}
The complex of $S$-modules $\koszul{\bsx}N$ is homologically finite,
and
 \[
\Supp_S\HH\big(\koszul{\bsx}N\big)=
\Supp_S(S/\bsx S)\cap\Supp_S\HH(N)
 \]
 \end{sublemma}

\begin{proof}
We may assume $\bsx=\{x\}$.  In view of the preceding observation,
finiteness follows from \eqref{koszul:les}.  By the same token, one gets
the first equality below
\begin{align*}
\Supp_S\HH\big(\koszul{x}N\big)
&=\Supp_S\big(\HH(N)/x\HH(N)\big)\cup
\Supp_S\big(\Ker\lambda^{\HH(N)}_x\big)\\
&=\Supp_S(S/\bsx S)\cap\Supp_S\HH(N)
\end{align*}
To get the second equality, remark that $\Supp_S\HH(N)\cap
\Supp_S(S/\bsx S)$ is equal to the support of $\HH(N)/x\HH(N)$
and contains that of $\Ker\lambda^{\HH(N)}_x$.
 \end{proof}

\begin{sublemma}
\label{finite betti}
For each $n\in\BZ$, both $\Tor^R_n(M,\koszul{\bsx}N)$ and
$\Ext_R^n(M,\koszul{\bsx}N)$ are finite $S$-modules annihilated
by $\bsx S+\hann RMS+\hann SN$
\end{sublemma}

\begin{proof}
The complex $\koszul{\bsx}N)$ is homologically finite by \eqref{finite
koszul}, so the assertion about finiteness follows from Lemma
\eqref{newfinite}.  The assertion about annihilation comes from
\eqref{hann koszul} and \eqref{hann2}.
 \end{proof}

\section{Dimensions}
\label{Dimensions}

We extend some well known theorems on modules of finite injective
dimension, namely, the Bass Equality, see \cite[(18.9)]{Ma}, and the
``Bass Conjecture'' proved by P.\ Roberts, see \cite[\S 3.1]{Rb}
and \cite[(13.4)]{Rb1}.  The novelty is that finiteness over $R$ is
relaxed to finiteness over $S$.  Khatami and Yassemi \cite[(3.5)]{KY},
and Takahashi and Yoshino \cite[(5.2)]{TY} have independently obtained
the first equality below.

\begin{theorem}
\label{bassroberts}
If $L$ is a finite $S$-module and $\id_RL<\infty$, then
 \[
\id_RL=\depth R =\dim R
 \]
\end{theorem}

This is proved at the end of the section, after we engineer a
situation where the original results apply.  We recall how
some classical concepts extend to complexes.

Let $M$ be a homologically bounded complex of $R$-modules.  Its \emph{flat
dimension} and its \emph{injective dimension} over $R$ are, respectively, the
numbers
\begin{gather*}
\fd_RM=\sup\left\{n\in\BZ\left\vert
\begin{gathered}Y_n\ne0 \text{ for some bounded complex}
\\
\text{of flat $R$-modules $Y$ with } Y\simeq M
\end{gathered}
\right\}\right.
\quad\text{and}
\\
\id_RM=-\inf\left\{n\in\BZ\left\vert
\begin{gathered}Y_n\ne0 \text{ for some bounded complex}
\\
\text{of injective $R$-modules $Y$ with } Y\simeq M
\end{gathered}
\right\}\right.
\end{gather*}
The $n$th \emph{Betti
number} $\betti nRM$ and the $n$th \emph{Bass number} $\bass nRM$
are defined to be
\begin{gather*}
\betti nRM = \rank_k \Tor_n^{R}(k,M)
\qquad\text{and}\qquad
  \bass nRM = \rank_k \Ext^n_{R}(k,M)
\end{gather*}

\begin{chunk}
\label{hdim:classical}
By \cite[(5.5)]{AF:huc}, for the homologically finite complex $N$, one has
\begin{gather*}
\fd_RN=\sup\{n\in\BZ\mid\betti nRN\ne0\}\\
\id_RN=\sup\{n\in\BZ\mid\bass nRN\ne0\}
\end{gather*}
 \end{chunk}

Next we present a construction of Avramov, Foxby, and B.~Herzog
\cite{AFH}, which is used in proofs of many theorems throughout the paper.

\begin{chunk}{\bf Cohen factorizations. }
\label{Cohen factorizations}
Let $\grave\vf\colon R\to \wh S$ denote the composition of $\vf$ with the
completion map $S\to\wh S$.  A \emph{Cohen factorization} of $\grave\vf$
is a commutative diagram
\begin{equation*}
\xymatrixrowsep{3pc}
\xymatrixcolsep{3pc}
\xymatrix{
   &R' \ar@{->}[dr]^{\textstyle\vf'} \\ 
R \ar@{->}[rr]^{\textstyle\grave\vf}
\ar@{->}[ur]^{\textstyle\dot\vf} && \wh S }
\end{equation*}
of local homomorphisms such that the map $\dot\vf$ is flat, the ring
$R'$ is complete, the ring $R'/\fm R'$ is regular, and the map $\vf'$ is
surjective.

Clearly, the homomorphisms and rings in the diagram above satisfy
 \[
\edim\vf\leq\edim\dot\vf=\dim(R'/\fm R')
 \]
When equality holds the factorization is said to be \emph{minimal}.
It is proved in \cite[(1.5)]{AFH} that the homomorphism $\grave\vf$
always has a minimal Cohen factorization.
\end{chunk}

For the rest of this subsection we fix a minimal Cohen factorization
of $\grave\vf$.  With the next result we lay the groundwork for several
proofs in this paper.

\begin{proposition}
\label{series:regular}
For each minimal set $\bsx$ of generators of $\fn$ modulo $\fm S$ there
are isomorphisms in the derived category of $S$-modules, as follows:
\begin{gather*}
\dtensor kR{\big(\koszul{\bsx}{N}\big)}
\simeq\koszul{\bsx}{\dtensor kR{N}}\simeq
\dtensor l{R'}{\wh N}
\\
\RHom Rk{\koszul{\bsx}N}
\simeq \koszul\bsx{\Rhom Rk{N}}
\simeq \shift^{\edim\vf}\,\Rhom {R'}l{\wh N}
\end{gather*}
 \end{proposition}

\begin{proof}
We give the proof for homomorphisms and omit the other, similar,
argument.

Considering $S$ as a subset of $\wh S$, lift $\bsx$ to a subset $\bsx'$
of $R'$, containing $\edim\vf$ elements.  As the ring $P=R'/\fm R'$
is regular, the image of $\bsx'$ in $P$ is a regular system of
parameters. Thus, the Koszul complex $\koszul{\bsx'}P$ is a free
resolution of $l$, that is, $\koszul{\bsx'}P\simeq l$.  This accounts
for the first isomorphism in the chain
\begin{align*}
\Rhom{R'}l{\wh N}
&\simeq\RHom{R'}{\koszul{\bsx'}P}{\wh N}\\
&\simeq\RHom{R'}{P}{\Rhom{R'}{\koszul{\bsx'}{R'}}{\wh N}}\\
&\simeq\RHom{R'}{{R'}\otimes_Rk}{\Hom_{R'}({\koszul{\bsx'}{R'},{\wh N}})}\\
&\simeq\RHom{R}k{\Hom_{R'}({\koszul{\bsx'}{R'},{\wh N}})}\\
&\simeq\RHom{R}k{\Hom_{\wh S}({\koszul{\bsx}{\wh S},{\wh N}})}\\
&\simeq\RHom{R}k{\Hom_S(\koszul{\bsx}S,{\wh N})}
 \end{align*}
The third isomorphism holds because $\koszul{\bsx'}{R'}$ is a finite
complex of free $R'$-modules, while the remaining ones are adjunctions.

The first isomorphism below holds because $\koszul{\bsx}S$ is a finite
free complex:
\begin{align*}
\RHom{R}k{\Hom_S(\koszul{\bsx}S,{\wh N})}
&\simeq\RHom{R}k{\Hom_S(\koszul{\bsx}S,N)\otimes_S\wh S}\\
&\simeq\RHom{R}k{\Hom_S(\koszul{\bsx}S,N)}\otimes_S\wh S
 \end{align*}
The second one does because $\wh S$ is a flat $S$-module,
$\Hom_S(\koszul{\bsx}S,N)$ has bounded homology, and $k$ has a resolution
by finite free $R$-modules.

As $\wh S$ is flat over $S$, and each
$\HH_n\RHom{R}k{\Hom_S(\koszul{\bsx}S,N)}$ has finite length over $S$
by \eqref{finite betti}, the canonical morphism
 \[
\RHom{R}k{\Hom_S(\koszul{\bsx}S,N)}\otimes_S\wh S
\longleftarrow\RHom{R}k{\Hom_S(\koszul{\bsx}S,N)}
 \]
is an isomorphism in $\dcat(S)$.  Furthermore, we have isomorphisms
\begin{align*}
\RHom{R}k{\Hom_S(\koszul{\bsx}S,N)}
&\simeq\Rhom{R}k{\shift^{-\edim\vf}{\koszul{\bsx}N}}\\
&\simeq\shift^{-\edim\vf}\Rhom{R}k{\koszul{\bsx}N}\\
&\simeq\shift^{-\edim\vf}\koszul\bsx{\Rhom Rk{N}}
 \end{align*}
obtained from \eqref{koszul dual}, the definition of $\shift$, and
\eqref{free koszul}, respectively.  Linking the chains above
we obtain the desired isomorphisms.
 \end{proof}

Our first application of the proposition above is to homological
dimensions.  The assertions on the flat dimension of $\wh N$ over
$R'$ are known, see \cite[(3.2)]{AFH}.

\begin{corollary}
\label{dim properties}
The following (in)equalities hold:
\begin{gather*}
\fd_{R}N \leq \fd_{R'}\wh N \leq \fd_{R}N +\edim\vf\\
\id_{R'}\wh N=\id_RN+\edim\vf
\end{gather*}
 \end{corollary}

\begin{proof}
Let $\bsx$ be a minimal generating set for $\fn$ modulo $\fm S$.
{}From the proposition and the expressions for homological dimensions in
\eqref{hdim:classical}, one gets
\begin{gather*}
\fd_{R'}\wh N=\sup\HH\big(\koszul\bsx{\dtensor kRN}\big)
\\
\id_{R'}\wh N=-\inf\HH\big(\koszul\bsx{\Rhom RkN}\big)+\edim\vf
\end{gather*}
As $\edim\vf=\card(\bsx)$, and the $S$-modules $\HH_n(\dtensor kRN)$
and $\Rhom RkN$ are finite for each $n$ by \eqref{newfinite},
Lemma \eqref{koszul:sup} gives the desired (in)equalities.
 \end{proof}

The interval for the flat dimension of $\wh N$ over $R'$ cannot
be narrowed, in general, as the following special case demonstrates.

\begin{example}
If $R$ is a field, the ring $S$ is complete and regular, and $N$ is an
$S$-module, then the statement above reduces to inequalities
$0\leq \fd_{S}N \leq\dim S$.
 \end{example}

\begin{proof}[Proof of Theorem {\em\eqref{bassroberts}}]
Let $R\to R'\to \wh S$ be a minimal Cohen factorization of
$\grave\vf$.  Corollary \eqref{dim properties} and the
minimality of the factorization yield
\begin{align*}
\id_R L  &= \id_{R'}{\wh L} - \edim\vf\\
&=\id_{R'}{\wh L} -\edim(R'/\fm R')
\end{align*}
In particular, $\id_{R'}{\wh L}$ is finite, so Bass' Formula gives the
first equality below. The other two are due to the
regularity of $R'/\fm R'$, and the flatness
of $R'$ over $R$:
\begin{align*}
\id_{R'}{\wh L} - \edim(R'/\fm R')&= \depth R'-\edim(R'/\fm R')\\
& = \depth R'-\depth(R'/\fm R')\\
&=\depth R
\end{align*}
Roberts' Theorem shows that $R'$ is Cohen-Macaulay; by flat descent,
so is $R$.
 \end{proof}

\section{Koszul invariants}
\label{Koszul invariants}

Koszul complexes on minimal sets of generators of the maximal ideal $\fn$
of $S$ appear systematically in our study.  This section is devoted to
establishing relevant properties of such complexes.  While several of them
are known, in one form or another, we have not been able to find in the
literature statements presenting the detail and generality needed below.
It should be noted that the invariants discussed in this section depend
only on the action of $S$ on $N$, while $R$ plays no role.

\begin{lemma}
\label{koszul:unique}
Set $s=\edim S$, let $\bsw$ be a minimal set of generators of $\fn$,
and let $\bsz$ be a set of generators of $\fn$ of cardinality $u$.

Let $A$ be the Koszul complex of a set consisting of $(u-s)$ zeros.
\begin{enumerate}[{\quad\rm(1)}]
\item There is an isomorphism of DG algebras
\(
\koszul{\bsz} S \cong \koszul{\bsw} S\otimes_S A
\)
\item
If $\bsz$ minimally generates $\fn$, then $\koszul{\bsz}S\cong \koszul{\bsw}S$.
\item The complexes of $S$-modules $\koszul{\bsz} N$ and
$\koszul{\bsw} N\otimes_S A$ are isomorphic.
\end{enumerate}
\end{lemma}

\begin{proof}
(1) Since $\bsz \subseteq (\bsw)$ we may write
  \[
z_j=\sum_{i=1}^s a_{ij} w_i \quad\text{with}\quad a_{ij}\in S
\quad\text{for}\quad j=1,\dots,u
  \]
Renumbering the elements of $\bsz$ if necessary, we may assume that
the matrix $(a_{ij})_{1\les i,j\les s}$ is invertible.  Choose bases
$\{e_1,\dots,e_s\}$ and $\{f_1,\dots,f_u\}$ of the degree $1$ components
of $\koszul{\bsw}S$ and $\koszul{\bsz}S$, respectively, such that
  \[
\dd(e_i)= w_i \text{ for } 1\le i\le s \quad\text{and}\quad \dd(f_j)=
z_j \text{ for }1\le j\le u
  \]
Let $\{e_{s+1},\dots,e_u\}$ be a basis for the free module $A_1$.
The graded algebra underlying $\koszul\bsz S$ is the exterior algebra on
$\{f_1,\dots,f_u\}$, so there exists a unique morphism $\varkappa\colon
\koszul{\bsz} S\to\koszul\bsw S\otimes_S A$ of graded algebras over $S$,
such that
  \[
\varkappa(f_j)= \begin{cases}
  \sum_{i=1}^s a_{ij}(e_i\otimes1)    &\text{for}\quad 1\le j\le s\\
  \sum_{i=1}^s a_{ij}(e_i\otimes1)-(1\otimes e_j)&\text{for}\quad s<
  j \le u
\end{cases}
  \]
Note that $\varkappa$ is an isomorphism in degree $1$. The graded
algebra underlying $\koszul\bsw S\otimes_S A$ is the exterior algebra
on $\{e_1\otimes1,\dots,e_s\otimes1, 1\otimes e_{s+1},\dots,1\otimes
e_u\}$, so $\varkappa$ is an isomorphism of graded algebras.  The formulas
above yield $\dd\varkappa(f_j)=\varkappa\dd(f_j)$ for $j=1,\dots,u$.  It
follows that $\varkappa$ is a morphism of DG algebras.

The assertions in (2) and (3) are consequences of (1).
 \end{proof}

\begin{chunk}
\label{koszul poly}
We let $\rkoszul SN$ denote the complex $\koszul{\bsw}N$ on a minimal
generating set of $\fn$. There is little ambiguity in the notation
because different choices of $\bsw$ yield isomorphic complexes; see
Lemma (\ref{koszul:unique}.2). We set $K^S = \rkoszul SS$.

Lemma \eqref{finite betti} shows that $\HH_n(\rkoszul SN)$ is
a finite $l$-vector space for each $n\in\BZ$, which is trivial
for all $|n|\gg0$.  Thus, one can form a Laurent polynomial
 \[
K^S_N(t)=\sum_{n\in\BZ}\kappa^S_n(N)\,t^n
\qquad\text{where}\qquad
\kappa^S_n(N)=\rank_l\HH_n(\rkoszul SN)
 \]
with non-negative coefficients.  We call it the \emph{Koszul
polynomial} of $N$ over $S$.
 \end{chunk}

Some numbers canonically attached to $N$ play a role in further
considerations.

\begin{chunk}
\label{depth}
The \emph{depth} and the \emph{type} of $N$ over $S$ are defined,
respectively, to be the number
\begin{gather*}
\depth_SN=\inf\{n\in\BZ\mid\Ext^n_S(l,N)\ne0\}\\
\type_SN=\rank_l\Ext^{\depth_SN}_S(l,N)
\end{gather*}
When $N$ is an $S$-module the expressions above yield the
familiar invariants.  By \cite[(6.5)]{Iy}, depth can be computed
from a Koszul complex, namely
 \[
\depth_SN=\edim S-\sup\HH\big(\rkoszul SN\big)
 \]
\end{chunk}

\begin{chunk}
\label{cohen presentation}
A \emph{Cohen presentation} of $\wh S$ is an isomorphism
$\wh S\cong T/\fb$ where $(T,\fr,l)$ is a regular local ring.
Cohen's Structure Theorem proves one always exists.  It can be taken
to be \emph{minimal}, in the sense that $\edim T=\edim S$: If not,
pick $x\in\fb\smallsetminus \fr^2$, note the isomorphism $\wh S\cong
(T/xT)/(\fb/xT)$ where $T/xT$ is regular, and iterate.
\end{chunk}

In the next result, and on several occasions in the sequel, we use the
theory of dualizing complexes, for which we refer to Hartshorne
\cite[V.2]{Ha}.

\begin{chunk}
\label{dualizing complex}
Let $Q$ be a local ring with residue field $h$.  A \emph{normalized
dualizing complex} $D$ for $Q$ is a homologically finite complex of
$Q$-modules $E$, such that
 \[
\Rhom{Q}{h}E\simeq h
 \]
These properties describe $E$ uniquely up to isomorphism in $\dcat(Q)$.
When $E$ is a normalized dualizing complex, for every complex of
$Q$-modules $X$ we set
 \[
X^\dagger=\RHom{Q}XE
 \]
If $X$ is homologically finite, then so is $X^\dagger$, and the
canonical morphism $X\to X^{\dagger\dagger}$ is an isomorphism
in $\dcat(Q)$.  When necessary, we write $X^\dagger_Q$ instead of
$X^\dagger$.

A local ring $P$ is Gorenstein if and only if $\shift^{\dim P}P$ is a
normalized dualizing complex for $P$. If this is the case and $P\to Q$
is a surjective homomorphism, then $\RHom{P}{Q}{\shift^{\dim P}P}$ is a
normalized dualizing complex for $Q$.  Thus, Cohen's Structure Theorem
implies that every complete local ring has a dualizing complex.
 \end{chunk}

\begin{chunk}
\label{dualizing koszul}
The complexes of $S$-modules $\rkoszul S{N^\dagger}{}^\dagger$ and
$\shift^{-\edim S}\rkoszul SN$ are isomorphic.

Indeed, if $D$ is a normalized dualizing complex for $S$, then one has
\begin{align*}
\Hom_S(\rkoszul S{N^\dagger},D)
&\cong\Hom_S(K^S\otimes_S{N^\dagger},D)\\
&\cong\Hom_S(K^S,N^{\dagger\dagger})\\
&\cong\Hom_S(K^S,N)\\
&\cong\shift^{-\edim S}\rkoszul SN
\end{align*}
 \end{chunk}

We set $\codim S=\edim S-\dim S$.  If $L$ is an $S$-module, then $\nu_SL$
denotes the minimal number of generators of $L$.  The following result
collects most of the arithmetic information required on the Koszul
polynomial.

\begin{theorem}
\label{koszul:poly}
Assume $\HH(N)\ne0$, and set
 \[
s=\edim S\,,\quad c=\codim S\,,\quad i=\inf\HH(N)\,,\quad
m=\nu_S\HH_i(N)\,,\quad g=\depth_SN\,.
 \]
The Koszul polynomial $K^S_N(t)$ has the following properties.
\begin{enumerate}[\quad\rm(1)]
\item
$\ord K^S_N(t)=i$ and $\kappa^S_{i}(N)=m$.
\item
$\deg K^S_N(t)=(s-g)$ and $\kappa^S_{s-g}(N)=\type_SN$.
\item
$K^S_N(t)=K^{\wh S}_{\wh N}(t)=\rP T{\wh N}$, where $\wh S\cong T/\fb$
is a minimal Cohen presentation.
\item
$K^{\wh S}_{{\wh N}^\dagger}(t)=t^{s}K^S_N(t^{-1})$.
\end{enumerate}
If $S$ is not regular, then the following also hold.
\begin{enumerate}[\quad\rm(1)]
\item[\rm(5)]
$\kappa^S_{i+1}(N)\ge \kappa^S_{i}(N)+c-1$.
\item[\rm(6)]
$\kappa^S_{s-g-1}(N)\ge \kappa^S_{s-g}(N)+c-1$.
\item[\rm(7)]
$K^S_N(t)=(1+t)\cdot L(t)$ where $L(t)$ is a Laurent polynomial of the
form
\[
L(t)=mt^i+\text{higher order terms with non-negative coefficients}
\]
\end{enumerate}
 \end{theorem}

 \begin{proof}
(1)  The equality $\ord K^S_N(t)=i$ comes from \eqref{koszul:sup}.

The equality $\kappa^S_{i}(N)=m$ results from the canonical
isomorphisms
 \[
\HH_i(\rkoszul SN)\cong\HH_i(K^S\otimes_SN)\cong
\HH_0(K^S)\otimes_S\HH_i(N)=l\otimes_S\HH_i(N)
 \]

(2) The equality $\deg K^S_N(t)=(s-g)$ comes from \eqref{depth}.

Because $\fn$ annihilates each $\HH_q(K^S)$, the standard spectral sequence
with
 \[
E_{pq}^2=\Ext^{-p}_S(\HH_q(K^S),N)\implies \HH_{p+q}\Hom_S(K^S,N)
  \]
has $E_{pq}^2=0$ for $p>-g$.  It has a corner, which yields
the first isomorphism below
 \[
\Ext^{g}_S(l,N)=E_{-g\,0}^2\cong\HH_{-g}\Hom_S(K^S,N)
\cong\HH_{s-g}(\rkoszul SN)
 \]
The second one is from \eqref{koszul dual}.  Thus, we get
$\kappa^S_{s-g}(N)=\type_SN$.

(3) Any minimal generating set $\bsw$ of $\fn$ minimally generates the
maximal ideal of $\wh S$, and the canonical map $\koszul {\bsw}{N}\to
\koszul{\bsw}{\wh N}$ is an isomorphism in $\dcat(S)$.

Having proved the first equality, for the rest of the proof we may
assume that $S$ is complete.  We choose an isomorphism $S\cong T/\fb$,
where $(T,\fr,l)$ is a regular local ring with $\dim T=\edim S$, and a
normalized dualizing $D$ complex for $S$.

For the second equality in (3), note that a
minimal generating set of $\fr$ maps to a minimal generating set of
$\fn$. As $T$ is regular, $K^T$ is a free resolution of $l$ over $T$.
Now invoke the isomorphism $K^T\otimes_TN\cong\rkoszul SN$.

(4) By (3), we may assume $S$ is complete. There is a spectral sequence with
 \[
E_{pq}^2=\Ext^{-p}_S(\HH_q(\rkoszul S{N^\dagger}),D)
\implies \HH_{p+q}\Hom_S(\rkoszul S{N^\dagger},D)
  \]
As $\HH(\rkoszul S{N^\dagger})$ is a finite dimensional
graded $l$-vector space, the defining property of $D$ implies
$E_{pq}^2=0$ for $p\ne0$, and yields for all $n$ isomorphisms
 \[
\Hom_l(\HH_{n}(\rkoszul S{N^\dagger}),l)\cong
\HH_{n}\Hom_S(\rkoszul S{N^\dagger},D)
 \]
of $l$-vector spaces.  On the other hand, \eqref{dualizing koszul} 
yields for all $n$ isomorphisms
 \[
\HH_{n}\Hom_S(\rkoszul S{N^\dagger},D)\cong
\HH_{s-n}(\rkoszul SN)
 \]
of $l$-vector spaces.  Comparison of ranks produces the desired
equalities.

As $N$ is homologically finite, it is isomorphic in $\dcat(T)$
to a finite free complex of $T$-modules $G$ that satisfies
$\dd(G)\subseteq\fr G$; see \cite[(II.2.4)]{Rb}.  {}From (3) we get
 \[
\rank_T(G_n)=\kappa^S_n(N)
\qquad\text{for all}\qquad
n\in\BZ
 \]
We use the complex $G$ in the proofs of the remaining statements.

(5) As $\HH(G)\cong\HH(N)$ and $\inf\HH(N)=i$, there is an exact sequence
of $T$-modules.
 \[
G_{i+1}\lra G_i\lra \HH_i(N) \lra 0
 \]
The Generalized Principal Ideal Theorem gives the first inequality below:
\begin{align*}
\kappa^S_{i+1}(N)-\kappa^S_{i}(N)+1&\ge\height\Ann_T(\HH_i(N))\\
&=\dim T-\dim_T\HH_i(N)\\
&\ge\dim T-\dim S\\
&=\edim S-\dim S
\end{align*}
see \cite[(13.10)]{Ma}.  The first equality holds because $T$ is regular,
and hence catenary and equidimensional, while the other relations
are clear.

(6) In view of (4), the inequality (5) for the complex $N^\dagger$
yields the desired inequality, $\kappa^S_{s-g-1}(N)\ge \kappa^S_{s-g}(N)+c-1$.

(7)  Let $U$ be the field of fractions of $T$. As $S$ is not
regular, $\fb$ is not zero, so $N\otimes_TU=0$. Also,
$G\otimes_TU\simeq N\otimes_TU$ holds in $\dcat(U)$, so
$G\otimes_TU$ is an exact complex of finite $U$-vector spaces, and
hence is the mapping cone of the identity map of a complex $W$
with trivial differential. For $L(t)=
\sum_{n\in\BZ}\rank_U(W_n)\,t^n$ one gets
\begin{align*}
K^S_N(t) =\rP TN  & = \sum_{n\in\BZ}\rank_T(G_n)\,t^n \\
         & = \sum_{n\in\BZ}\rank_U(G\otimes_TU)_n\,t^n \\
         & = \sum_{n\in\BZ}\rank_U(W_n)\,t^n \cdot (1+t)
\end{align*}
It follows that $L(t)$ has order $i$ and starts with $mt^i$.
 \end{proof}

An important invariant of Koszul complexes is implicit in a result of Serre.

\begin{chunk}
\label{serre number}
Let $L$ be a finite $S$-module.  Set $K=\koszul \bsw L$, where
$\bsw=\{w_1,\dots,w_s\}$ mini\-ma\-lly generates $\fn$, and $\fn^j=S$
for $j\le0$.  As $\dd(\fn^jK_n)\subseteq\fn^{j+1}K_{n-1}$ for all
$j$ and $n$, for each $i\in\BZ$ one has a complex of $S$-modules
  \[
I_L^i =\quad 0\lra \fn^{i-s} K_s \lra \fn^{i-s+1}K_{s-1} \lra \cdots
\lra \fn^{i-1}K_1 \lra \fn^{i}K_0\to 0
  \]
Lemma (\ref{koszul:unique}.2) shows that it does not depend on the choice
of $\bsw$, up to isomorphism.  We define the \emph{spread} of $N$ over
$S$ to be the number
 \[
\serre_SL=\inf\{i\in\BZ\mid \HH(I_L^j)=0\text{ for all }j\ge i\}
 \]
Serre \cite[(IV.A.3)]{Se} proves that it is finite. We write $\serre S$
in place of $\serre_SS$.
 \end{chunk}

In some cases the new invariant is easy to determine.

\begin{example}
\label{useless}
If $\fn^vL=0\ne\fn^{v-1}L$, then $\serre_SL=\edim S+v$.

Indeed, for $s=\edim S$ one has $I^{s+v}_L=0$ and
$\HH_{s}(I^{s+v-1}_L)\supseteq\fn^{v-1}K_s\ne 0$.
 \end{example}

To deal with more involved situations we interpret $\serre_SN$ in terms
of associated graded modules and their Koszul complexes.

\begin{chunk}
Given an $S$-module $L$ we set $\fn^i L=L$ for all $i\le0$ and
write $\grass\fn L$ for the graded abelian group associated
to the $\fn$-adic filtration $\{\fn^i L\}_{i\in\BZ}$; it is
a graded module over the graded ring $\grass\fn S$.  Note that
$\bsw^*=\{w_i+\fn^2,\dots,w_s+\fn^2\}\subseteq\grass\fn S\mathstrut_1$
minimally generates $\grass\fn S$ as an algebra over $\grass\fn
S\mathstrut_0=l$.

The Koszul complex $\koszul{\bsw^*}{\grass\fn S}$ becomes a complex
of graded $\grass\fn S$-modules by assigning bidegree $(1,1)$ to each
generator of $\koszul{\bsw^*}{\grass\fn S}_1$.  If $\bsw'$ also is a
minimal generating set for $\fn$, then the complexes of $S$-modules
$\koszul{\bsw}S$ and $\koszul{\bsw'}S$ are isomorphic, see Lemma
(\ref{koszul:unique}.2); any such isomorphism induces an isomorphism
$\koszul{\bsw'{}^*}{\grass\fn S}\cong\koszul{\bsw^*}{\grass\fn S}$ of
complexes of graded $\grass\fn S$-modules.  We let $K^{\grass\fn S}$
denote any such complex, and form the complex of graded $\grass\fn
S$-modules
 \[
\rkoszul{\grass\fn S}{\grass\fn L}= K^{\grass\fn S}\otimes_{\grass\fn
S}{\grass\fn L}
 \]
  \end{chunk}

In the next proposition we refine Serre's result \eqref{serre number}.
Rather than using the result itself, we reinterpret some ideas of its
proof in our argument.

\begin{proposition}
\label{serre}
The following equality holds:
\[
\serre_SL=
\sup\{i\in\BZ\mid \HH_n\big(\rkoszul{\grass\fn S}{\grass\fn L}\big)_i \ne 0
               \,\text{for some $n\in\BZ$}\}
\]
\end{proposition}

\begin{proof}
Set $I=I_L$.  By construction, we have isomorphisms
\[
\bigoplus_{i\in\BZ}I^i/I^{i+1}\cong\koszul{\bsw^*}{\grass\fn L}
\]
of complexes of graded $\ell$-vector spaces.  Thus, to prove the
proposition it suffices to show that the following claim holds for
each integer $a$.

\begin{Claim}
$\HH(I^i)=0$ for $i\ge a$ if and only if $\HH(I^i/I^{i+1})=0$ for
$i\ge a$.
\end{Claim}

Indeed, the exact sequence of complexes \( 0\to I^{i+1} \to I^i \to
I^i/I^{i+1}\to 0 \) shows that if $\HH(I^i)$ vanishes for $i\ge a$, then
$\HH(I^i/I^{i+1})=0$ for $i\ge a$.  Conversely, assume that the latter
condition holds.  Induction on $j$ using the exact sequences
 \[
0\to I^{a+j}/I^{a+j+1} \to I^{a}/I^{a+j+1} \to I^{a}/I^{a+j}\to 0
 \]
of complexes yields $\HH(I^a/I^{a+j})=0$.  Since the inverse system of complexes
 \[
\cdots\to I^{a}/I^{a+j+1}\to I^{a}/I^{a+j}\to\cdots\to I^{a}/I^{a+2}
\to I^{a}/I^{a+1}\to 0
 \]
is surjective, one gets $\HH\big(\varprojlim_{j}\, {I^{a}/I^{a+j}}\big) = 0$.
Now $(I^{a+j})_n = \fn^{a+j-n}K_n$ for each $n$, so
 \(
\varprojlim_{j}({I^{a}/I^{a+j}})_n
 \)
is the $\fn$-adic completion of the $S$-module ${\fn^{a-n}K_n}$.
This module is finite, so its completion is isomorphic to
$(\fn^{a-n}K_n)\otimes_S\wh S$. The upshot is that
 \[
\textstyle{\varprojlim_{j}}\big(\, {I^{a}/I^{a+j}}\big)
\cong I^a \otimes_S\wh S
 \]
as complexes of $S$-modules.  {}From here and
the flatness of $\wh S$ over $S$ one gets
 \[
\HH(I^a)\otimes_S\wh S\cong\HH(I^a\otimes_S\wh S)=
\HH\big(\textstyle{\varprojlim_{j}}\,{I^{a}/I^{a+j}}\big)=0
 \]
Since the $S$-module $\wh S$ is also faithful, we deduce $\HH(I^a)=0$,
as claimed.
 \end{proof}

There is another interpretation of the numbers in the preceding proposition.

\begin{remark}
Let $B=\{B_j\}_{i\ges 0}$ be a graded ring with $B_0\cong l$, finitely
generated as a $B_0$-algebra, and let $M$ be a graded and finite
$B$-module.  The number
 \[
\astar_B(M) = \sup\{i\in\BZ \mid \lch n{\fb}M_i \ne 0\, \text{for some
$n\in\BZ$} \}
 \]
where $\fb$ is the irrelevant maximal ideal $B_{\ges 1}$ of $B$, and
$\lch n{\fb}M$ is the $n$th local cohomology module of $M$, has been
studied by Trung \cite{Tr} and others.  We claim: if $b$ is the minimal
number of homogeneous generators of the $l$-algebra $B$, then
 \[
\sup\{i\in\BZ \mid \HH_n(\rkoszul BM)_i \ne 0\, \text{for some $n\in\BZ$}
\} =\astar_B(M) + b
 \]

This preceding identity can be verified by using a spectral sequence with
 \[
E^{pq}_2=\HH_p\big(\rkoszul B{\lch{-q}{\fb} M}\big) \implies
\HH_{p+q}(\rkoszul BM)
  \]
that lies in the fourth quadrant and has differentials $d^{pq}_r\colon
E^{pq}_r\to E^{p-r,q+r-1}_r$.
 \end{remark}

\section{Sequences}
\label{Sequences}

The flat dimension and the injective dimension of $N$ over $R$ can be
determined from the \emph{number} of non-vanishing groups $\Tor_n^R(k,N)$
and $\Ext^n_R(k,N)$, as recalled in \eqref{hdim:classical}.  When this
number is infinite, we propose to use the \emph{sizes} of these
groups as a measure of the homological intricacy of the complex $N$
over $R$.  Each such group carries a canonical structure of \emph{finite
$S$-module}, so its size is reflected in a number of natural invariants,
such as minimal number of generators, multiplicity, rank, or length.
Among these it is the length over $S$, denoted $\ell_S$, that has the best
formal properties, but it is of little use unless some extra hypothesis (for
instance, that the ring $S/\fm S$ is artinian) guarantees its finiteness.

To overcome this problem, we take a cue from Serre's approach to multiplicities and
replace $N$ by an appropriate Koszul complex.  As there is no canonical choice of such a
complex, the question arises whether the new invariants are well defined.  In this section
we prove that indeed they are, and describe them in alternative terms.

\subsection{Betti numbers and Bass numbers}
\label{bettibass}
Let $\bsx$ be a minimal generating set for $\fn$ modulo $\fm S$.
Lemma \eqref{finite betti} shows that $\Tor^R_n(k,\koszul{\bsx}N)$
and $\Ext_R^{n}(k,\koszul{\bsx}N)$ are finite $S$-modules annihilated
by $\fn$, hence they are $l$-vector spaces of finite rank.

For each $n\in\BZ$ we define the $n$th \emph{Betti number $\betti n\vf N$
of $N$ over $\vf$} and the $n$th \emph{Bass number $\bass n\vf N$ of $N$
over $\vf$} by the formulas
\begin{gather*}
  \betti n\vf N = \rank_l\Tor^R_n(k,\koszul{\bsx}N)  \\
  \bass n\vf N = \rank_l\Ext_R^{n-\edim\vf}(k,\koszul{\bsx}N)
\end{gather*}
These numbers are invariants of $N$:  pick a minimal Cohen
factorization $R\to R'\to\wh S$ of $\grave\vf$ and apply Proposition
\eqref{series:regular} to get isomorphisms of $S$-modules
 \[
\tag{4.1.0}
\begin{gathered}
\Tor^R_n(k,\koszul{\bsx}N)\cong\Tor^{R'}_n(\ell,\wh N)\\
\Ext^{n-\edim\vf}_R(k,\koszul{\bsx}N)\cong\Ext_{R'}^n(\ell,\wh N)
\end{gathered}
 \]
The vector spaces on the right do not depend on $\bsx$, so neither do
their ranks.  Such an independence may not be taken for granted, because
$\koszul{\bsx}N$ is not defined uniquely up to isomorphism in $\dcat(S)$;
even its homology may depend on $\bsx$.

\begin{subexample}
Let $k$ be a field, set $R=k[[u]]$ and $S=k[[x,y]]/(xy)$, and let
$\vf\colon R\to S$ be the homomorphism of complete $k$-algebras defined
by $\vf(u)=y$.  The Koszul complexes on $\bsx=\{x\}$ and $\bsx'=\{x+y\}$
satisfy
 \[
\HH_n(\koszul{\bsx}S)\cong
\begin{cases}
S/(x)&\text{for }n=0,1\\
0 &\text{otherwise}
\end{cases}
\qquad
\HH_n(\koszul{\bsx'}S)\cong
\begin{cases}
S/(x+y)&\text{for }n=0\\
0 &\text{otherwise}
\end{cases}
 \]
 \end{subexample}

Next we show that Betti or Bass numbers over $\vf$ are occasionally
equal to the eponymous numbers over $R$, but may differ even when $N$
is a finite $R$-module.

\begin{subremark}
\label{surjective}
If the map $k\to S/\fm S$ induced by $\vf$ is bijective (for instance,
if $\vf$ is the inclusion of $R$ into its $\fm$-adic completion $\wh
R$, or if $\vf$ is surjective), then
 \[
\betti n\vf N = \betti nRN \quad\text{and}\quad \bass n\vf N = \bass nRN
 \]
Indeed, in this case the set $\bsx$ is empty, and ranks over $l$ are
equal to ranks over $k$.
  \end{subremark}

\begin{subexample}
\label{rod ex1}
Set $S=R[[x]]$ and let $\vf\colon R\to S$ be the canonical inclusion.
The $S$-module $N=S/xS$ is finite over $R$, and elementary computations yield
\begin{gather*}
\betti n{\vf}N=
\begin{cases}
\betti nRN +1& \text{ for }n=1 \\
\betti nRN &  \text{ for } n\ne1
\end{cases}
\\
\bass n{\vf}N = \bass{n-1}{R}R+\bass n{R}R
\quad\text{for each}\quad n\in\BZ
\end{gather*}
 \end{subexample}

In general, it is possible to express invariants over $\vf$
in terms of corresponding invariants over a ring, but a new ring
is necessary:

\begin{subremark}
\label{series:cohen}
If $R\xra{\dot\vf} R'\xra{\vf'}\wh S$ is a minimal Cohen
factorization of $\grave\vf$, then for each $n\in\BZ$ the following
equalities hold:
\begin{gather*}
\betti n{\vf}{N}=\betti n{\grave\vf}{\wh N}
=\betti n{\wh\vf}{\wh N}=\betti n{R'}{\wh N}
\\
\bass n{\vf}{N}=\bass n{\grave\vf}{\wh N}
=\bass n{\wh\vf}{\wh N}=\bass n{R'}{\wh N}
\end{gather*}

Indeed, the isomorphisms of $S$-modules in (4.1.0) yield $\betti
n{\vf}{N}=\betti n{R'}{\wh N}$.  The same argument applied to $\grave\vf$
yields $\betti n{\grave\vf}{\wh N}=\betti n{R'}{\wh N}$, and applied
to $\wh\vf$ gives $\betti n{\wh\vf}{\wh N}=\betti n{R'}{\wh N}$.
Bass numbers are treated in a similar fashion.
 \end{subremark}

It is often possible to reduce the study of sequences of Betti numbers
to that of sequences of Bass numbers, and vice versa.

\subsection{Duality}
\label{Duality}
Let $D$ be a dualizing complex for $S$, see \eqref{dualizing complex},
and set
\[
N^\dagger_{S}=\RHom SND
\]

\begin{sublemma}
\label{duality lemma} If $S$ is complete and $R'\to S$ is a
surjective homomorphism, then for each $n\in\BZ$ there is an
isomorphism of $l$-vector spaces
 \[
\Ext^n_{R'}(l,N^\dagger_{S})\cong\Hom_l(\Tor^{R'}_n(l,{N}),l)
 \]
\end{sublemma}

\begin{proof}
Let $D$ be a normalized dualizing complex for $S$. The
isomorphisms
 \[
\Rhom{S}{(\dtensor l{R'}{N})}D \simeq
\Rhom{R'}l{\Rhom{S}{N}{D}}=\Rhom{R'}l{N^\dagger_{S}}
 \]
yield a spectral sequence with
 \[
E^{pq}_2=\Ext^{p}_S(\Tor^{R'}_{q}(l,{N}),D)
\implies \Ext^{p+q}_{R'}(l,N^\dagger_{S})
  \]
The $S$-module $\Tor^{R'}_q(l,{N})$ is a direct sum of copies of $l$,
so $E^{pq}_2=0$ for $p\ne0$ and $E^{0q}\cong\Hom_l(\Tor^{R'}_{q}(l,{N}),l)$
by the defining property of $D$; the assertion follows.
 \end{proof}

\begin{subtheorem}
\label{duality}
For each $n\in\BZ$ the following equalities hold:
 \[
\bass n{\vf}N = \betti n{\grave \vf}{{\wh N}^\dagger_{\wh S}}
          \qquad\text{and}\qquad
  \betti n{\vf}N = \bass n{\grave \vf}{{\wh N}^\dagger_{\wh S}}
 \]
\end{subtheorem}

\begin{proof}
By Remark \eqref{series:cohen} we may assume that $S$ is complete.  As $N$
and $N^{\dagger}_S{}^{\dagger}_S$ are isomorphic, it is enough to justify
the second equality.  If $R\to R'\to S$ is a minimal Cohen factorization
of $\vf$ then the preceding lemma gives the middle equality below:
 \[
\bass n{\vf}{N^\dagger_S} =\bass n{R'}{N^\dagger_S} =\betti
n{R'}{N} =\betti n{\vf}N
 \]
The equality at each end is given by Remark \eqref{series:cohen}.
 \end{proof}

\subsection{Poincar\'e series and Bass series}
\label{pseries}
To study Betti numbers or Bass numbers we often use their generating
functions.  We call the formal Laurent series
 \[
\rP{\vf}{N} = \sum_{n\in\BZ} \betti n\vf N t^n
  \qquad\text{and}\qquad
  \rI{\vf}{N} = \sum_{n\in\BZ} \bass n\vf N t^n
 \]
the \emph{Poincar\'e series} of $N$ over $\vf$ and the \emph{Bass series}
of $N$ over $\vf$, respectively.  When $\vf=\idmap^R$, we speak of the
Poincar\'e series and the Bass series of the complex $N$ over the ring $R$,
and write $\rP{R}{N}$ and $\rI{R}{N}$, respectively.

It is often convenient to work with sets of generators of $\fn$ modulo
$\fm$ that are not necessarily minimal.  The next results provides
the necessary information.

\begin{subproposition}
\label{pseries:extra}
If $\bsy=\{y_1,\dots,y_q\}$ generates $\fn$ modulo $\fm S$, then
\begin{gather*}
\rP{\vf}{N}\cdot(1+t)^{q-\edim\vf}=
\sum_{n\in\BZ}\rank_l\big(\Tor^R_n(k,{\koszul \bsy N})\big)\,t^n
\\
\rI{\vf}{N}\cdot(1+t)^{q-\edim\vf}=
\sum_{n\in\BZ}\rank_l\big(\Ext_R^n(k,{\koszul \bsy N})\big)\,t^{q+n}
\end{gather*}
 \end{subproposition}

\begin{proof}
The two formulas admit similar proofs.  We present the first one.

To simplify notation, for any finite subset $\bsz$ of $\fn$ we set
 \[
\lP{}{\bsz,N}
=\sum_{n\in\BZ}\ell_S\big(\Tor^R_n(k,{\koszul \bsz N})\big)\,t^n
 \]

Let $\bsv=\{v_1,\dots,v_r\}$ be a minimal set of generators of $\fm$.
In $\dcat(S)$ one has
 \[
\koszul{\bsy\sqcup\vf(\bsv)}N\simeq
\koszul{\bsy}N\otimes_S\koszul{\bsv}S\simeq
\dtensor{\koszul{\bsy}N}R{\koszul{\bsv}R}
 \]
These isomorphisms induce the first link in the chain
\begin{align*}
\dtensor kR{\big(\koszul{\bsy\sqcup\vf(\bsv)}N\big)}
&\simeq\dtensor kR{\big(\dtensor{\koszul{\bsy}N}R{\koszul{\bsv}R}\big)}\\
&\simeq\dtensor{\koszul{\bsy}N}R{\big(k\otimes_R{\koszul{\bsv}R}\big)}\\
&\simeq\dtensor{\koszul{\bsy}N}R{\big(\twedge_k\shift k^r\big)}\\
&\simeq\big(\dtensor kR{\koszul{\bsy}N}\big)\otimes_k
{\big(\twedge_k\shift k^r\big)}
\end{align*}
The third one holds because $\bsv\subseteq\fm$ implies that
$k\otimes_R{{\koszul{\bsv}R}}$ has trivial differential; the other two
are standard.  Taking homology and counting ranks over $l$ we get
 \[
\lP{}{{\bsy\sqcup\vf(\bsv)},N} =
\lP{}{\bsy,N}\cdot (1+t)^r
 \]

Let $\bsx=\{x_1,\dots,x_e\}$ be a minimal set of generators of
$\fn$ modulo $\fm S$; clearly, one has $e\le q$.  Let $\bsu$ be
generating set of $\fm$ consisting of $\bsv$ and $q-e$ additional
elements, all equal to $0$.  The formula above then yields
 \[
\lP{}{{\bsx\sqcup\vf(\bsu)},N} =
\lP{}{\bsx,N}\cdot (1+t)^{r+q-e}
 \]

Applying Lemma (\ref{koszul:unique}.3) to $\bsz={\bsy\sqcup\vf(\bsv)}$
and to $\bsz={\bsx\sqcup\vf(\bsu)}$ one gets
$\koszul{\bsy\sqcup\vf(\bsv)}N\cong \koszul{\bsx\sqcup\vf(\bsu)}N$,
so $\lP{}{{\bsy\sqcup\vf(\bsv)},N}=\lP{}{{\bsx\sqcup\vf(\bsu)},N}$.
 \end{proof}

Additional elbow room is provided by the use of arbitrary
Cohen factorizations.

\begin{subproposition}
\label{extra:cohen}
If $R\xra{\dot\vf} R'\xra{\vf'}\wh S$ is any Cohen factorization of
$\grave\vf$, then
\begin{gather*}
\rP{\vf}{N}\cdot(1+t)^{\edim\dot\vf-\edim\vf}=\rP{R'}{\wh N}
\\
\rI{\vf}{N}\cdot(1+t)^{\edim\dot\vf-\edim\vf}=
\rI{R'}{\wh N}
\end{gather*}
 \end{subproposition}

\begin{proof}
By \eqref{series:cohen}, we may assume that $S$ is complete.  In view
of Theorem \eqref{duality} and Lemma \eqref{duality lemma}, it suffices
to prove the expression for Poincar\'e series.

Let $\fm'$ be the maximal ideal of $R'$.  Choose a set $\bsy$ that
minimally generates $\fm'$ modulo $\fm R'$, set $q=\card\bsy$, and
note that $\edim\dot\vf=q$.
The ring $P=R'/\fm R'$ is regular, so $\koszul{\bsy}{P}$ is a resolution
of $l$, hence $\koszul{\bsy}{P}\simeq l$.  {}From this isomorphism and
the associativity of derived tensor products we get
\begin{align*}
\dtensor{l}{R'}N
&\simeq \dtensor{\koszul{\bsy}{P}}{R'}N\\
&\simeq \dtensor{P}{R'}{\koszul{\bsy}N}\\
&\simeq \dtensor{(\dtensor k{R}{R'})}{R'}{\koszul{\bsy}N}\\
&\simeq \dtensor k{R}{\koszul{\bsy}N}
\end{align*}
Passing to homology, we obtain the first equality below:
 \[
\rP{R'}N=
\sum_{n\in\BZ}\rank_l\big(\Tor^{R}_n(k,\koszul{\bsy}N)\big)\,t^n
=\rP{\vf}N\cdot(1+t)^{q-\edim\vf}
 \]
The second equality comes from Proposition \eqref{pseries:extra}.
 \end{proof}

\section{Illustrations}
\label{Illustrations}

Betti numbers or Bass numbers over the identity map of $R$ are equal
to the corresponding numbers over the ring $R$, and decades of research
have demonstrated that ``computing'' the latter invariants is in general
a very difficult task.  Nonetheless, drawing on techniques developed for
the absolute case, or singling out interesting special classes of maps,
it is sometimes possible to express invariants over $\vf$ in closed form,
or to relate them to better understood entities.  In this section we
present several such results, some of which are used later in the paper.

\subsection{Trivial actions}
We record a version of the K\"unneth formula.

\begin{subremark}
\label{kunneth}
If $V$ is a complex of $S$-modules satisfying $\fm V=0$, then there
exist isomorphisms of graded $S$-modules
 \[
\Tor^R(k,V)\cong \Tor^R(k,k)\otimes_k\HH(V)
 \]

Indeed, $V$ is naturally a complex over $S/\fm S$.  If $F$ is a free
resolution of $k$ over $R$, and $\overline F=F/\fm F$, then there are
isomorphisms of complexes of $S$-modules
 \[
F\otimes_RV \cong \overline F\otimes_k V
 \]
Expressing their homology via the K\"unneth formula, one gets
the desired assertion.
 \end{subremark}

\begin{subexample}
\label{series:trivial}
If $\fm N=0$, then setting $s=\edim S$ one has
 \begin{gather*}
(1+t)^{s-\edim\vf}\cdot\rP{\vf}N=\rP{R}k\cdot K^S_N(t)
\\      
(1+t)^{s-\edim\vf}\cdot\rI{\vf}N=\rP{R}k\cdot t^{s}\,K^S_N(t^{-1})
 \end{gather*}
where  $K^S_N(t)$ is the polynomial defined in \eqref{koszul poly}.
Remark \eqref{kunneth}, applied to the complex $V=\rkoszul SN$
yields the expression for $\rP{\vf}N$.  It implies the one for
$\rI{\vf}N$, via Theorems \eqref{duality} and
(\ref{koszul:poly}.4).
 \end{subexample}

\subsection{Regular source rings}
When one of the rings $R$ or $S$ is regular, homological invariants
over the other ring dominate the behavior of the Betti numbers of
$N$ over $\vf$ and its Bass numbers over $\vf$.  Here we consider
the case when the base ring is regular; the case of a regular target
is dealt with in Corollary \eqref{separated:regular}.

\begin{subproposition}
\label{regular source}
If the ring $R$ is regular, then
\[
\rP{\vf}N=K^S_N(t)\cdot(1+t)^{m}
\qquad\text{and}\qquad
\rI{\vf}N=K^S_N(t^{-1})\cdot t^s\,(1+t)^{m}
\]
where $m=\edim R-\edim S+\edim\vf$ and $s=\edim S$.
 \end{subproposition}

\begin{proof}
The complex $K^R$ is a free resolution of $k$ over
$R$, so $\Tor^R(k,\rkoszul SN)\cong\HH(K^R\otimes_R\rkoszul SN)$.
This isomorphism gives the second equality below:
\begin{align*}
\rP{\vf}N\cdot (1+t)^{\edim S-\edim\vf}
&=\sum_{n\in\BZ}\ell_S\big(\Tor^R_n(k,\rkoszul SN)\big)\,t^n\\
&=\sum_{n\in\BZ}\ell_S\big(\HH(K^R\otimes_R\rkoszul SN)\big)\,t^n\\
&=K^S_N(t)\cdot(1+t)^{\edim R}
\end{align*}
The first one comes Proposition \eqref{pseries:extra}, the last
from Lemma (\ref{koszul:unique}.3).

We have proved the expression for $\rP{\vf}N$.  The one for
$\rI{\vf}N$ is obtained from it, by applying Theorems
\eqref{duality} and (\ref{koszul:poly}.4).
 \end{proof}

\begin{subremark}
Recall that $\dcat^{\mathrm f}(S)$ is the derived category of homologically
finite complexes of $S$-modules.  If $R$ is regular, then the proposition
shows that $\rP{\vf}{X}$ and $\rI{\vf}{X}$ are Laurent polynomials
for every $X\in\dcat^{\mathrm f}(S)$.  Conversely, if $\rP{\vf}{V}$
or $\rI{\vf}{V}$ is a Laurent polynomial for some $V\in\dcat^{\mathrm
f}(S)$ with $\HH(V)\ne0$ and $\fm V=0$, then Remark \eqref{series:trivial}
implies that $\rP Rk$ is a polynomial, hence that $R$ is regular.
 \end{subremark}

\subsection{Complete intersection source rings}
We start with some terminology.

\begin{subchunk}
\label{ci rings}
The ring $R$ is \emph{complete intersection} if in some Cohen presentation
$\wh R\cong Q/\fa$, see \eqref{cohen presentation}, the ideal $\fa$
can be generated by a regular set.  When this is the case, the defining
ideal in each (minimal) Cohen presentation of $\wh R$ is generated by a
regular set (of cardinality $\codim R$).  Complete intersection
rings of codimension at most $1$ are also called \emph{hypersurface rings}.
 \end{subchunk}

In the next result we use the convention $0!=1$.

\begin{subtheorem}
\label{betti:ci}
Assume $R$ is complete intersection and set $c=\codim R$.

When $\fd_RN=\infty$ there exist polynomials $b_\pm(x)\in\BQ[x]$ of the
form
 \[
b_\pm(x)= \frac{b_N}{2^c(d-1)!}\cdot x^{d-1}+\text{\ lower order terms\ }
 \]
with an integer $b_N>0$ and $1\le d\le c$, such that
 \[
\betti n{\vf}N=
\begin{cases}
b_+(n) &\text{for all even }n\gg0\,;\\
b_-(n) &\text{for all odd }n\gg0\,.
\end{cases}
 \]
In particular, if $d=1$, then $\betti n{\vf}N=\betti {n-1}{\vf}N
$ for all $n\gg0$.

If $d\ge2$, then there exist polynomials $a_\pm(x)\in\BQ[x]$ of
degree $d-2$ with positive leading coefficients, such that
 \[
\betti {n-1}{\vf}N+a_+(n)\ge\betti n{\vf}N\ge\betti {n-1}{\vf}N+a_-(n)
\quad\text{for all}\quad n\gg0
 \]

The corresponding assertions for the Bass numbers $\bass n{\vf}N$ hold
as well.
 \end{subtheorem}

\begin{proof}
Let $R\to R'\to \wh S$ be a minimal Cohen factorization of
$\grave\vf$. Choose an isomorphism $F'\simeq\wh N$, where $F'$ is a
complex of finite free $R'$-modules with $F'_n=0$ for all $n\le i$,
where $i=\inf\HH(N)$. Set $L=\dd(F_s)$, where $s=\sup\HH(N)$.  The complex
$F_{i\geqslant s}$ is a free resolution of $L$ over $R'$, hence the
second equality below:
 \[
\betti{j+s}{\vf}{N}=\betti{j+s}{R'}{\wh N}=\betti j{R'}L
\quad\text{for all}\quad j>0
 \]
The first equality comes from Remark \eqref{series:cohen}.

As $R'$ is flat over the complete intersection ring $R$, and the ring
$R'/\fm R'$ is regular, $R'$ is complete intersection with $\codim R'=c$.
The properties of the Betti numbers of $L$ over $R$, which translate
into the desired properties of the Betti numbers of $N$ over $\vf$,
are given by \cite[(8.1)]{AGP}.  That the Bass numbers of $M$ have the
corresponding properties is now seen from Theorem \eqref{duality}.
 \end{proof}

\subsection{Weakly regular homomorphisms}
Let $\psi\colon (Q,\fl,h)\to(R,\fm,k)$ be a local homomorphism.

\begin{subchunk}
\label{weakly reg} We say that $\psi$ is \emph{weakly regular} at $\fm$
if the map $\grave\psi$ has a Cohen factorization $Q\to (Q',\fl',k)\xra{\psi'}\wh R$
with $\Ker\psi'$ generated by a \emph{superficial regular set} $\bsf'$,
that is, one whose image in $\fl'/\fl'{}^2$ spans a $k$-subspace of rank
$\card\bsf'$. If $\psi$ is weakly regular at $\fm$, then in each Cohen
factorization of $\grave\psi$ the kernel of the surjection is generated
by a superficial regular set: this follows from the Comparison Theorem
\cite[(1.2)]{AFH}.
 \end{subchunk}

\begin{subtheorem}
\label{flat ass}
If $\psi$ is weakly regular, then
 \[
\rP{\vf}N\cdot (1+t)^{\dim Q-\dim R}=
\rP{\vf\circ\psi}N\cdot(1+t)^{\edim\vf-\edim(\vf\circ\psi)}
 \]
A similar equality holds for the Bass series of $N$.
 \end{subtheorem}

One ingredient of the proof is the following lemma, which extends
\cite[(3.3.5)]{barca}.

\begin{sublemma}
\label{superficial:prop1}
If $f\in\fl\smallsetminus\fl^2$ is a
regular element and $R=Q/(f)$, then
 \[
\rP{\vf\circ\psi}N =\rP{\vf}N \cdot (1+t)
 \]
where $\psi\colon Q\to R$ is the canonical map.
\end{sublemma}

\begin{proof}
We start by fixing some notation.  Let $\bsx$ be a minimal set of
generators of $\fn$ modulo $\fm S$; it also minimally generates $\fn$
modulo $\fl S$, since $\psi$ is surjective.

Let $\koszul fQ$ be the Koszul complex of $f$, and let $E$ denote
the underlying exterior algebra.  Since $f\not\in\fl^2$, one can
extend $\koszul fQ$ to an acyclic closure $Y$ of $k$ over $Q$;
see \cite[(6.3.1)]{barca}.  Moreover, there exists a derivation
$\theta$ on $Y$, compatible with its divided powers structure,
such that $\theta(e_f)=1$, see \cite[(6.3.3)]{barca}. The induced
derivation $\theta\otimes R$ on the DG algebra $Y\otimes_QR$
satisfies $(\theta\otimes R)(e_f\otimes 1)=1$; in addition,
$\dd(e_f\otimes1)=0$. Under these conditions, Andr\'e \cite[Proposition
6]{An2} shows that there exists a complex of $R$-modules $X$ that appears
in an isomorphism of complexes of $R$-modules
 \[
Y\otimes_QR \cong X\oplus \shift X
 \]
Note that $\HH(Y\otimes_QR)\cong \Tor^Q(k,R)\cong k\oplus \shift k$,
so $\HH(X)\cong k$.  In addition, since $Y$ is complex of free
$Q$-modules, the complex of $R$-modules $X$ is free.  Thus, $X$ is a
free resolution of $k$ over $R$. Tensoring the isomorphism above with
$\koszul \bsx N$ yields
\begin{align*}
  Y\otimes_Q\koszul \bsx N &= (Y\otimes_QR)\otimes_R\koszul\bsx N \\
  &\cong \big(X\otimes_R \koszul \bsx N\big) \oplus
  \shift\big(X\otimes_R \koszul \bsx N\big)
\end{align*}
Taking homology one obtains an isomorphism of graded $k$-vector spaces
 \[
\Tor^Q(k,\koszul \bsx N) \cong \Tor^R(k,\koszul \bsx N)\oplus\shift
\Tor^R(k,\koszul \bsx N)
 \]
which yields the desired equality of Poincar\'e series.
 \end{proof}

Also needed for the theorem is a construction from \cite[(4.4)]{AF},
see also \cite[(5.9)]{IS}.

\begin{subconstruction}
\label{big-diagram}
Assume that the rings $R$ and $S$ are complete.

Let $Q\xrightarrow{\Dot{\psi}} Q'\xrightarrow{\psi'}R$ and
$Q'\xrightarrow{\ddot{\varkappa}} Q''\xrightarrow{\varkappa'}S$ be
minimal Cohen factorizations of $\psi$ and $\vf\psi'$, respectively.
The ring $R'=R\otimes_{Q'}Q''$ and the maps
 \[
\psi''=\psi'\otimes_{Q'}Q''\,,\qquad
\dot\vf=R\otimes_{Q'}\ddot\varkappa\,,\qquad
\dot\varkappa=\ddot\varkappa\dot\psi\,,
 \]
fit into a commutative diagram of local homomorphisms
\begin{equation*}
\xymatrixrowsep{2.5pc}
\xymatrixcolsep{2.5pc}
\xymatrix {
  && Q'' \ar@{->}[dr]^{\psi''}  \ar@/^2.5pc/[ddrr]^{\varkappa'}  \\
  & Q' \ar@{->}[dr]^{\psi'} \ar@{->}[ur]^{\ddot\varkappa} && R'
\ar@{->}[dr]^{\vf'}  \\
  Q \ar@{->}[rr]^{\psi} \ar@/^2.5pc/[uurr]^{\dot\varkappa}
\ar@{->}[ur]^{\dot\psi} && R
  \ar@{->}[rr]^{\vf} \ar@{->}[ur]^{\dot\vf} && S }
\end{equation*}
where $\vf'$ is the induced map.  It contains the
following Cohen factorizations:
\begin{enumerate}[\rm(1)]
\item
$R\xra{\dot\vf}R'\xra{\vf'}S$ of $\vf$, which is minimal;
\item
$Q\xra{\dot\varkappa} Q''\xra{\varkappa'} S$ of $\vf\psi$, in which
$\edim\dot{\varkappa} =\edim\vf+\edim\psi$.
\end{enumerate}
\end{subconstruction}

\begin{proof}
[Proof of Theorem \emph{\eqref{flat ass}}]
By Remark \eqref{series:cohen} we may assume $R$ and $S$ complete.
We form the diagram of the construction above and adopt its notation.

Since $\psi$ is weakly regular at $\fm$, the kernel of $\psi'$ is generated
by a superficial regular set $\bsf'$, see \eqref{weakly reg}.  Since
$\psi''=Q''\otimes_{Q'}\psi'$, and $Q''$ is faithfully flat over $Q'$,
the image $\bsf''$ of $\bsf'$ in $Q''$ is a regular set of $\card\bsf'$
elements, and generates the kernel of $\psi''$.  Let $\fl'$ and $\fl''$
be the maximal ideals of $Q'$ and $Q''$, respectively.  As
$Q''/\fl'Q''$ is regular, any minimal set of generators of the ideal
$\fl'Q''=\Ker(Q''\to Q''/\fl'Q'')$ extends to a minimal set of generators
of $\fl''$.  It follows that $\bsf''$ is superficial.

The first equality below holds because $\varkappa'$ is surjective,
see \eqref{surjective}, the second from iterated applications of Lemma
\eqref{superficial:prop1}:
  \[
\rP{Q''}N=\rP{\varkappa'}N=\rP{\vf'}N\cdot (1+t)^{\dim Q-\dim R+\edim\psi}
  \]
Property (2) of Construction \eqref{big-diagram} shows that
$\varkappa'\dot\varkappa$ is a Cohen factorization of
$\vf\,\psi$, with $\edim\dot\varkappa=\edim\vf+\edim\psi$.
Thus, from Proposition \eqref{extra:cohen} we get
\begin{align*}
\rP{Q''}N
&=\rP{\vf\,\psi}N\cdot(1+t)^{\edim\dot\varkappa-\edim(\vf\,\psi)}\\
&=\rP{\vf\,\psi}N\cdot(1+t)^{\edim\vf+\edim\psi-\edim(\vf\,\psi)}
\end{align*}
Comparison of the two expressions for $\rP{Q''}N$ finishes the
proof for Poincar\'e series.

The equality for Bass series is a formal consequence, due to Theorem
\eqref{duality}.
 \end{proof}

\section{Separation}
\label{Separation}

This section introduces one of the main themes of the paper---the
comparison of homological invariants of $N$ over $\vf$ with homological
invariants of $k$ over $R$.

We start by establishing upper bounds on the Poincar\'e series and on the
Bass series of $N$.  The bounds involve series with specific structures:
input from each of the rings $R$ and $S$ and from the map $\vf$ appear
as separate factors.  We say that $N$ is (injectively) separated if the
upper bound is reached, and we turn to the natural problem of identifying
cases when this happens.  Several sufficient conditions are obtained, in
terms of the structure of $N$ over $R$, the properties of the ring $R$,
those of the ring $S$, or the way that $\vf(\fm)$ sits inside $\fn$.
Finally, we obtain detailed information on Betti sequences and Bass
sequences of (injectively) separated complexes.  To do this, we use
the specific form of their Poincar\'e series or Bass series along with
extensive information on $\rP Rk$, available from earlier investigations.

\subsection{Upper bound}
\label{Upper bound}
In the theorem below $K^S_N(t)$ denotes the Laurent polynomial defined
in \eqref{koszul poly}.  The first inequality and its proof are closely
related to results of Lescot in \cite[\S2]{Le}; see Remark \eqref{lescot}.

The symbols $\preccurlyeq$ and $\succcurlyeq$ denote coefficientwise
inequalities of formal Laurent series.

\begin{subtheorem}
\label{ceiling}
Assume $\HH(N)\ne0$, and set $s=\edim S$ and $e=\edim\vf$.

There is an inequality of formal Laurent series
 \[
(1+t)^{s-e}\cdot\rP{\vf}N\preccurlyeq\rP{R}k \cdot K^S_N(t)
 \]
and on both sides the initial term is $(\nu_S\HH_i(N))\,t^i$, where
$i=\inf\HH(N)$.

There is an inequality of formal Laurent series
 \[
(1+t)^{s-e}\cdot\rI{\vf}N
\preccurlyeq\rP{R}k \cdot K^S_N(t^{-1})\,t^{s}
 \]
and on both sides the initial term is $(\type_SN)\,t^{g}$, where
$g=\depth_S N$.
 \end{subtheorem}

\begin{proof}
Let $F$ be a free resolution of $k$ over $R$.  Filter the complex
$F\otimes_R\rkoszul SN$ by its subcomplexes
$F\otimes_R\big(\rkoszul SN_{\les q}\big)$ to get a spectral
sequence
  \[
E_{pq}^2=\Tor_p^R\left(k,\HH_q\rkoszul SN\right)
\implies \Tor_{p+q}^R\left(k,\rkoszul SN\right).
  \]
The $R$-module $\HH_q\rkoszul SN$ is annihilated by $\fm$, so
one has an isomorphism
  \[
\Tor_p^R\left(k,\HH_q\rkoszul SN\right)
  \cong
\Tor_p^R(k,k)\otimes_k\HH_q\rkoszul SN
  \]
It is responsible for the last equality below
 \begin{align*}
\rP{\vf}N\cdot (1+t)^{s-e}
&=
\sum_{n\in\BZ}\rank_l\Tor_{p+q}^R\left(k,\rkoszul SN\right) t^n\\
&\preccurlyeq
\sum_{n\in\BZ}\bigg(\sum_{p+q=n}\rank_l E^2_{p,q}\bigg) t^n\\
&=\rP{R}k\cdot K^S_N(t)
 \end{align*}
while the first equality comes from Lemma \eqref{pseries:extra}.
The convergence of the spectral sequence yields the inequality.

To identify the leading terms, note that Theorem (\ref{koszul:poly}.1)
yields $\HH_q(N)=0$ for $q<i$.  Thus, in the spectral sequence above
$E_{pq}^2=0$ when $p<0$ or $q<i$, hence
 \[
\Tor_{n}^R\left(k,\rkoszul SN\right)\cong
\begin{cases}
0 & \text{for } n<i \\
\HH_i(\rkoszul SN) & \text{for } n=i
\end{cases}
 \]
This shows that $\rank_l\HH_i(\rkoszul SN)\,t^i$ is the initial term
of both series under consideration, and its coefficient is identified
by Theorem (\ref{koszul:poly}.1).

The inequality concerning the Bass series of $N$ and the equality of
initial terms follow from the corresponding statements about the
Poincar\'e series of $N^\dagger$, in view of the equalities
$\rI{\vf}{N}=\rP{\vf}{N^\dagger}$ from Theorem \eqref{duality},
and $K^S_{N^\dagger}(t)=t^{s}K^S_{N}(t^{-1})$ from Theorem
(\ref{koszul:poly}.4).  To obtain the leading term itself use
Theorem (\ref{koszul:poly}.2).
 \end{proof}

The theorem gives a reason to consider the following
classes of complexes.

\begin{subchunk}
\label{separated complexes}
We say that $N$ is \emph{separated} over $\vf$ if
 \[
\rP{\vf}N=\rP{R}k\cdot\frac{K^S_N(t)}{(1+t)^{\edim S}}\cdot(1+t)^{\edim\vf}
 \]
Separation is related to Lescot's \cite{Le} notion of inertness, see
Remark \eqref{lescot}. 

We say that $N$ is \emph{injectively separated} over $\vf$ if
 \[
\rI{\vf}N
=\rP{R}k\cdot \frac{K^S_N(t^{-1})\,t^{\edim S}}{(1+t)^{\edim S}}
\cdot(1+t)^{\edim\vf}
 \]
\end{subchunk}

Comparing the formulas above with those in Remark \eqref{series:trivial}
and Proposition \eqref{regular source}, we exhibit the first instances of
separation.

\begin{subexample}
\label{separated:examples}
If $\fm N=0$, or if the ring $R$ is regular, then $N$ is separated and
injectively separated over $\vf$.
 \end{subexample}

To obtain further examples we introduce a
concept of independent interest.

\subsection{Loewy lengths}
\label{loewy}
Let $X$ be a complex of $S$-modules.  The number
 \[
\lol_S(X) =\inf\{i\in\BN\mid \fn^i\cdot X =0\}
 \]
is the \emph{Loewy length} of $X$ over $S$.
To obtain an invariant in $\dcat(S)$, we introduce
 \[
\hell SX=\inf\{\lol_S(V)\mid V\in\dcat(S)\text{ with }V\simeq X\}
 \]
as the \emph{homotopical Loewy length} of $X$ over $S$.
Evidently, there are inequalities
 \[
\lol_S(\HH(X))\le\hell SX\le\lol_S(X)
 \]
Equalities hold when $X$ is an $S$-module, but not in general; see
Corollary \eqref{separated:regular}.

\begin{sublemma}
\label{hell:product}
If $X$ and $Y$ are complexes of $S$-modules, then
 \begin{gather*}
\hell S{\dtensor XSY} \leq \min\{\hell SX, \hell SY\}\\
\hell S{\Rhom SXY} \leq \min\{\hell SX, \hell SY\}
 \end{gather*}
\end{sublemma}

\begin{proof}
We may assume $\hell SX=i<\infty$, hence $X\simeq V$ and $\fn^i\cdot
V=0$.  Replace $Y$ with a $K$-projective resolution  $F$ to get $\dtensor
XSY\simeq V\otimes_SF$ and $\fn^i\cdot(V\otimes_SF)=0$.  This gives
$\hell S{\dtensor XSY}\le\hell SX$.  The inequality $\hell S{\dtensor
XSY}\le\hell SY$ follows by symmetry.  A similar argument yields the second set of inequalities.
 \end{proof}

Our interest in homotopical Loewy lengths is due to the next result, which
also involves the invariant $\serre S$ introduced in \eqref{serre number}.

\begin{subtheorem}
\label{hell:separated}
The following inequalities hold:
\[
\hell S{\rkoszul SN} \leq \hell S{K^S}\le \serre S
\]

The complex $N$ is separated and injectively separated over $\vf$ if
\[
\vf(\fm)\subseteq\fn^q
\quad\text{where}\quad
q=\hell S{\rkoszul SN}
\]
  \end{subtheorem}

\begin{proof}
Lemma \eqref{hell:product} gives the first one of the desired inequalites.
On the other hand, Proposition \eqref{serre} shows that, for $a=\serre S$,
the complex
 \[
I^a =\quad 0\lra \fn^{a-s} K^S_s \lra \fn^{a-s+1}K^S_{s-1} \lra \cdots
\lra \fn^{a-1}K^S_1 \to \fn^{a}K^S_0\to 0
 \]
is exact.  Thus, $K^S\simeq K^S/I^a$; note that $\fn^a\cdot (K^S/I^a) =0$, and so
$q\le a$.

Assume now $\vf(\fm)\subseteq\fn^q$.  Let $V$ be a complex of $S$-modules
with $\fn^q V=0$ and $V\simeq\rkoszul SN$ in $\dcat(S)$.  This explains
the first and the last isomorphisms below
\begin{align*}
\Tor^R(k,\rkoszul SN) &\cong \Tor^R(k,V) \\
                        &\cong \Tor^R(k,k)\otimes_k\HH(V)\\
                        &\cong \Tor^R(k,k)\otimes_k \HH(\rkoszul SN)
\end{align*}
The second is Remark \eqref{kunneth}.  Computing ranks and using
Proposition \eqref{pseries:extra} one gets separation over $\vf$.
A similar argument yields injective separation.
 \end{proof}

As a corollary, we complement Proposition \eqref{regular source},
see also Example \eqref{separated:examples}.

\begin{subcorollary}
\label{separated:regular}
When $S$ is regular $N$ is separated and injectively separated
over $\vf$.  Also, $S$ is regular if and only if
$\hell S{K^S}\le1$, if and only if $\serre S\le1$.
 \end{subcorollary}

\begin{proof}
Assume first $\hell S{K^S}\le 1$, that is, $K^S\simeq V$ for a complex of
$S$-modules $V$ modules with $\fn\cdot V=0$.  Thus, $V$ is a complex
of $l$-vector spaces, hence $\HH_0(V)=l$ is isomorphic in $\dcat(S)$
to a direct summand of $V$.  One sees that $S$ is regular from:
\[
\fd_Sl\le\fd_SV=\fd_SK^S=\edim S<\infty
\]

Assume next $S$ is regular.  In this case $\grass\fn S$ is a polynomial
algebra, so $\HH(K^{\grass\fn S})\cong l$, and hence $\serre S\le1$ by
Proposition \eqref{serre}.  Now Theorem \eqref{hell:separated} yields an
inequality $\hell S{K^S}\le1$ and the assertions about $N$.
 \end{proof}

\subsection{Properties}
The first non-zero Betti number and the first non-zero Bass number
of $N$ are computed in Theorem \eqref{ceiling}.  For separated
complexes detailed information is also available on subsequent
numbers in each sequence.

\begin{subtheorem}
\label{separated}
Assume $\HH(N)\ne0$, and set $i=\inf\HH(N)$ and $m=\nu_S\HH_i(N)$.

If $N$ is separated over $\vf$, then its Betti numbers have the following
properties.
\begin{enumerate}[\quad\rm(1)]
\item
When $R$ is not a hypersurface ring there are inequalities
\begin{align*}
\betti {n}{\vf}N\ge
\begin{cases}
\betti{i}{\vf}N&\text{for }n=i+1\\
\betti {n-1}{\vf}N+m&\text{for }n\ge i+2
\end{cases}
\end{align*}
If $\codim S\ge2$ or $\edim R+\edim\vf>\edim S$, then also
 \[
\betti{i+1}{\vf}N\ge \betti{i}{\vf}N+m
 \]
\item
When $R$ is not complete intersection there is a real number
$b>1$ such that
 \[
\betti {n}{\vf}N\ge m\,b^{n-i} \quad\text{for all}\quad n\ge i+2
 \]
\end{enumerate}

If $N$ is injectively separated over $\vf$, then its Bass numbers over
$\vf$ satisfy similar inequalities, obtained from those
above by replacing $\betti{n}{\vf}N$ with $\bass{n}{\vf}N$, the number
$i=\inf\HH(N)$ with $\depth_SN$, and $m=\nu_S\HH_i(N)$ with $\type_SN$.
 \end{subtheorem}

The theorem is proved at the end of this section.

\begin{subremark}
\label{euler} An inequality  $\edim R+\edim\vf\ge\edim S$ always
holds, and there is equality if and only if some (respectively,
each) minimal Cohen factorization $R\to R'\to\wh S$ of $\grave\vf$
satisfies $\edim R'=\edim S$.

Indeed, if $\bsv$ minimally generates $\fm$ and $\bsx$ minimally
generates $\fn$ modulo $\fm S$, then $\vf(\bsv)\cup\bsx$ generates
$\fn$, hence the inequality. In a minimal Cohen factorization
$\edim R'=\edim R+\edim\vf$, so equality holds if and only if
$\edim R'=\edim S$.
 \end{subremark}

Next we show that the results in the theorem are optimal.

\begin{subexample}
\label{emmy}
Let $k$ be a field and set $R=k[[x,y]]/(x^2,xy)$.  One then has
 \[
\rP Rk=\frac{(1+t)^2}{1-2t^2-t^3}=\frac{1+t}{1-t-t^2}
 \]
see \cite[(5.3.4)]{barca}.  Set $S=R$ and let $N$ be the $S$-module
$S/(x)$.  Set $\vf=\idmap^R$, so $\rP\vf N=\rP RN$.  The
exact sequence $0\to k\to R\to N\to 0$ yields
 \[
\rP RN=1+t\rP Rk=\frac{1}{1-t-t^2}
 \]
Furthermore, one has $K^S_N(t)=1+t$, e.g.\ from
Theorem (\ref{koszul:poly}.3) with $T=k[[x,y]]$.

The computations above show that $N$ is separated. However $\betti
1\vf N=\betti0\vf N=1$, so equality holds in (1), and the
inequality in (2) fails for $n=i+1$.
 \end{subexample}

\begin{subremark}
\label{cohen:separation}
If $R\xra{\dot\vf}R'\xra{\vf'}\wh S$ is a minimal Cohen factorization
of $\wh S$, then $N$ is separated over $\vf$ if and only if $\wh N$
is separated over $\vf'$.

Indeed, set $e=\edim\vf$.  As $\dot\vf$ is flat and $R'/\fm R'$ is regular,
it is not hard to show $\rP{R'}l=\rP Rk(1+t)^{e}$.
The equality $\rP{\vf}N=\rP{\vf'}{\wh N}$ yields our assertion.
 \end{subremark}

We are ready to compare separation with inertness, as defined by
Lescot \cite{Le}.

\begin{subremark}
\label{lescot}
An $S$-module $L$ is \emph{inert} if $(1+t)^{\edim S}\cdot\rP SL=\rP
Sk\cdot K^S_L(t)$, see \cite[(2.5)]{Le}.  This is precisely the condition
that $L$ is separated over $\idmap^S$.

Suppose $\edim S=\edim R+\edim\vf$.  Remark \eqref{euler} shows that
$\grave\vf$ has a minimal Cohen factorization with $\edim R'=\edim S$,
so $K^S_L(t)=K^{R'}_L(t)$.  By the preceding remark, $L$ is separated
over $\vf$ if and only if it is inert over $R'$.

When $\vf$ induces the identity on $k$ and the ring $S/\fm S$ is artinian,
Lescot says that $L$ is \emph{inert through} $\vf$ if $\rP Sk\cdot\rP
RL=\rP SL\cdot\rP Rk$.  This condition is very different from separation
over $\vf$: it can be shown that $L\ne0$ has both properties if and only
if it is inert over $R$, the map $\vf$ is flat, and $\fn=\fm S$.
 \end{subremark}

Next we recall background information on the homology of $k$.

\begin{subchunk}
\label{background}
Set $r=\edim R$, and let $p$ denote minimal number of generators of
the ideal $\fa$ in a minimal Cohen presentation $\wh R\cong Q/\fa$.
The series $\rP{R}k$ can be written as
 \[
\rP{R}k
 =\frac{(1+t)^r}{(1-t^2)^p}\cdot F(t)
\quad\text{where}\quad
F(t)=\frac{\prod_{i=1}^\infty(1+t^{2i+1})^{\varepsilon_{2i+1}(R)}}
{\prod_{i=1}^\infty(1-t^{2i+2})^{\varepsilon_{2i+2}(R)}}
 \]
for uniquely determined integers $\varepsilon_n(R)\ge0$ for $n\ge3$,
see \cite[(7.1.4)]{barca}.

\begin{enumerate}[\rm(1)]
\item
The ring $R$ is regular if and only $p=0$; in this case $F(t)=1$.
\item
The ring $R$ is complete intersection if and only $\varepsilon_3(R)=0$,
if and only if $F(t)\!=\!1.$
\item
If $R$ is not complete intersection, then there exists an
integer $s_1\ge2$ and a sequence of integers $(i_j)_{j=1}^\infty$
with $2\le i_j\le \edim R+1$ for all $j$, such that
\[
s_{j+1}=i_j(s_j-1)+2
\quad\text{satisfy}\quad
\varepsilon_{s_j}(R)\ge a^{s_j}
\]
for some real number $a>1$.

Of the assertions above (1) is clear, (2) is due to Assmus and
Tate, see \cite[(7.3.3)]{barca}, and (3) is due to Avramov, see
\cite[(8.2.3)]{barca}.
 \end{enumerate}
 \end{subchunk}

To prove the next lemma we abstract an argument from the proof of
\cite[(8.2.1)]{barca}.

\begin{sublemma}
\label{calculus} Assume $R$ is not a complete intersection and let
$F(t)$ be as above.
\begin{enumerate}[\quad\rm(1)]
\item
The following inequalities hold:
 \[
\frac{F(t)}{(1-t^2)}\succcurlyeq\frac{1+t^3}{(1-t^2)}=
1+\sum_{h=2}^\infty t^h
 \]
\item
There exists a real number $b>1$, such that the following
inequality holds:
 \[
\frac{F(t)}{(1-t)(1-t^2)}\succcurlyeq 1+t+\sum_{h=2}^\infty
b^h\,t^h
 \]
\end{enumerate}
\end{sublemma}

\begin{proof}
(1) is clear from (\ref{background}.2).

(2) Set $\sum_{h\in\BZ}a_ht^h=F(t)/(1-t^2)$ and
$\sum_{h\in\BZ}b_ht^h=F(t)/((1-t)(1-t^2))$.

Since $b_h=a_h+\cdots+a_{0}$ for each $h\ge0$, we get $b_1\ge
b_0\ge1$ and
\begin{gather*}
b_{h+1}=\bigg(\sum_{j=0}^{h+1}a_j\bigg)>\bigg(\sum_{j=0}^{h}a_j\bigg)=b_h
\qquad\text{for all}\qquad h\ge 1
\end{gather*}
Let $s_1,s_2,\dots$ be the numbers from (\ref{background}.3) and
set $r=\edim R$. The number
  \[
b=
\min\bigg\{\,\sqrt[r+1]{\mathstrut a}\,,\, \sqrt{b_2}
\,,\dots\,,\,\sqrt[s_1] b_{s_1}\bigg\}
  \]
satisfies $b_h\ge b^h>1$ for $s_1\ge h\ge2$. If $s_{j+1}\ge h>s_j$
with $j\ge1$, then
  \[
b_h>b_{s_j}\ge\varepsilon_{s_j}(R)\ge a^{s_j}\ge b^{(r+1)s_j}>
b^{s_{j+1}}>b^h
  \]
so we see that $b_h\ge b^h$ holds for all $h\ge 2$.  This is the
desired inequality.
\end{proof}

\begin{proof}[Proof of Theorem \emph{\eqref{separated}}]
We adopt the notation of Lemma \eqref{calculus}, set
 \[
d=\edim R-\edim S+\edim\vf
 \]
Remark \eqref{euler} and the hypothesis that $R$ is not a
hypersurface yield
 \[
d\ge0 \qquad\text{and}\qquad p\ge2 \tag{$*$}
 \]
Theorem (\ref{koszul:poly}.7) provides a Laurent polynomial
$L(t)\succcurlyeq m\,t^i$, such that
 \[
K^S_N(t)(1+t)^{d}=L(t)(1+t)^{d+1} \tag{$\dagger$}
 \]

(1)  From the discussion above we obtain the relation
\[
\sum_{n=i}^\infty\big(\betti n{\vf}N-\betti{n-1}{\vf}N\big)\,t^n
=(1-t)\cdot\rP{\vf}N= \frac{L(t)(1+t)^{d}}{(1-t^2)^{p-1}}\cdot
F(t)
\]

If $\codim S\ge2$, then Theorem (\ref{koszul:poly}.5) yields
$L(t)\succcurlyeq m\,t^i\,(1+t)$, in particular,
 \[
L(t)(1+t)^{d}\succcurlyeq m\,t^i\,(1+t)
 \]
The same inequality holds when $d\ge1$. Thus, in these cases we get
\begin{align*}
\sum_{n=i}^\infty\big(\betti n{\vf}N-\betti{n-1}{\vf}N\big)\,t^n
&=\frac{L(t)(1+t)^{d}}{(1-t^2)}\cdot\frac{F(t)}{(1-t^2)^{p-2}}\\
&\succcurlyeq \frac{m\,t^i\,(1+t)}{(1-t^2)}\\
&=m\sum_{n=i}^\infty t^n
\end{align*}
This implies $\betti {n}{\vf}N\ge\betti{n-1}{\vf}N+m$ for all
$n\ge i+1$, as desired.

Now we assume $\codim S=1$ and $d=0$.  By Remark \eqref{euler},
$\grave\vf$ has a minimal Cohen factorization $R\to R'\to\wh S$ with
$\edim R'=\edim S$, hence
 \[
\codim R=\codim R'=\edim R'-\dim R'\le\edim S-\dim S=\codim S= 1
 \]
Since $R$ is not a hypersurface by hypothesis, the inequality above
bars it from being complete intersection. Using Theorem
(\ref{koszul:poly}.1) and Lemma (\ref{calculus}.1) we obtain
 \[
\sum_{n=i}^\infty (\betti n{\vf}N-\betti{n-1}{\vf}N)\,t^n
=\frac{L(t)}{(1-t^2)^{p-2}}\cdot\frac{F(t)}{(1-t^2)}
\succcurlyeq m\,t^i\cdot(1+\sum_{h=2}^\infty t^h)
 \]
This yields $\betti{i+1}\vf N\ge\betti i\vf N$ and $\betti
{n}{\vf}N>\betti{n-1}{\vf}N)+m$ for all $n\ge i+2$.

(2) Formula ($\dagger$) gives the equality below:
\begin{align*}
\rP{\vf}N
&=L(t)\cdot\frac{(1+t)^{d}}{(1-t^2)^{p-2}}
\cdot\frac{F(t)}{(1-t)(1-t^2)}\\
&\succcurlyeq m\,t^i\cdot\bigg(1+t+\sum_{h=2}^\infty b^h\,t^h\bigg)
\end{align*}
The inequality, involving a real number $b>1$, comes from Theorem
(\ref{koszul:poly}.7), formula ($*$), and Lemma
(\ref{calculus}.2).  Thus, $\betti n{\vf}N\ge m\,b^{n-i}$ holds
for all $n\ge i+2.$
 \end{proof}

\section{Asymptotes}
\label{Asymptotes}

The data encoded in the sequences of Betti numbers and of Bass numbers
of $N$ are often too detailed to decipher.  Many  results suggest that
it is the \emph{asymptotic behavior} of these sequence that carries an
understandable algebraic significance.  In this section we introduce and
initiate the study of a pair of numerical invariants that evaluate the
(co)homological nature of $N$ over $R$ by measuring the asymptotic growth
of the sequences of its Betti numbers and of its Bass numbers over $\vf$.

We use two scales to measure the rate at which these numbers grow: a
polynomial one, leading to the notion of complexity, and an exponential
one, yielding that of curvature.  One reason for restricting to these
scales is that these numbers grow at most exponentially, and there
is no example with rate of growth intermediate between polynomial and
exponential.  A second reason comes from Theorem \eqref{ceiling}:
it is natural to compare the homological sequences of $N$ over $\vf$ to
$(\betti nRk)_{n\ges0}$, and this sequence is known to have either polynomial
or exponential growth.

\subsection{Complexities and curvatures}
The \emph{complexity} of $N$ over $\vf$ is the number
\begin{gather*}
\cxy_{\vf}N =
\inf\left\{d \in\BN \left|
\begin{gathered}
\text{there exists a number $c\in\BR$ such that}\\
\betti n{\vf}N\leq c n^{d-1}\text{ for all $n\gg0$}
\end{gathered}
\right\}\right.
\end{gather*}
The \emph{curvature} of $N$ over $\vf$ is the number
 \[
\curv_{\vf}N = \limsup_n \sqrt[n]{\betti n{\vf}N}
 \]
Similar formulas, in which Betti numbers are replaced by Bass numbers,
define the \emph{injective complexity} $\injcxy_{\vf}N$ and the
\emph{injective curvature} $\injcurv_{\vf}N$ of $N$ over $\vf$.

When $\vf=\idmap^R$ one speaks of complexities and curvatures of $N$
over $R$, modifying the notation accordingly to $\cxy_R(N)$, etc.,
see \cite{Av:extremal}.  {}From Proposition \eqref{extra:cohen} one gets:

\begin{subremark}
\label{cx cohen}
If $R\xra{\dot\vf}R'\xra{\vf'}\wh S$ is any Cohen factorization of
$\vf'$, then
\begin{gather*}
\cxy_{\vf}{N} = \cxy_{\grave\vf}{\wh N} =
\cxy_{\wh\vf}{\wh N}=\cxy_{R'}(\wh N) \\
\curv_{\vf}{N} = \curv_{\grave\vf}{\wh N} =
\curv_{\wh\vf}{\wh N}=\curv_{R'}(\wh N)
\end{gather*}
The corresponding injective invariants satisfy analogous relations.
\end{subremark}

The equalities below often allow one to restrict proofs to projective
invariants.

\begin{subremark}
\label{cx duality} If $D$ is a dualizing complex of $\wh S$ and
$\wh N{}^\dagger= \RHom{\wh S}{\wh N}D$, then
 \[
\injcxy_{\vf}{N} = \cxy_{\grave\vf}{\wh N^\dagger}
\qquad\text{and}\qquad
\injcurv_{\vf}{N} = \curv_{\grave\vf}{\wh N^\dagger}
 \]
Indeed, this is an immediate consequence of Theorem \eqref{duality}.
\end{subremark}

We list some basic dependencies among homological invariants of $N$.

\begin{subproposition}
\label{dependencies}
The numbers $\cxy_\vf N$ and $\curv_\vf N$ satisfy the relations below.
\begin{enumerate}[\rm\quad(1)]
\item
$\fd_R N<\infty$ if and only if $\cxy_{\vf}N=0$, if and only
if $\curv_\vf N=0$.
\item
If $\fd_R N=\infty$, then $\cxy_{\vf}N\ge1$ and $\curv_\vf N\geq1$.
\item
If $1\le\cxy_{\vf}N<\infty$, then $\curv_\vf N=1$.
\item
If $\curv_\vf N>1$, then there exist an infinite sequence
$n_1<n_2<\dots$ of integers and a real number $b>1$, such that
$\betti{n_i}\vf N\ge b^{n_i}$ holds for all $i\ge1$.
\item
$\curv_\vf N\le\curv_Rk<\infty$.
\end{enumerate}
The corresponding injective invariants satisfy analogous relations.
\end{subproposition}

\begin{Remark}
Part (5) shows, in particular, that there exists a positive real number
$c$, such that $\betti{n}\vf N\le c^{n}$ and $\bass{n}\vf N\le c^{n}$ hold
for all $n\ge1$.
 \end{Remark}

\begin{proof}
By Remark \eqref{cx duality} we may restrict to projective invariants.
Part (1) comes from Remark \eqref{series:cohen} and Corollary \eqref{dim
properties}.  Parts (2), (3), and (4) are clear.  In (5) the first
inequality comes from Theorem \eqref{ceiling}, see \eqref{polyfactor}.
The second inequality is known, see e.g.\ \cite[(4.2.3)]{barca}; a
self-contained proof is given in Corollary \eqref{finite curv}.
 \end{proof}

Complexity may be infinite, even when $\vf$ is the identity map.

\begin{subexample}
\label{trivial extension}
For the local ring $R=k[x,y]/(x^2,xy,y^2)$ and the $R$-module $N=R/(x,y)$
one has $\betti nRN=2^n$, so $\cxy_RN=\infty$.
 \end{subexample}

We determine when \emph{all} homologically finite complexes have finite
complexity.  Recall that $\dcat^{\rm f}(S)$ denotes the derived category
of homologically finite complexes.

\pagebreak

\begin{subtheorem}
\label{generalcichar}
Set $c=\codim R$.  The following conditions are equivalent.
\begin{enumerate}[{\quad\rm(i)}]
\item $R$ is a local complete intersection ring.
\item
$\cxy_\vf{X}\leq c$ for each $X\in\dcat^{\rm f}(S)$.
\item
$\cxy_\vf l=c$
\item
$\curv_\vf l\leq 1$
\end{enumerate} They are also equivalent to those obtained by replacing
Poincar\'e series, complexity, curvature with Bass series, injective
complexity, injective curvature, respectively.
 \end{subtheorem}

\begin{proof}
As $l$ is separated over $\vf$ by Example \eqref{separated:examples}, 
from \eqref{polyfactor} one gets
\[
\cxy_{\vf}l=\cxy_Rk
 \qquad\text{and}\qquad
\curv_{\vf}l=\curv_Rk
\]

(i) $\implies$ (iii).  As $\rP Rk=(1+t)^{\edim R}/(1-t^2)^c$, see
(\ref{background}.2), one has $\cxy_Rk=c$ by (\ref{polyfactor}.4), 
so the equalities above
yield $\cxy_{\vf}l=c$.

(iii) $\implies$ (ii) holds by Theorem \eqref{ceiling}.

(ii) $\implies$ (iii) $\implies$ (iv) hold by definition.  

(iv) $\implies$ (i) is a consequence of the formulas above and
(\ref{background}.3).
 \end{proof}

Our results on complexities and curvatures are often deduced from
coefficientwise inequalities of formal Laurent series.  We collect some
simple rules of operation.

\begin{subchunk}
 \label{polyfactor}
Let $a(t)=\sum a_n t^n$ be a formal Laurent series with non-negative
integer coefficients.  Extending notation, let $\cxy(a(t))$ denote the
least integer $d$ such that $a_n\le cn^{d-1}$ holds for some $c\in\BR$
and all $n\gg0$, and set $\curv(a(t))=\limsup_n\sqrt[n]{a_n}$.

Let $b(t)$ be a formal Laurent series with non-negative integer
coefficients.
\begin{enumerate}[\quad\rm(1)]
\item
If $a(t)\preccurlyeq b(t)$, then
 \[
\cxy(a(t))\le\cxy(b(t))
\quad\text{and}\quad
\curv(a(t))\le\curv(b(t))
 \]
\item
The following inequalities hold:
\begin{align*}
\min\{\cxy(a(t))\,,\,\cxy(b(t))\}
&\le\cxy(a(t)\cdot b(t))
\le\cxy(a(t))+\cxy(b(t))
\\
\min\{\curv(a(t))\,,\,\curv(b(t))\}
&\le\curv(a(t)\cdot b(t))\\
&\le\max\{\curv(a(t))\,,\,\curv(b(t))\}
\end{align*}
\end{enumerate}
When $a(t)$ represents a rational function the following hold as well.
\begin{enumerate}[\quad\rm(1)]
\item[(3)]
$\curv(a(t))$ is finite.
\item[(4)]
$\cxy(a(t))$ is finite if and only if $a(t)$ converges in the unit circle;
when this is the case $\cxy(a(t))$ is equal to the order of the pole
of $a(t)$ at $t=1$;
\item[(5)]
If $b(t)$ represents a rational function, and $\cxy(a(t))$,
$\cxy(b(t))$ are finite, then
 \[
\cxy(a(t)\cdot b(t))=\cxy(a(t))+\cxy(b(t))
 \]
\end{enumerate}

Indeed, (1) and the first three inequalities in (2) follow from
the definitions.  Since $\curv(a(t))$ is the inverse of the radius
of convergence of the series $a(t)$, the last inequality in (2) is
a reformulation of a well known property of power series, and (3) is
clear. For (4), see e.g.\ \cite[(2.4)]{Av:msri}; (5) follows from (4).
 \end{subchunk}

\subsection{Reductions}

The determination of complexity or curvature can sometimes be simplified
using Koszul complexes.  The following result is a step in that direction.

\begin{subproposition}
\label{cxy koszul.1}
If $\bsv$ is a finite subset of $\fn$, then
 \[
\cxy_\vf N = \cxy_\vf(\koszul\bsv N)
\qquad\text{and}\qquad
\curv_\vf N = \curv_\vf(\koszul\bsv N)
 \]
Similar equalities hold for the corresponding injective invariants.
\end{subproposition}

\begin{proof}
In view of the isomorphism of complexes $\koszul{\bsx}{\koszul{\bsv}N}
\cong\koszul{\bsx,\bsv}N$, Proposition \eqref{pseries:extra} applied
to the set $\bsy=\bsx\sqcup\bsv$ yields
 \[
\rP{\vf}{\koszul{\bsv}N}=\rP{\vf}N\cdot(1+t)^{\card\bsv}
 \]
The equalities for complexity and curvature now result from
\eqref{polyfactor}, and those for their injective counterparts
follow via \eqref{dualizing koszul} and Theorem \eqref{duality}.
 \end{proof}

\begin{subcorollary}
\label{cxy koszul.2}
If $L$ is an $S$-module and $\bsv$ is an $L$-regular subset of $S$, then
 \[
\cxy_\vf(L/\bsv L) = \cxy_\vf L
\qquad\text{and}\qquad
\curv_\vf(L/\bsv L) = \curv_\vf L
 \]
Similar equalities hold for the corresponding injective invariants.
\end{subcorollary}

\begin{proof}
In this case $L/\bsv L$ is isomorphic to $\koszul\bsv L$ in the derived
category of $S$.
 \end{proof}

As a special case of the next theorem, the Koszul complex on any system of
parameters for the $S$-module $\HH(N)/\fm\HH(N)$ can be used to determine
the complexity or the curvature of $N$ over $\vf$.

\begin{subtheorem}
\label{ass:finite}
Let $\bsv$ be a finite subset of $\fn$.  If the $S$-module
 \[
\frac{\HH(N)}{\bsv\HH(N) +\fm \HH(N)}
 \]
has finite length, then
\begin{gather*}
\cxy_{\vf} N =
\cxy\bigg(\sum_{n\in\BZ}\ell_S \Tor^R_n(k,\koszul{\bsv}N)\,t^n \bigg) \\
\curv_{\vf} N =
\limsup_n\sqrt[n]{\ell_S \Tor^R_n(k,\koszul{\bsv}N)}
\end{gather*}
Similar equalities involving $\ell_S \Ext_R^n(k,\koszul{\bsv}N)$
express $\injcxy_{\vf}N$ and $\injcurv_{\vf}N$.
\end{subtheorem}

\begin{subremark}
\label{ass:finite_remark}
When $\ell_S\big(\HH(N)/\fm \HH(N)\big)$ is finite one may choose
$\bsv=\varnothing$, so that $\koszul{\bsv}N=N$; this applies, in
particular, when the ring $S/\fm S$ is artinian.
 \end{subremark}

\begin{proof}[Proof of Theorem \emph{\eqref{ass:finite}}]
The proofs of the formulas for complexity and curvature are similar to
those for their injective counterparts; for once, we verify the latter.

Set $K=\koszul\bsv N$.  By Proposition \eqref{cxy koszul.1}, it
suffices to prove that the expressions on the right hand sides of the
desired equalities are equal to $\injcxy_\vf K$ and $\injcurv_\vf K$,
respectively.  Note the equalities of supports of $S$-modules
\begin{align*}
\Supp_S\bigg(\frac{\HH(N)}{\bsv\HH(N) +\fm \HH(N)}\bigg)
&=\Supp_S(S/\fm S)\cap\Supp_S(S/\bsv S)\cap\Supp_S\HH(N)\\
&=\Supp_S(S/\fm S)\cap\Supp_S\HH(K)\\
&=\Supp_S\big(\!\HH(K)/\fm\HH(K)\big)
\end{align*}
where the middle one comes from Lemma \eqref{finite koszul} and the other
two are standard.  Therefore, the length of the $S$-module $\HH(K)/\fm
\HH(K)$ is finite.

We claim that there exists a positive integer $v$ for which there are
inequalities:
 \[
\rI {\vf}K t^{\edim\vf} \preccurlyeq
\sum_{n\in\BZ}\ell_S\big(\Ext_R^n(k,K)\big)\,t^n\cdot (1+t)^{\edim\vf} \preccurlyeq
\rI{\vf}K\cdot (vt)^{\edim\vf}
 \]
Indeed, let $\bsx=\{x_1,\dots,x_e\}$ be a minimal generating set
of $\fn$ modulo $\fm S$.  Set $K^{(j)}=\koszul{x_1,\dots,x_j}K$ for
$j=0,\dots,e$, and note that $K^{(j+1)}$ is the mapping cone of the
morphism $\lambda^{(j)}\colon K^{(j)}\to K^{(j)}$ given by multiplication
with $x_{j+1}$, cf.\ \eqref{cone koszul}.

By \eqref{primary}, our hypothesis implies that there is an integer
$v\ge1$ such that $\fn^v$ is contained in $\fm S+\hann SK$.  This ideal
is contained in $\fm S+\hann S{K^{(j)}}$ for each $j$, so by Lemma
\eqref{finite betti} each $S$-module $\HH(\dtensor kR{K^{(j)}})$ is
finite and annihilated by $\fn^v$.  Furthermore, there is a triangle
\[
\dtensor kR{K^{(j)}}\xra{\dtensor kR{\lambda^{(j)}}}\dtensor kR{K^{(j)}}
\lra \dtensor kR{K^{(j+1)}}\lra
\]
In view of Lemma \eqref{rod:sandwich}, the formal  Laurent series
 \[
I^{(j)}(t)=\sum_{n\in\BZ}\ell_S(\Ext_R^n(k,K^{(j)})\big)\,t^n
 \]
satisfy coefficientwise inequalities
 \[
I^{(j+1)}(t) \preccurlyeq I^{(j)}(t)\cdot(1+t) \preccurlyeq
I^{(j+1)}(t)\cdot v
 \]
It remains to note that $I^{(e)}(t)= \rI{\vf}K\cdot t^{\edim\vf}$.
 \end{proof}

\section{Comparisons}
\label{Comparisons}

The source ring $R$ and the target ring $S$ of the homomorphism $\vf$ play
different roles in the definitions of Betti numbers and Bass numbers,
and hence in the characteristics of $N$ introduced above: The impact of
the $R$-module structure of $N$ is through the entire resolution of $k$,
while that of its $S$-module structure is limited to the Koszul complex
$\rkoszul SN$.  We take a closer look at the role of $S$.

For the rest of the section we fix a second local homomorphism
 \[
{\wt\vf}\colon (R,\fm,k)\to({\wt S},{\wt\fn},{\wt l})
 \]
and we assume that $N$ also has a structure of a homologically finite
complex of ${\wt S}$-modules, and that the action of $R$ through ${\wt\vf}$
coincides with that through $\vf$.

Our main conclusion is that if the actions of $S$ and ${\wt S}$ on $N$ commute,
then the (injective) complexity or curvature of $N$ over $\vf$ is equal to
that over ${\wt\vf}$.  This may be surprising, as it is in general impossible
to express the Betti numbers or the Bass numbers over one of the maps
in terms of those over the other.  As an application, we show that in
the presence of commuting actions important invariants of $N$ over $S$,
such as depth and Krull dimension, are equal to those over ${\wt S}$.

\subsection{Asymptotic invariants}
We consider the following natural questions.

\begin{subquestion}
Do the equalities below always hold:
 \[
\cxy_{\vf} N=\cxy_{{\wt\vf}}N
\qquad\text{and}\qquad
\curv_{\vf} N=\curv_{{\wt\vf}}N\ ?
 \]
\end{subquestion}

\begin{subquestion}
What about the corresponding injective invariants?
 \end{subquestion}

The next result settles an important special case.  While we do not know
the answer in general, Example \eqref{rod:projdim} raises the possibility
of a negative answer.

\begin{subtheorem}
\label{bimodules:asymptotes}
If the actions of $S$ and ${\wt S}$ on $N$ commute, then
 \[
\cxy_{\vf} N=\cxy_{{\wt\vf}}N
\qquad\text{and}\qquad
\curv_{\vf} N=\curv_{{\wt\vf}}N
 \]
Similar equalities also hold for the injective invariants of $N$.
 \end{subtheorem}

The special case where $\vf'=\idmap^R$ is of independent interest.

\begin{subcorollary}
\label{ass:overR}
If $\HH(N)$ is finite over $R$, then
 \[
\cxy_\vf N=\cxy_RN
\qquad\text{and}\qquad
 \curv_{\vf} N =\curv_R N
 \]
Similar equalities for $\injcxy N$ and $\injcurv N$ also hold.\qed
 \end{subcorollary}

The theorem is a consequence of the relations of formal Laurent series
established in the proposition below.  Indeed, \eqref{polyfactor}
yields $\cxy_\vf N\le\cxy_{{\wt\vf}}N$ and $\curv_{\vf} N \le\curv_{{\wt\vf}}
N$, and the converse inequalities follow by symmetry.

\begin{subproposition}
\label{bimodules:series} If the actions of $S$ and ${\wt S}$ commute,
then there exists a positive integer $w$, such that
\begin{gather*}
\rP{\vf}N\cdot(1+t)^{\edim{\wt\vf}}
\preccurlyeq\rP{{\wt\vf}}N\cdot w\,(1+t)^{\edim\vf}\\
\rI{\vf}N\cdot(1+t)^{\edim{\wt\vf}}
\preccurlyeq\rI{{\wt\vf}}N\cdot w\,(1+t)^{\edim\vf}
\end{gather*}
\end{subproposition}

The proposition is proved at the end of the section.  Next we use it to
compare patterns of vanishing of Betti numbers and Bass numbers over $\vf$
and over ${\wt\vf}$.

\subsection{Homological dimensions}
We define the \emph{projective dimension} of $N$ over $\vf$ and the
\emph{injective dimension} of $N$ over $\vf$, respectively, to be
the numbers
\begin{gather*}
  \pd_\vf N=\sup\{n\in\BZ \mid\betti n\vf N\ne0\}-\edim\vf\\
  \id_\vf N=\sup\{n\in\BZ \mid\bass n\vf N\ne0\}-\edim\vf
\end{gather*}
These expressions have been chosen in view of the characterizations
of homological dimensions over $R$ recalled in \eqref{hdim:classical}.
The shifts by $\edim\vf$ appear because the modified invariants appear
to carry a more transparent algebraic meaning.

The equalities below come from Remark \eqref{series:cohen}, and reconcile
our notion of projective dimension over $\vf$ with that defined in
\cite[(3.5)]{IS} via Cohen factorizations.

\begin{subchunk}
\label{dim cohen}
Let $\grave\vf\colon R\to \wh S$ be the composition of $\vf$ with
the completion $S\to \wh S$. If $R\to R'\to \wh S$ is a minimal Cohen
factorization of $\grave\vf$, then
 \begin{gather*}
\pd_\vf N=\pd_{\grave\vf}\wh N=\pd_{\wh\vf}\wh N=\fd_{R'}\wh N-\edim\vf \\
\id_\vf N=\id_{\grave\vf}\wh N=\id_{\wh\vf}\wh N=\id_{R'}\wh N-\edim\vf
 \end{gather*} \end{subchunk}

Combining these equalities with Corollary \eqref{dim properties}, we obtain

\begin{subcorollary}
\label{phi dim properties}
The following (in)equalities hold:
\begin{xxalignat}{3}
&\hphantom{\square}
&\fd_{R}N-\edim\vf
&\leq \pd_{\vf} N \leq \fd_{R}N
&&\hphantom{\square}\\
&\hphantom{\square}
&\id_{\vf}N
&=\id_RN
&&{\square}
\end{xxalignat}
 \end{subcorollary}
Corollary \eqref{phi dim properties}, applied to $\vf$ and ${\wt\vf}$, implies
that $N$ has finite projective dimension over $\vf$ if and only if it
does over ${\wt\vf}$. However, unlike injective dimension, the \emph{value}
of $\pd_{\vf}N$ may depend on the map.

\begin{subexample}
\label{rod:projdim}
Let $R$ be a field.  Set $S=R[x]_{(x)}$ and let $\vf\colon R\to S$ be
the canonical injection.  Set ${\wt S}=R(x)$ and let ${\wt\vf}\colon R\to {\wt S}$
be the canonical injection.  Decomposition into partial fractions shows
that $S$ and ${\wt S}$ have the same rank as vector spaces over $R$, namely,
$c=\max\{\aleph_0,\card(R)\}$.  Pick an $R$-vector space $N$ of rank $c$.
Choosing $R$-linear isomorphisms $S\cong N$ and ${\wt S}\cong N$, endow $N$
with structures of free module of rank $1$ over $S$ and over ${\wt S}$;
both actions restrict to the original action of $R$ on $N$.  Directly,
or from \eqref{dim cohen}, one gets
 \[
\pd_{\vf}N=-\edim S=-1\ne 0=-\edim {\wt S}=\pd_{{\wt\vf}}N
 \] \end{subexample}

The module in the example is \emph{not} an $S$-${\wt S}$-bimodule.  There is
a reason:

\begin{subtheorem}
\label{bimodules:pd}
If the actions of $S$ and ${\wt S}$ on $N$ commute, then
 \[
\pd_{\vf} N = \pd_{{\wt\vf}} N
 \]
 \end{subtheorem}

\begin{proof}
Comparing the degrees of the Laurent series in Proposition
\eqref{bimodules:series}, one gets
 \[
(\pd_{\vf} N+\edim\vf)+\edim{\wt\vf}\le (\pd_{{\wt\vf}} N+\edim{\wt\vf})+\edim\vf
 \]
Thus, $\pd_{\vf}N\le\pd_{{\wt\vf}}N$; the converse inequality holds by
symmetry.
 \end{proof}

\subsection{Depth and Krull dimension}
\label{dim:definition}
The definition of depth of $N$ over $S$ is recalled in \eqref{depth}.
When $\pd_{\vf}N$ is finite, it appears in the following equality of
Auslander-Buchsbaum type, proved in \cite[(4.3)]{IS}:
\[
\depth_{S}N=\depth R-\pd_{\vf}N
\]

Foxby \cite[(3.5)]{Fo} defines the (\emph{Krull}) \emph{dimension} of
$N$ over $S$ to be
 \[
\dim_S N = \sup\{\dim_S\HH_n(N)-n\mid n\in\BZ\}
 \]
Clearly, this number specializes to the usual concept when $N$ is an
$S$-module.

The module $N$ in Example \eqref{rod:projdim} has $\dim_SN=\depth_SN=1$
and $\dim_{{\wt S}}N=\depth_{{\wt S}}N=0$.  This explains the interest of the
following theorem.  Unlike most results in this paper, it is \emph{not}
about invariants over a homomorphism; in particular, we are not assuming
that $Q$ is an $R$-algebra.

\begin{subtheorem}
\label{rod:depth}
Let $Q$ be a local ring such that $N$ has a structure of homologically
finite complex of $Q$-modules.  If the actions of $Q$ and $S$ on $N$
commute, then
 \[
\dim_QN = \dim_SN
\qquad\text{and}\qquad
\depth_QN = \depth_SN
 \]
\end{subtheorem}

\begin{proof}
We may assume $\HH_n(N)=H\ne 0$ for some $n\in\BZ$.  Let $q$ be the
characteristic of the residue field of $Q$, and $p$ be that of $l$. We
claim $q=p$.  There is nothing to prove unless $p$ or $q$ is positive,
and then we may assume $q>0$.  As $V=H/\fn H$ is a finite $Q$-module,
and $q$ is contained in the maximal ideal of $Q$, Nakayama's Lemma
yields $qV\ne V$.  Since $V$ is a $l$-vector space, this implies $qV=0$,
hence $q=p$.

Set $P=\BZ_{(p)}$ and let $\eta\colon P\to Q$ and $\eta'\colon P\to S$
be the canonical maps.

First, we verify the assertion on depths: since the actions of $P$ on $N$
through $\eta$ and through $\eta'$ coincide, $\pd_{\eta}N=\pd_{\eta'}N$ by
Theorem \eqref{bimodules:pd}.  The local ring $P$ is regular, so $\fd_PN$
is finite, and hence, by \eqref{phi dim properties}, both $\pd_{\eta}N$
and $\pd_{\eta'}N$ are finite.  The Auslander-Buchsbaum
formula, see above, yields
 \[
\depth_QN=\depth R-\pd_{\eta}N=\depth R-\pd_{\eta'}N=\depth_SN
 \]

Now we turn to dimensions. It suffices to check that $\dim_Q\HH_n(N)=
\dim_S\HH_n(N)$ for each $n\in\BZ$, so we may assume that $N$ is a module,
concentrated in degree zero.  The actions of $Q$ and $S$ on $N$ commute,
so one obtains a homomorphism of rings
 \[
\tau\colon Q\otimes_PS\lra\Hom_P(N,N)
 \quad\text{where}\quad
(q\otimes s)\longmapsto\big(n\mapsto qsn\big)
 \]
Set $U=\Image(\tau)$: This is a commutative subring of $\Hom_P(N,N)$,
so $N$ has a natural $U$-module structure.  The actions of $U$ and $Q$
commute, so we have inclusions
 \[
U\subseteq\Hom_Q(N,N)\subseteq \Hom_P(N,N)
 \]
where the first one is $Q$-linear.  Since $N$ is a finite $Q$-module,
the same is true of $\Hom_Q(N,N)$, and hence of $U$.  Therefore, the
composed homomorphism of rings $Q\to Q\otimes_PS\to U$ is module finite.
Thus, $U$ is noetherian and $\dim_Q N = \dim_U N$.  By symmetry, $\dim_S
N = \dim_U N$, hence $\dim_Q N=\dim_SN$, as desired.
 \end{proof}

The preceding theorem allows us to complement a result in \cite{IS}.
In that paper a notion of Gorenstein dimension of $N$ over $\vf$,
denoted $\gd_\vf N$, is defined and studied.  In particular, it
is proved in \cite[(7.1)]{IS} that $\gd_\vf N$ and $\gd_{{\wt\vf}} N$ are
simultaneously finite.  However, they can differ, as shown by an example
in \cite[(7.2)]{IS}; there, as in Example \eqref{rod:projdim}, the actions
of $S$ and ${\wt S}$ do not commute.

\begin{subremark}
\label{bimodules:gd}
If the actions of $S$ and ${\wt S}$ on $N$ commute, then
 \[
\gd_{\vf} N = \gd_{{\wt\vf}} N
 \]

Indeed, it is clear from the preceding discussion that we may assume
both Gorenstein dimensions are finite.  In that case, we have
 \[
\gd_{\vf} N=\depth R-\depth_SN=\depth R-\depth_{{\wt S}}N=\gd_{{\wt\vf}} N
 \]
where the equalities on both ends are given by \cite[(3.5)]{IS}, and the
one in the middle comes from Theorem \eqref{rod:depth}.
 \end{subremark}

\subsection{Poincar\'e series}
Now we turn to the proof of Proposition \eqref{bimodules:series}.
The argument hinges on the following construction.

\begin{subconstruction}
\label{square} Each complex $Y$ of modules over the ring
$T=S\otimes_R{\wt S}$ defines a commutative diagram of homomorphisms of
rings
\begin{equation*}
\xymatrixrowsep{3pc}
\xymatrixcolsep{2.5pc}
\xymatrix{
\Ext^0_{T}(Y,Y) \ar@{->}[rrr]
\ar@{->}[ddd] &&& \Ext^0_S(Y,Y) \ar@{->}[ddd]
\\
&
T
\ar@{->}[ul]_-{\eta_{T}}
&S
\ar@{->}[ur]^-{\eta_{S}}
\ar@{->}[l]
\\
&
{\wt S}
\ar@{->}[u]
\ar@{->}[dl]_-{\eta_{{\wt S}}}
&
R
\ar@{->}[u]_{\vf}
\ar@{->}[l]^{{\wt\vf}}
\ar@{->}[dr]^-{\eta_{R}}
\\
\Ext^0_{{\wt S}}(Y,Y) \ar@{->}[rrr] &&& \Ext^0_{R}(Y,Y) }
\end{equation*}
where the slanted arrows land into central subrings.

It is assembled is follows.  The inner square is the canonical
commutative diagram associated with a tensor product of $R$-algebras.
The outer square is obtained from it by functoriality, hence it commutes.
The slanted arrows refer to the maps described in \eqref{Annihilators},
so they are central and the trapezoids commute.

Let $\tau\colon T\to\Ext^0_{R}(Y,Y)$ be the homomorphism of provided by
the diagram above.  Setting $U=\Image(\tau)$, one obtains a commutative
diagram
\begin{equation*}
\xymatrixrowsep{3pc}
\xymatrixcolsep{3pc}
\xymatrix{
R
\ar@{->}[r]^-{{\wt\vf}}
\ar@{->}[d]_{\vf}
&
{\wt S}
\ar@{->}[d]^{\sigma'}
\\
S
\ar@{->}[r]_-{\sigma}
&
U
}
\end{equation*}
of homomorphisms of commutative rings.
 \end{subconstruction}

\begin{sublemma}
\label{radical}
With the notation of Construction \eqref{square}, let $\fr$ denote the
Jacobson radical of the ring $U$, set $V_n=\Tor^R_n(k,Y)$, respectively,
$V_n=\Ext^{-n}_R(k,Y)$ for all $n\in\BZ$, and assume that $Y$ is
homologically finite over $S$ and over ${\wt S}$.
\begin{enumerate}[\rm\quad(1)]
\item
The ring $U$ is finite as an $S$-module, $\fr=\rad(\fn U)$, and
 \[
\inf\{j\in\BN\var\fr^j\subseteq\fn U\}=v<\infty
 \]
\item
For $u=\ell_S(U/\fr)$ the following inequalities hold:
 \[
\ell_T(V_n)
\le\ell_S(V_n)
\le u\ell_T(V_n)
 \]
\item[\rm(3)]
If $\ell_S(V_j)$ is finite for all $j\in\BZ$, then for each ${\wt x}\in{\wt\fn}$ 
the $S$-modules ${\wt V}_n=\Tor^R_n(k,\koszul{{\wt x}}Y)$, respectively,
${\wt V}_n=\Ext^{-n}_R(k,\koszul{{\wt x}}Y)$, satisfy 
 \[
\frac1{v}\big(\ell_S(V_n)+\ell_S(V_{n-1})\big)
\le\ell_S({\wt V}_n)
\le\ell_S(V_n)+\ell_S(V_{n-1})
 \]
\end{enumerate}
 \end{sublemma}

\begin{proof}
(1) The commutative diagram in Construction \eqref{square} shows
that the maps 
 \[
T\to\Ext^0_{S}(Y,Y)\to\Ext^0_{R}(Y,Y)
 \]
are $S$-linear. Thus, $U$ is an $S$-submodule of the image of
$\Ext^0_{S}(Y,Y)$. The last $S$-module is finite by \eqref{newfinite},
hence so is $U$. The Going-up Theorem now implies $\fr=\rad(\fn U)$.
Since $U/\fn U$ is a finite $l$-algebra, its Jacobson radical $\fr/\fn U$
is nilpotent, hence $\fr^v\subseteq\fn$ for some $v$.

(2) It follows from \eqref{homotopies} that $T$ acts on $V_n$ through
the surjective map $\tau$.  Since for every $U$-module $V$ one has
$\ell_T(V)=\ell_U(V)$, it suffices to prove the inequalities
 \[
\ell_S(U/\fr)<\infty
 \qquad\text{and}\qquad
\ell_U(V)\le\ell_S(V)\le\ell_U(V)\cdot\ell_S(U/\fr)
 \]
The first one is an immediate consequence of (1).  When $V$ is simple $\fr
V=0$ and $\ell_U(V)=1$, so the last two inequalities are clear in this case.
The general case follows by computing lengths of subquotients of the
filtration $\{\fr^iV\}_{i\ge0}$.

(3) For this argument we identify $\koszul{{\wt x}}Y$ with the mapping cone of
$\lambda_{{\wt x}}^Y$, see \eqref{cone koszul}, and set $W=\dtensor kRY$. This
leads to an equality of morphisms
\[
\dtensor kR{\lambda_{{\wt x}}^Y}=\lambda_{{\wt x}}^W\colon W\lra W
\]
and to an identification of the complex $\dtensor
kR{\koszul{{\wt x}}Y}$ with the mapping cone of $\lambda_{{\wt x}}^W$.
For $y=\sigma'({\wt x})\in U$ we now obtain
 \[
\HH(\lambda_{{\wt x}}^W)=\lambda_{{\wt x}}^{\HH(W)}= \lambda_{y}^{\HH(W)}
 \]
Applied first with ${\wt S}$ in place of $S$, Part (1) yields
$\fr=\rad({\wt\fn}U)$, hence $y\in\fr$; applied then to $S$, it gives
$y^v\in\fn U$.  By \eqref{finite betti}, the ideal $\fn$ annihilates
$V_n$, so the equalities above yield $\HH(\lambda_{{\wt x}}^W)^v=0$.
The desired result now follows from Lemma \eqref{rod:sandwich}.
 \end{proof}

\begin{proof}[Proof of Proposition \emph{\eqref{bimodules:series}}]
We only give the proof for Poincar\'e series.

Let $\bsx=\{x_1,\dots,x_{e}\}$ and ${\wt\bsx}=\{{\wt x}_1,\dots,{\wt x}_{{g}}\}$
be minimal sets of generators of $\fn$ modulo $\fm S$ and ${\wt\fn}$
modulo $\fm {\wt S}$, respectively.  For any finite set $\bsy$ in $T$,
set
 \[
\lP S{\bsy}=
\sum_{n\in\BZ}\ell_S\big(\Tor^R_n(k,\koszul{\bsy}N)\big)\,t^n
 \]
By definition, one has  $\lP{S}{\bsx}=\rP{\vf}N$ and
$\lP{{\wt S}}{{\wt\bsx}}=\rP{\wt\vf}N$, so the desired formula follows
from the chain of inequalites
\begin{align*}
\lP{S}{\bsx}\cdot(1+t)^{{g}} &\preccurlyeq
\lP{S}{\bsx\sqcup\{{\wt x}_1\}}\cdot(1+t)^{{g}-1}\,v_1
\preccurlyeq\cdots\\
&\preccurlyeq
\lP{S}{\bsx\sqcup{\wt\bsx}}\cdot v_1\cdots v_{{g}}\\
&\preccurlyeq
\lP{T}{\bsx\sqcup{\wt\bsx}}\cdot v_1\cdots v_{{g}}\cdot u\\
&\preccurlyeq
\lP{{\wt S}}{\bsx\sqcup{\wt\bsx}}\cdot v_1\cdots v_{{g}}\cdot u\\
&\preccurlyeq
\lP{{\wt S}}{\bsx\smallsetminus\{x_e\}\sqcup{\wt\bsx}}\cdot(1+t)\,v_1\cdots
v_{{g}}\cdot u
\preccurlyeq\cdots\\
&\preccurlyeq\lP{{\wt S}}{{\wt\bsx}}\cdot (1+t)^{e}v_1\cdots v_{{g}}\cdot u
\end{align*}
obtained as follows.  The inequalities in the first two rows come from
the left hand side of the sandwich in Lemma (\ref{radical}.3),
applied successively to the complexes of $S$-modules
$\koszul{\bsx\sqcup\{{\wt x}_1,\dots,{\wt x}_{j}\}}N$ for $j=0,\dots,{g}-1$.
The next pair of inequalities are provided by Lemma (\ref{radical}.2),
applied to the complex $\koszul{\bsx\sqcup{\wt\bsx}}N$ considered first
over $S$, then over ${\wt S}$.  The final string of $e$ inequalities
is obtained from the right hand side of the sandwich in Lemma
(\ref{radical}.3), applied successively to the complexes of ${\wt S}$-modules
$\koszul{\{x_1,\dots,x_{e-j}\}\sqcup{\wt\bsx}}N$ for $j=1,\dots,e$.
 \end{proof}

\section{Composition}
\label{Composition}

In this section we study how the sequences of Betti numbers and Bass
numbers over $\vf$ react to various changes of rings. Most of the results
take the form of coefficientwise equalities and inequalities involving
their generating functions.  This format is well adopted to study the
finiteness and the asymptotic behavior of the Betti numbers and the Bass
numbers.  In particular, it does not depend on the choice of scale used
to measure asymptotic growth: this fact might acquire importance should
future investigations discover modules or complexes whose Betti numbers
have superpolynomial, but subexponential rates of growth.

Throughout the section we fix a local homomorphism
 \[
\psi\colon (Q,\fl,h)\to(R,\fm,k)
 \]

\subsection{Upper bounds for compositions}
We bound (injective) complexities and curvatures over $\vf\circ\psi$
in terms of invariants over the maps $\psi$ and $\vf$.

\begin{subtheorem}
\label{ass composition}
The following relations hold.
\begin{enumerate}[{\quad\rm(1)}]
\item
The following inequalities hold.
\begin{align*}
\cxy_{\vf\circ\psi}N &\leq \cxy_{\vf}N + \cxy_{\psi}R\\
\curv_{\vf\circ\psi}N &\leq \max\{\curv_{\vf}N, \curv_{\psi}R \}
\end{align*}
\item
If $\HH(N)$ is finite over $R$, then
 \[
\cxy_{\vf\circ\psi}N = \cxy_{\psi}N
\quad\text{and}\quad
\curv_{\vf\circ\psi}N = \curv_{\psi}N
 \]
\end{enumerate}
In addition, the relations obtained by replacing complexities and
curvatures of $N$ with their injective counterparts also hold.
\end{subtheorem}

\begin{proof}
In view of \eqref{polyfactor}, Part (1) is a consequence of
Proposition \eqref{transitivity} below.  Part (2) follows from Theorem
\eqref{bimodules:asymptotes}, applied to $\psi$ and $\vf\circ\psi$.
 \end{proof}

We introduce notation that will be used in several arguments.

\begin{subchunk}
\label{comp:notation}
 Let $\bsu$ be a minimal generating set of
$\fm$ modulo $\fl R$ and $\bsx$ a minimal generating set of $\fn$
modulo $\fm S$, and set $\bsy = \vf(\bsu)\sqcup \bsx$.  The
following number is non-negative:
 \[
d=\edim\psi+\edim\vf - \edim(\vf\circ\psi)
 \]
Indeed, the isomorphism $\fn/\fm S \cong (\fn/\fl S)/(\fm S/\fl
S)$ implies that the set $\bsy$ generates $\fn/\fl S$. Thus,
$\edim(\vf\circ\psi)\leq \edim\vf+\edim\psi$; that is to say,
$d\geq 0$.
 \end{subchunk}

\begin{subproposition}
\label{transitivity} The Poincar\'e series of $\psi$, $\vf$, and
$\vf\circ\psi$ satisfy
 \[
\rP{\vf\circ\psi}N\cdot (1+t)^d \preccurlyeq
\rP{\vf}N\cdot\rP{\psi}R
  \]
\end{subproposition}

\begin{proof}
In the derived category of $S$, the isomorphism
  \[
\koszul\bsy N\simeq\koszul\bsu R\otimes_R\koszul\bsx N
  \]
combined with the associativity formula for derived tensor products yields
  \[
\dtensor{\left(\dtensor{h}Q{\koszul\bsu R}\right)}R{\koszul\bsx N}
\simeq \dtensor {h}Q{\koszul\bsy N}
  \]
Thus, one has a standard spectral sequence
  \[
E_{pq}^2=\Tor_p^R\left(\Tor_q^Q(h,\koszul{\bsu}R),\koszul{\bsx}N\right)
\implies \Tor_{p+q}^Q\left(h,\koszul{\bsy}N\right).
  \]
As the $R$-module $\Tor^Q(h,\koszul{\bsw}R)$ is annihilated by $\fm$,
one has an isomorphism
  \[
\Tor_p^R\left(\Tor_q^Q(h,\koszul{\bsu}R),\koszul{\bsx}N\right)
  \cong
\Tor_p^R(k,\koszul{\bsx}N) \otimes_k \Tor_q^Q(h,\koszul{\bsu}R)
  \]
The desired inequality is a formal consequence of this isomorphism
and the convergence of the spectral sequence; see the proof of
Theorem \eqref{ceiling}.
 \end{proof}

\subsection{Upper bound for right factors}
When $\psi$ is surjective and $\vf=\idmap^R$ the following theorem
specializes to a known result for $\rP{R}N$, see \cite[(3.3.2)]{barca}.

\begin{subproposition}
\label{golod bound}
If $\psi\colon Q\to R$ is a local homomorphism, then
 \[
\rP{\vf}N \preccurlyeq {\rP{\vf\circ\psi}N}
\cdot\frac 1 {1-t\big(\rP{\psi}R-1\big)} \cdot {(1+t)^{\edim R +\edim\vf - \edim S}}
 \]
A similar inequality also holds for Bass series: it is obtained from the
one above by replacing $\rP{\vf}N$ and ${\rP{\vf\circ\psi}N}$ with
$\rI{\vf}N$ and ${\rI{\vf\circ\psi}N}$, respectively.
 \end{subproposition}

\begin{Remark}
As $\rP{\psi}R$ is a formal power series, the formal power series
$1-t\big(\rP{\psi}R-1\big)$ has an inverse, which is itself a formal
power series. It is given by the formula
 \[
\frac{1}{1-t\big(\rP{\psi}R-1\big)} =\sum_{i=0}^\infty
t^i\big(\rP{\psi}R-1\big)^i
 \]
which shows that it has nonnegative integer coefficients and constant
term $1$.
 \end{Remark}

\begin{proof}
In view of Theorem \eqref{duality}, it suffices to deal with Poincar\'e
series.

Let $E$ be a DG algebra resolution of $h$ over $Q$, see
\cite[(2.1.10)]{barca}. Form the DG algebra $A=E\otimes_Q\koszul \bsu R$
and the DG $A$-module $Y=A\otimes_R\koszul{\bsx}N$.  The augmentations
$E\to h$ and $\koszul{\bsu}R\to k$ are morphisms of DG algebras, where
$h$ and $k$ are concentrated in degree $0$.  They yield a morphism of
DG algebras $A\to k$.

Note the standard isomorphisms in the derived category of $S$-modules:
  \[
\dtensor kAY \simeq \dtensor kA{\big(\dtensor AR{\koszul \bsx N}\big)}
              \simeq \dtensor kR{\koszul\bsx N}.
  \]
For each $n$ they induce $l$-linear isomorphisms $\Tor^A_n(k,Y)\cong
\Tor^R_n(k,\koszul{\bsx}N)$, hence
  \[ \sum_{n\in\BZ}
\rank_l\Tor_n^A(k,Y)\,t^n = \rP{\vf}N
  \] By the choice of $E$, there
is an isomorphism $\HH_n(A) \cong \Tor^Q_n(h,\koszul{\bsu}R)$ and hence
  \[
\sum_{n\in\BZ} \rank_k \HH_n(A)\, t^n = \rP{\psi}R
  \]
Likewise, $\HH_n(Y) \cong \Tor^Q_n(h,\koszul{\bsy}N)$, as $Y\cong
E\otimes_Q\koszul\bsy N$, so we get
  \[
\sum_{n\in\BZ} \rank_l \HH_n(Y)\, t^n = \rP{\vf\circ\psi}N \cdot (1+t)^e
  \]
from Proposition \eqref{pseries:extra}. Thus, the inequality we seek
translates to
 \begin{equation}
\sum_{n\in\BZ} \rank_l \Tor_n^A(k,Y)\, t^n \preccurlyeq
\frac{\sum_{n\in\BZ}\rank_l\HH_n(Y)\,t^n}
{1-t\big(\sum_{n\in\BZ}\rank_k\HH_n(A)\,t^n-1\big)} \tag{$*$}
 \end{equation}

Now, as $\HH_n(A)=0$ for $n<0$ and $\HH_n(Y)=0$ for $n\ll 0$,
there exists a strongly convergent Eilenberg-Moore spectral sequence
 \[
E_{pq}^2= \Tor_p^{\HH(A)}(k,\HH(Y))_q \implies \Tor_{p+q}^A(k,Y)
 \]
see \cite[Chapter 7]{Mc}.  In the light of this, there are (in)equalities
 \[
\rank_l\Tor_n^A(k,Y)\leq \sum_{p+q=n} \rank_l E_{pq}^2
= \sum_{p+q=n} \rank_l \Tor_p^{\HH(A)}(k,\HH(Y))_q
 \]
They can be rewritten as an inequality of formal power series
 \[
\sum_{n\in\BZ} \rank_l\Tor_n^A(k,Y)\,t^n \preccurlyeq
\sum_{n\in\BZ}\bigg(\sum_{p+q=n}
\rank_l \Tor_p^{\HH(A)}(k,\HH(Y))_q\bigg)t^n
 \]
By using a standard resolution of $k$ over $\HH(A)$, one can compute
the graded $l$-vector space $\Tor_p^{\HH(A)}(k,\HH(Y))$ as the $p$th
homology of a complex of the form
 \[
\cdots \lra \HH_{\ges 1}(A)^{\otimes n}\otimes_k{\HH(Y)}
\lra \cdots \lra \HH_{\ges 1}(A)\otimes_k\HH(Y) \lra \HH(Y) \lra 0
 \]
There is thus an inequality of formal Laurent series
 \[
\sum_{n\in\BZ}\bigg(\sum_{p+q=n}\rank_l
\Tor_p^{\HH(A)}(k,\HH(Y))_q\bigg)t^n
\preccurlyeq
\frac{\sum_{n\in\BZ}\rank_l\HH_n(Y)\,t^n}
{1-t\big(\sum_{n\in\BZ}\rank_k\HH_n(A)\,t^n-1\big)}
 \]
The desired inequality $(*)$ is obtained by combining the inequalities
above.
 \end{proof}

\begin{subcorollary}
\label{finite curv}
The curvature and the injective curvature of $N$ over $\vf$ are finite.
 \end{subcorollary}

\begin{proof}
By Remark \eqref{cx duality} it suffices to deal with $\curv_\vf
N$. Let $p$ denote the characteristic of $k$, set $Q=\BZ_{(p)}$ and
let $\psi\colon Q\to R$ be the structure map. Since $Q$ is regular,
both $\rP{\psi}N$ and $\rP{\vf\circ\psi}N$ are Laurent polynomials,
see \eqref{dim properties}.  Thus, the coefficientwise upper bound for
$\rP{\vf}N$ in Proposition \eqref{golod bound} is a formal Laurent series
that represents a rational function.  Now apply (\ref{polyfactor}.1)
and (\ref{polyfactor}.3).
 \end{proof}

\subsection{Complete intersection homomorphisms}
\label{CI homomorphisms}
When $\psi$ has  a `nice' property the invariants of $N$ over $\vf\circ\psi$
do not differ much from those over $\vf$.

\begin{subchunk}
\label{ci maps}
As in \cite{Av:lci}, we say that $\psi$ is \emph{complete intersection}
at $\fm$ if in some Cohen factorization $Q\to (Q',\fl',k)\xra{\psi'}R$ of
$\grave\psi$ the ideal $\Ker\psi'$ is generated by a regular set $\bsf'$.
When this is the case, the Comparison Theorem \cite[(1.2)]{AFH} shows
that $\Ker\psi''$ is generated by a regular set whenever $\psi''\ddot\psi$
is a Cohen factorization of $\grave\psi$.  If $\psi$ is weakly regular at
$\fm$, see \eqref{weakly reg}, then it is complete intersection at $\fm$.
 \end{subchunk}

\begin{subtheorem}
\label{ass ci}
If $\psi$ is complete intersection at $\fm$, then
\begin{gather*}
\cxy_{\vf\circ\psi}N
\leq\cxy_{\vf}N\leq\cxy_{\vf\circ\psi}N +\dim Q-\dim R+\edim\psi \\
\curv_{\vf\circ\psi}N \le \curv_\vf N
\le\max\{\curv_{\vf\circ\psi}N,1\}
\end{gather*}
Analogous inequalities also hold for $\injcxy$ and $\injcurv$.
\end{subtheorem}

The next result complements Lemma \eqref{superficial:prop1}.

\begin{sublemma}
\label{superficial:prop2} If $f\in Q$ is a regular element and
$R=Q/Qf$, then
 \[
\rP{\vf\circ\psi}N\preccurlyeq\rP{\vf}N\cdot(1+t)
\qquad\text{and}\qquad \rP{\vf}N \preccurlyeq {\rP{\vf\circ\psi}
N}\cdot\frac1{(1-t^2)}
 \]
where $\psi\colon Q\to R$ is the canonical map. If, moreover, $f$
is in $\fl\Ker(\vf\circ\psi)$, then
 \[
\rP{\vf}N=\rP{\vf\circ\psi} N\cdot\frac1{(1-t^2)}
 \]
\end{sublemma}

\begin{proof}
Since $\rP{\psi}R=1+t$, the inequality on the left follows from
Proposition \eqref{transitivity} and the one on the right from Theorem
\eqref{golod bound}.  The proof of \cite[(3.3.5.2)]{barca} applies
verbatim to yield the equality.
 \end{proof}

\begin{proof}[Proof of Theorem \emph{\eqref{ass ci}}]
By Theorem \eqref{duality} we restrict to projective invariants.

In view of \eqref{cx cohen} we may assume that $R$ and $S$ are
complete. We use Construction \eqref{big-diagram} and adopt the notation
introduced there. By \eqref{ci maps} the kernel of the homomorphism
$\psi'$ is generated by a regular set $\bsf'$; clearly, one has
 \[
\card\bsf'=\dim Q'-\dim R=(\dim Q+\edim\psi)-\dim R
 \]
Since $\psi''=Q''\otimes_{Q'}\psi'$, and $Q''$ is faithfully flat
over $Q'$, the image $\bsf''$ of $\bsf'$ in $Q''$ is a regular set of
$\card\bsf'$ elements, and generates the kernel of $\psi''$.

The complexity and the curvature of $N$ over $\vf\circ\psi$ equal those
over $\vf\circ\psi'$ by Theorem \eqref{flat ass}. Thus, to prove the
theorem we may assume that $\psi$ is surjective with kernel generated by
a regular set. In this case the desired assertions result from repeated
applications of the preceding lemma through \eqref{polyfactor}.
 \end{proof}

\subsection{Flat base change}
We consider complexes of $S$-modules induced from $R$.

\begin{subtheorem}
\label{base change}
If $\vf$ is flat and $S/\fm S$ is artinian, then for every homologically
finite complex of $R$-modules $M$ the following equalities hold:
 \[
\cxy_{\vf\circ\psi}(S\otimes_RM)=\cxy_{\psi}M \quad\text{and}\quad
\curv_{\vf\circ\psi}(S\otimes_RM)=\curv_{\psi}M
 \]
\end{subtheorem}

As usual, this follows from more precise relations involving Poincar\'e
series.

 \begin{subproposition}
\label{series:extended}
Let $M$ be a homologically finite complex of $R$-modules.

If $\vf$ is flat and $\fn^v\subseteq \fm S$ for some integer $v$, then
 \[
\rP{\psi}{M}\cdot{v^{-e}}(1+t)^e\preccurlyeq
\rP{\vf\circ\psi}{M\otimes_RS}\cdot u^{-1}(1+t)^d \preccurlyeq
\rP{\psi}{M}\cdot(1+t)^e
 \]
where $u=\ell(S/\fm S)$, $e =\edim{\vf}$, and $d=\edim R+e-\edim
S$.
 \end{subproposition}

\begin{proof}
Let $\bsu$ be a minimal set of generators of $\fm$ modulo $\fl R$
and set
 \[
X = \dtensor{h}Q{\koszul{\vf(\bsu)}{S\otimes_RM}}
 \]
The isomorphism $\koszul{\vf(\bsu)}{(S\otimes_RM)}\simeq
\koszul{\bsu}M\otimes_RS$ of complexes of $S$-modules and the
flatness of $S$ over $R$ yield isomorphisms
 \begin{align*}
\HH(X)\cong\Tor^Q(h,\koszul{\bsu}M)\otimes_RS\cong
\Tor^Q(h,\koszul{\bsu}M)\otimes_k(S/\fm S)
 \end{align*}
They produce the following equality of Poincar\'e series
 \[
\rP{\psi}M \cdot u=\sum_{n\in\BZ}\ell_S{\HH_n(X)}\,t^n
 \]
Let $\bsx$ be a minimal set of generators of $\fn$ modulo $\fm$, and
note that the set $\bsy=\vf(\bsu)\sqcup\bsx$ generates $\fn$
modulo $\fl S$. From our hypothesis and \eqref{hann koszul} we get
 \[
\fn^v\subseteq\fm S=\fl S+\bsu S\subseteq\hann S{X}
 \]
Iterated applications of Lemma \eqref{rod:sandwich} now give
 \[
\sum_{n\in\BZ}\ell_S{\HH_n(X)}\,t^n\cdot v^{-e}(1+t)^e\preccurlyeq
\sum_{n\in\BZ} \ell_S \HH_n(\koszul{\bsx}X)\,t^n
\preccurlyeq\sum_{n\in\BZ}\ell_S{\HH_n(X)}\,t^n\cdot(1+t)^e
 \]
The isomorphism $\koszul{\bsx}X\simeq\dtensor{h}Q{\koszul{\bsy}{S\otimes_RM}}$
and Proposition \eqref{pseries:extra} now yield
 \[
\sum_{n\in\BZ} \ell_S \HH_n(\koszul{\bsx}X)\,t^n =
\rP{\vf\circ\psi}{S\otimes_RM}(1+t)^d
 \]
Putting together the formulas above we obtain the desired inequalities.
 \end{proof}

\section{Localization}
\label{Localization}

A basic and elementary result asserts that Betti numbers of finite modules
over local rings do not go up under localization:  If $M$ is a finite
$R$-module and $\fp$ is a prime ideal of $R$, then for each $n\in\BZ$
there is an inequality $\betti n{R_\fp}{M_\fp}\le \betti nRM$.

A notion of localization is also available, and has been systematically
used, for complexes over local homomorphisms.  Namely, for each prime
ideal $\fq$ of $S$ the complex $N_\fq$ is homologically finite over the
induced local homomorphism $\vf_\fq\colon R_\fp\to S_\fq$, where $\fp=
\fq\cap R$.  However, $\betti n{\vf_\fq}{N_\fq}$ may exceed $\betti
n{\vf}N$, even when $\vf$ is flat.

\begin{example}
  Let $R$ be a one-dimensional local domain $R$ whose completion $\wh R$ has a minimal
  prime ideal $\fq$, such $\wh R_\fq$ is not a field; see Ferrand and Raynaud
  \cite[(3.1)]{FR}.  Set $S=\wh R$, let $\vf\colon R\to S$ be the completion map, and note
  that $\fq\cap R=(0)$; thus, $R_{(0)}$ is a field.  For the $ S$-module $N=S$ 
Proposition (\ref{regular source}) then yields 
 \[
\betti 1{\vf_\fq}{N_\fq}=\kappa^{S_\fq}_1(N_\fq)>0= \betti 1{\vf}{N}
 \]
  \end{example}

Nor is there an analog over maps of the inequality
$\fd_{R_\fp}M_\fp\le\fd_RM$.

\begin{example}
\label{dim localization}
Let $R$ be a field, let $S$ be a Cohen-Macaulay ring of positive
dimension, and let $\fq$ be a non-maximal prime of $S$.  In view of
\eqref{depth}, one then has
\[
\pd_{\vf_\fq}\!{S_\fq}=-\depth S_\fq>-\depth S=\pd_{\vf}{S}
\]
 \end{example}

Thus, it is noteworthy that asymptotic invariants over homomorphisms
localize as expected.  Unlike the corresponding result for complexity and
curvature over rings, the theorem below needs a fairly involved argument.

\begin{theorem}
\label{ass localization}
For every prime ideal $\fq$ in $S$ there are inequalities
 \[
\cxy_{\vf_\fq}{N_\fq} \leq \cxy_{\vf}{N} \quad\text{and}\quad
\curv_{\vf_\fq}{N_\fq} \leq \curv_{\vf}{N}
 \]

In particular, if $\pd_{\vf}{N}$ is finite, then so is
 $\pd_{\vf_\fq}\!{N_\fq}$.
 \end{theorem}

\begin{Remark}
If $D$ is a normalized dualizing complex for $S$, then an appropriate
shift of $D_\fq$ is a normalized dualizing complex for $S_\fq$.
Using this fact together with Theorem \eqref{duality}, one sees that
the preceding theorem has a counterpart for injective complexities and
curvatures, \emph{provided} $S$ has a dualizing complex.  We do not know
whether the last condition is necessary.
 \end{Remark}

Once again, the theorem follows from a more precise result on
Poincar\'e series.

\begin{proposition}\label{localization}
For every prime ideal $\fq$ in $S$  there
exists a polynomial $q(t)\in\BZ[t]$ with non-negative coefficients such
that the following inequality holds
 \[
\rP{\vf_\fq}{N_\fq}\preccurlyeq q(t)\cdot\rP{\vf}{N}
 \]
\end{proposition}

\begin{proof}
  Form a commutative diagram of local homomorphisms
\begin{equation*}
  \xymatrixrowsep{3pc} 
  \xymatrixcolsep{3pc} 
  \xymatrix{ R
    \ar@{->}[r]^{\dot\vf} \ar@{->}[d]_{\vf} \ar@{->}[dr]^{\grave\vf}
    &R' \ar@{->}[d]^{\vf'} \\ S \ar@{->}[r]_{\sigma} &\wh S }
\end{equation*}
where $\sigma$ is the completion map of $S$ in the $\fn$-adic
topology and the upper triangle is a minimal Cohen factorization of
$\grave\vf$. By faithful flatness, choose $\wt\fq$ in $\Spec\wh S$ so that
$\wt\fq\cap S=\fq$ and $\ell\big(\wh S_{\wt\fq}/\fq S_{\wt\fq}\big)=u
<\infty$. Setting $\wt\fp=\wt\fq\cap R'$ and $\fp=\wt\fq\cap R$ and
localizing, one gets a commutative diagram of local homomorphisms
\begin{equation*}
\begin{gathered}
  \xymatrixrowsep{3pc}
  \xymatrixcolsep{3pc} 
   \xymatrix{ R_\fp
    \ar@{->}[r]^{\dot\vf_{\wt\fp}} \ar@{->}[d]_{\vf_\fq}
    \ar@{->}[dr]^{\grave\vf_{\wt\fq}} &R'_{\wt\fp}
    \ar@{->}[d]^{\vf'_{\wt\fq}} \\ S_\fq \ar@{->}[r]_{\sigma_{\wt\fq}}
    &\wh S_{\wt\fq} }
\end{gathered}
\qquad\text{that we redraw as}\qquad
\begin{gathered}
  \xymatrixrowsep{3pc} 
  \xymatrixcolsep{3pc} 
  \xymatrix{ P
    \ar@{->}[r]^{\dot\pi} \ar@{->}[d]_{\varkappa} \ar@{->}[dr]^{\pi} &P'
    \ar@{->}[d]^{\pi'} \\ Q \ar@{->}[r]_{\tau} &T }
\end{gathered}
 \end{equation*}
The homomorphism $\tau$ is flat with artinian closed fiber.  Indeed,
$\tau$ is a localization of the flat homomorphism $\sigma$, and its
closed fiber is the ring $\wh S_{\wt\fq}/\fq S_{\wt\fq}$.

Set $L=N_\fq$. Proposition \eqref{series:extended} yields the inequality
 \[
\rP{\varkappa}{L}\cdot(1+t)^s \preccurlyeq
\rP{\tau\circ\varkappa}{L\otimes_QT}\cdot(1+t)^e(v^s/u)
 \]
where $s$, $e$, $u$, and $v$ are non-negative integers.
As $\tau\circ\varkappa=\pi=\pi'\circ\dot\pi$, one has
 \[
\rP{\tau\circ\varkappa}{L\otimes_QT}= \rP{\pi'\circ\dot\pi}{L\otimes_QT}
 \]
Set $p(t)=\rP{\dot\pi}{P'}$. Since $\dot\pi$ is flat, $p(t)$ is a
polynomial: see, for example, \eqref{phi dim properties}.  Thus,
Proposition \eqref{transitivity} provides a non-negative integer $g$
such that
 \[
\rP{\pi'\dot\pi}{L\otimes_QT}\cdot(1+t)^g
\preccurlyeq\rP{\pi'}{L\otimes_QT}\cdot p(t)
 \]
Next we note that over the ring $T=\wh S_{\wt\fq}$ there are isomorphisms
of complexes
 \[
L\otimes_QT\cong(N\otimes_S Q)\otimes_{Q}T\cong N\otimes_ST\cong
(N\otimes_S\wh S)\otimes_{\wh S}T=\wh N\otimes_{\wh S}T\cong\wh N_{\wt\fq}
 \]
They give the first equality in the sequence below, where the second one
comes from \eqref{surjective}, the last from Remark \eqref{series:cohen},
and the inequality is classical:
 \[
\rP{\pi'}{L\otimes_QT}=\rP{\pi'}{\wh N_{\wt\fq}}=\rP{R_{\wt\fp}}{\wh
N_{\wt\fq}}\preccurlyeq\rP{R'}{\wh N}=\rP{\vf}{N}
 \]
Putting together the comparisons above, one obtains that
 \[
\rP{\varkappa}{L}\cdot(1+t)^{(s+g)}\preccurlyeq \rP{\vf}N\cdot p(t)(v^s/u)
 \]
The inequality we seek is contained in the one above, because
$1\preccurlyeq(1+t)^{(s+g)}$.
 \end{proof}

\section{Extremality}
\label{Extremality}

A recurrent theme in local algebra is that the homological properties of the $R$-module
$k$ carry a lot of information on the structure of $R$.  Thus, modules or complexes
`homologically similar' to $k$ provide test objects for properties of $R$.  We focus on
rates of growth of Betti numbers, which controls regularity and the complete intersection
property.  As in \cite{Av:extremal}, we say that a finite $R$-module is
(\emph{injectively}) \emph{extremal} if its (injective) complexity and curvature are equal
to those of $k$.

Proposition \eqref{dependencies} shows that these same numbers are also
upper bounds for (injective) complexities and curvatures over $\vf$.
This leads to obvious extensions of the notions of extremality.
One reason to study them may not be obvious a priori: it provides many new
classes of extremal modules over $R$, see Remark \eqref{extremal modules}.

We say that the complex $N$ is \emph{extremal over $\vf$} if
 \[
\cxy_\vf N=\cxy_Rk\quad\text{and}\quad \curv_\vf N=\curv_Rk
 \]
It is \emph{injectively extremal over $\vf$} if
 \[
\injcxy_\vf N=\cxy_Rk\quad\text{and}\quad \injcurv_\vf N=\curv_Rk
 \]
Using $\injcxy_Rk$ and $\injcurv_Rk$ to define injective extremality
yields the same result: computing with a minimal free resolution of $k$,
one gets $\bass nRk=\betti nRk$.

The homological characterization of complete intersections in Theorem
\eqref{generalcichar} shows that the notion of extremality has two
distinct aspects.

\begin{remark}
\label{dichotomy}
If $R$ is complete intersection, then $N$ is (injectively) extremal over
$\vf$ if and only if its (injective) complexity is equal to $\codim R$.

If $R$ is not complete intersection, then $N$ is (injectively)
extremal over $\vf$ if and only if its (injective) curvature is equal
to $\curv_Rk$.
 \end{remark}

\begin{remark}
\label{extremality:trivial}
The definition of separation and \eqref{polyfactor} yield: If $N$
is (injectively) separated over $\vf$ and $\HH(N)\ne 0$, then it is
(injectively) extremal over $\vf$.
  \end{remark}

Separation is a much stronger condition than extremality:

\begin{example}
Assume that $R$ is not regular and set $\vf=\idmap^R$.

(1) If $R$ is complete intersection of codimension $\ge2$, and $M_n$
is the $n$th syzygy of an extremal $R$-module $M$, then it is clear
that each $M_n$ is extremal; however, it follows from \cite[(6.2)]{AGP}
that $M_n$ can be separated for at most one value of $n$.

(2) Remark \eqref{separated:sums} below shows that $k\oplus R$ is extremal, but not
separated.
 \end{example}

\begin{remark}
\label{extremal modules}
When $N$ is homologically finite over $R$, its Betti numbers and Bass
numbers may be easier to compute over $\vf$ than over $R$.  On the
other hand, Theorem \eqref{ass:finite} shows that the corresponding
asymptotic invariants over $R$ and over $\vf$ coincide.  This can be used
to identify new classes of extremal complexes over $R$.  For example,
Remark \eqref{extremality:trivial} shows that if the local homomorphism
$\vf\colon R\to S$ is module finite, the ring $S$ is regular, and
$\HH(N)\ne0$, then $N$ is extremal over $R$.
 \end{remark}

\begin{remark}
\label{separated:sums}
Let $X$ be a homologically finite complex of $S$-modules.  The complex
of $S$-modules $N\oplus X$ is separated over $\vf$ if and only if
both $N$ and $X$ are separated, while it is extremal over $\vf$ if and
only if one of $N$ or $X$ is extremal.

Indeed, the claim on extremality is clear.  The claim on separation is
verified by using Proposition \eqref{ceiling} along with the equalities
 \[
K^S_{N\oplus X}(t) = K^S_{N}(t) + K^S_{X}(t)
 \quad\text{and}\quad
\rP{\vf}{N\oplus X} =\rP{\vf}{N} + \rP{\vf}{X}
\]
 \end{remark}

The main result in this section quantifies and significantly generalizes
a theorem of Koh and Lee, see \cite[(2.6.i)]{KL}.  The idea to use socles
to locate nonzero homology classes is inspired by an argument in their
proof of \cite[(1.2.i)]{KL}.

Let $L$ be an $S$-module. Recall that its \emph{socle} is the
$S$-submodule
 \[
\Soc_S(L) = \{a\in L\mid \fn a=0\}
 \]
Note that $\Soc_R(L)$ is also defined; it is an $S$-submodule of $L$
and contains $\Soc_S(L)$.

\begin{theorem}
\label{socR}
Let $L$ be a module. If $\bsv$ is an $L$-regular set in $S$ such that
\begin{enumerate}[{\quad\rm(a)}]
\item
$\Soc_S(L/\bsv L)\nsubseteq \fm (L/\bsv L)$, or
\item
$\Soc_R(L/\bsv L)\nsubseteq \fm (L/\bsv L)$ and the ring $S/\fm S$
is artinian,
\end{enumerate}
then $L$ is extremal and injectively extremal over $\vf$.
 \end{theorem}

The proof is based on a simple sufficient condition for extremality.

\begin{lemma}
\label{extremal:complex}
If $X$ is a homologically finite complex of $S$-modules such that
 \[
0<\ell_S\big(\HH(X)/\fm\HH(X)\big)<\infty
\quad\text{and}\quad
\Soc_R\big(\Ker(\dd^X_j)\big)\not\subseteq
\fm X_j + \dd(X_{j+1})
 \]
for some $j\in\BZ$, then for all $n\in\BZ$ the following inequalities hold
 \[
\ell_S\big(\Tor^R_n(k,X)\big)\ge \betti{n-j}Rk
 \quad\text{and}\quad
\ell_S\big(\Ext^n_R(k,X)\big)\ge \betti{n+j}Rk
 \]
In particular, the complex $X$ is extremal and injectively extremal
over $\phi^i$.
 \end{lemma}

\begin{proof}
Set $X'_n=X_n$ for all $n\ne j$ and $X'_j=\fm X_j+\dd(X_{j+1})$.  It
is easy to see that $X'$ is a subcomplex of $X$.  We form a
commutative diagram
\begin{equation*}
\xymatrixrowsep{3pc}
\xymatrixcolsep{1pc}
\xymatrix{
{\shift^j\Soc_R\big(\Ker(\dd^X_j)\big)}
\ar@{=}[r]
&V
\ar@{->}[d]_{\zeta}
\ar@{->}[rrr]^{\pi}
&&& W \ar@{->}[d]^{\iota}
\ar@{=}[r]&{\displaystyle\frac{V}{V\cap X'}}
\\
&X
\ar@{->}[rrr]^{\rho}
&&& Y\ar@{=}[r]&{\displaystyle\frac{X}{X'}\qquad}
}
\end{equation*}
of complexes of $S$-modules, where $V$, $W$, and $Y$ are concentrated
in degree $j$. The ring $R$ acts on them
through $k$, so $\Tor^R(k,-)$ applied to the diagram above yields a
commutative diagram of graded $S$-modules
\begin{equation*}
\xymatrixrowsep{3pc}
\xymatrixcolsep{6pc}
\xymatrix{
\Tor^R(k,k)\otimes_k V
\ar@{->}[d]_{\Tor^R(k,\zeta)}
\ar@{->}[r]^-{\Tor^R(k,k)\otimes_k\pi}
&
\Tor^R(k,k)\otimes_k W
\ar@{->}[d]^{\Tor^R(k,k)\otimes_k\iota}
\\
\Tor^R(k,X) \ar@{->}[r]^{\Tor^R(k,\rho)} & \Tor^R(k,k)\otimes_kY}
\end{equation*}
By construction, the map $\pi$ is surjective and the map $\iota$ is injective, so the
image of $\Tor^R(k,\rho)$ contains an isomorphic copy of $\Tor^R(k,k)\otimes_kW$.  Thus,
 \[
\ell_S\big(\Tor^R_n(k,X)\big)\ge\betti {n-j}Rk\cdot \ell_S(W)
 \]
for all $n\in\BZ$.  Since $W\ne 0$, by hypothesis, extremality follows from Remark
\eqref{ass:finite_remark}.

Similar arguments yield the statements concerning injective invariants.
 \end{proof}

\begin{proof}[Proof of Theorem \emph{\eqref{socR}}]
Corollary (\ref{cxy koszul.2}) shows that $L$ and $L/\bsv L$ are
extremal simultaneously, so replacing $L$ with $L/\bsv L$, we
assume one of the conditions:
\begin{enumerate}[{\quad\rm(a)}]
\item
$\Soc_S(L)\not\subseteq\fm L$, or
\item
$\Soc_R(L)\not\subseteq\fm L$ and $S/\fm S$ is artinian.
 \end{enumerate}

(a) Let $\bsy$ be a finite set of $s$ generators of $\fn$, and
set $X=\koszul\bsy L$.  Note that $\rank_l\HH(X)$ is finite, that
$\Ker(\dd^X_s) = \Soc_S(L)$, and that $X_{s+1}=0$.  Thus, Lemma
\eqref{extremal:complex} shows that $X$ is extremal and Proposition
(\ref{cxy koszul.1}) completes the proof.

(b) Apply Lemma \eqref{extremal:complex} to the complex $X=L$.
 \end{proof}

We wish to compare the hypotheses of various theorems yielding
extremality.

\begin{example}
Let $R=S=k[x]/(x^3)$ where $k$ is a field of characteristic $2$,
let $\phi_R$ be the Frobenius endomorphism of $R$, and set $L=S$.
One then has
 \[
\Soc_R(L)=(x)\,,\qquad \Soc_S(L)=(x^2)=\fm L
 \qquad\text{and}
\qquad\serre_SL=4\,,
 \]
see Remark \eqref{useless}.  Thus, Theorem (\ref{socR}.b) shows that $S$
is extremal over $\vf$, but neither Theorem (\ref{socR}.a) nor Theorem
\eqref{hell:separated} can be applied.
  \end{example}

\begin{remark}
The conditions `extremal' and `injectively extremal' are independent
in general, even over $\vf=\idmap^R$.  For instance, Example \eqref{trivial
extension} yields
\begin{gather*}
  \cxy_R k = \cxy_R E = \injcxy_R R = \infty \\
  \curv_Rk = \curv_RE = \injcurv_RR   = 2
\end{gather*}
for $R=k[x,y]/(x^2,xy,y^2)$ and the injective hull $E$ of $k$ over $R$.
Thus, $E$ is extremal but not injectively extremal, while $R$ is
injectively extremal, but not extremal.  \end{remark}

\section{Endomorphisms}
\label{Endomorphisms}

Let $\phi$ be a local endomorphism, that is, a local homomorphism
 \[
\phi\colon (R,\fm,k) \to (R,\fm,k)
 \]
and let $N$ be a homologically finite complex of $R$-modules.

\subsection{Contractions}
We say that the homomorphism $\phi$ is \emph{contracting}, or that
$\phi$ is a \emph{contraction}, if for each element $x$ in $\fm$
the sequence $(\phi^i(x))_{i\ges1}$ converges to $0$ in the $\fm$-adic
topology.  Observe that $\phi$ is contracting if $\phi^j$ is contracting
for some integer $j\geq 1$, and only if $\phi^j$ is contracting for
every integer $j\geq 1$.

One motivation for considering contracting homomorphisms comes from

\begin{subexample}
\label{archetype}
The archetypal contraction is the \emph{Frobenius endomorphism} $\phi_R$
of a ring $R$ of prime characteristic $p$, defined by $\phi_R(x)=x^p$
for all $x\in R$.
 \end{subexample}

Interesting contractions are found in every characteristic:

\begin{subexample}
Let $k$ be a field, let $G$ be an additive semigroup without torsion, let
$k[G]$ denote the semigroup ring of $G$ over $k$, and let $\fp$ denote
the maximal ideal of $k[G]$ generated by the elements of $G$.  For each
non-negative integer $s$, the endomorphism $\sigma$ of the semigroup $G$,
given by $\sigma(g)= s\cdot g$, defines an endomorphism $k[\sigma]$
of the ring $k[G]$.  It satisfies $k[\sigma](\fp)\subseteq\fp$, and
so it induces an endomorphism $\phi$ of the local ring $R=k[G]_\fp$.
The endomorphism $\phi$ is contracting when $s\geq 2$.
 \end{subexample}

Many complexes are separated over contractions.  Theorem
\eqref{hell:separated} implies the next result, which uses the homotopy
Loewy length $\hell SN$ defined in \eqref{loewy}.

\begin{subtheorem}
\label{contra artin}
Assume $\phi$ is a contraction and $\HH(N)\ne0$.  If $j\ge1$ and $q\ge2$
satisfy $\phi^j(\fm)\subseteq\fm^q$, then $N$ is separated and injectively
separated over $\phi^i$ for all
\[
i\ge j\log_q(\hell S{\rkoszul SN})
\]
In particular, $N$ is extremal and injectively extremal over $\phi^i$. \qed
 \end{subtheorem}

The following alternative description of contractions shows that the
numbers $j$ and $q$ in the hypothesis of the theorem always exist,
and gives bounds for them.  Example \eqref{ci example} shows that these
bounds cannot be improved in general.

\begin{sublemma}
The map $\phi$ is contracting if and only if $\phi^{\edim R}
(\fm)\subseteq\fm^2$.
 \end{sublemma}

\begin{proof}
The `if' part is clear, so we assume $\phi$ is contracting.

Let $\delta\colon\fm/\fm^2\to\fm/\fm^2$ be the map induced by $\phi$. It
is a homomorphism of abelian groups, which defines on $\fm/\fm^2$
a filtration
 \[
0=\Ker(\delta^0)\subseteq\Ker(\delta^1)\subseteq
\cdots\subseteq\Ker(\delta^i)\subseteq\Ker(\delta^{i+1})\subseteq\cdots
 \]
We need to prove that for $r=\edim R$ one has $\Ker(\delta^r)=\fm/\fm^2$,
that is, $\delta^r=0$.

Each subgroup $\Ker(\delta^i)$ is a $k$-vector subspace of
$\fm/\fm^2$, by a direct verification.  We show next that
if $\Ker(\delta^i)=\Ker(\delta^{i+1})$ for some $i\ge1$, then
$\Ker(\delta^i)=\Ker(\delta^{j})$ for all $j\ge i$.  By induction,
it suffices to do it for $j=i+2$.  If $\delta^{i+2}(x)=0$, then
 \[
\delta^{i+2}(x)=\delta^{i+1}\big(\delta(x)\big)=0
 \]
implies $\delta(x)\in\Ker(\delta^{i+1})=\Ker(\delta^{i})$.  Therefore,
$\delta^{i+1}(x)=0$, as desired.

Since $\rank_k\fm/\fm^2=r$, the properties of $\Ker(\delta^i)$ that
we have established imply equalities $\Ker(\delta^j)=\Ker(\delta^r)$
for all $j\ge r$, that is,
 \[
\Ker(\delta^r)=\bigcup_{i=0}^\infty\Ker(\delta^i)
 \]
Our hypothesis means that the union above is all of $\fm/\fm^2$.
 \end{proof}

We can do better than Theorem \eqref{contra artin} when $N=R$, at least
for curvature.

\begin{subtheorem}
\label{endomorphism}
If $\phi\colon R\to R$ is a contracting endomorphism, then
 \[
\curv_{\phi}R =\curv_R k\qquad
 \]

If, in addition, the ring $R$ is Gorenstein, then
 \[
\injcurv_{\phi}R=\curv_Rk
 \]
 \end{subtheorem}

\begin{proof}
By the preceding theorem, there is an integer $j$ giving the equality
below:
 \[
\curv_Rk = \curv_{\phi^j}R \leq \curv_{\phi}R \leq \curv_R k
 \]
The inequalities come from Theorem (\ref{ass composition}.1) and from
Proposition (\ref{dependencies}.5).

If $R$ is Gorenstein, then $R^\dagger\simeq R$, see \eqref{dualizing
complex}, hence the middle equality below:
 \[
\injcurv_{\phi}R = \curv_{\phi}R^\dagger=\curv_\phi R = \curv_Rk
 \]
The other two are given by Remark \eqref{cx duality} and the first part
of the theorem.
 \end{proof}

Theorem \eqref{endomorphism} has no counterpart for complexities.

\begin{subexample}
\label{ci example} Let $k$ be a field, set
$R=k[[x_1,\dots,x_r]]/(x_1^2,\dots,x_r^2)$, and let $\phi\colon
R\to R$ be the $k$-algebra homomorphism given by
 \[
\phi(x_1)=0\quad\text{and}\quad \phi(x_i)=x_{i-1}
\quad\text{for}\quad 2\leq i\leq r
 \]
Then $\cxy_{\phi^j}R = \min\{j,r\}$ for each integer $j\ge0$.

Indeed, fix an integer $j\ge1$ and form the ring
 \[
S=k[[x_{j+1},\dots,x_r]]/(x_{j+1}^2,\dots,x_r^2)
 \]
The map $\phi^j$ factors as $R\xra{\pi} S\xra{\iota} R$, where $\pi$ is
the canonical surjection with kernel $(x_{1},\dots,x_j)$ and $\iota$ is
the $k$-algebra homomorphism with $\iota(x_i)=x_{i-j}$ for $j<i\leq r$.
As $\iota$ is flat and $R$ is artinian, Theorem \eqref{base change}
yields $\cxy_{\phi^j}R = \cxy_{\pi}S$. Set
 \[
Q=k[[x_1,\dots,x_r]]/(x_{j+1}^2,\dots,x_r^2)
 \]
 and let $\psi\colon Q\to R$ be the canonical surjection with kernel
 $(x_1^2,\dots,x_j^2)$. With $m=\min\{j,r\}$, repeated application of Lemma
 \eqref{superficial:prop2} yields the first equality in the sequence
 \[
\rP{\pi}{S} = \frac{\rP{\pi\circ\psi}{S}}{(1-t^2)^m} =
\frac{\rP{Q}{S}}{(1-t^2)^m}= \frac{(1+t)^m}{(1-t^2)^m}=\frac1{(1-t)^m}
 \]
 The second equality by Remark \eqref{surjective}; the third holds because
 $\koszul{\{x_1,\dots,x_m\}}Q$ is a minimal resolution of $S$ over $Q$.  Now invoke
 \eqref{polyfactor}.
 \end{subexample}

\subsection{Frobenius endomorphisms}
\label{frobenius}
In this subsection $R$ has prime characteristic $p$ and $\phi_R\colon R\to
R$ is its Frobenius endomorphism.  This is a contracting endomorphism,
so the results from the preceding subsection apply. One noteworthy
additional feature is that they can be interpreted entirely in terms of
classically defined invariants.  Indeed, Theorem \eqref{ass:finite} 
validates the following

\begin{subremark}
Let $\up{i}N$ denote $N$ as an $R$-$R$-bimodule with the left action 
through $\phi_R^i$ and the right action the usual one.
Each $\Tor^R_n(k,\up{i}N)$ is an $R$-module where the action of $R$ 
is induced from the right action of $R$ on $\up{i}N$;  
by (\ref{hann2}) and Lemma \eqref{newfinite} it has finite length. 

For each integer $i\ge 1$ the following equalities hold:
\begin{gather*}
\cxy_{\phi_R^i}N=
\inf\left\{d \in\BN \left|
\begin{gathered}
\text{there exists a number $c\in\BR$ such that}\\
\ell_R \Tor^R_n(k,\up{i}N)\le c n^{d-1}\text{ for all $n\gg0$}
\end{gathered}
\right\}\right.
\\
\curv_{\phi_R^i} N=
\limsup_n\sqrt[n]{\ell_R \Tor^R_n(k,\up{i}N)}
\end{gather*}
\end{subremark}

This subsection is organized around the

\begin{subquestion}
Is $N$ separated (respectively, extremal) over $\phi_R^i$ for all $i\ge1$?
 \end{subquestion}

A lot of evidence points to a positive answer.

\begin{subremark}
When $R$ is not complete intersection, Theorem \eqref{contra artin}
shows that $N$ is separated and injectively separated over $\phi^i_R$
for all $i\ge\log_p(\hell RN)$.

When $R$ is complete intersection $\rP Rk=(1+t)^d/(1-t)^{c}$, see
(\ref{background}.2), so the next theorem asserts that $N$ is separated
and injectively separated over $\phi^i_R$.
 \end{subremark}

\begin{subtheorem}
\label{frob:ci}
If $R$ is complete intersection, $d=\dim R$, and $c=\codim R$, then for
each $i\ge1$ the following equalities hold:
 \[
\rP{\phi^i_R}N=K^R_N(t)\cdot\frac{(1+t)^d}{(1-t)^{c}}
\qquad\text{and}\qquad
\rI{\phi^i_R}N=K^R_N(t)\cdot\frac{(1+t)^dt^c}{(1-t)^{c}}
 \]
In particular, $N$ is extremal and injectively extremal over $\phi^i_R$.
 \end{subtheorem}

\begin{subcorollary}
\label{extremal frob}
For every ring $R$ of positive characteristic, the module $R$ is extremal
over $\phi_R^i$ for each integer $i\geq 1$.
 \end{subcorollary}

\begin{proof}
In view of Remark \eqref{dichotomy}, Theorem \eqref{endomorphism}
establishes the assertion when $R$ is not complete intersection.  When it
is, Theorem \eqref{frob:ci} applies.
 \end{proof}

The proof of the theorem is given at the end of this section.  We approach
it through an explicit description of minimal Cohen factorizations of
powers of the Frobenius endomorphisms of complete rings.

\begin{subconstruction}
\label{frob factorization}
Let $\bsv=\{v_1,\dots,v_r\}$ be a minimal set of generators for $\fm$.
Identifying $R$ with its image in $\wh R$ under the completion map,
note that $\bsv$ minimally generates the maximal ideal of $\wh R$.

Let $\bsx=\{x_1,\dots,x_r\}$ be a set of formal indeterminates and
set $Q=k[[\bsx]]$.  With the standard abbreviation $\bsx^J=x_1^{j_1}\cdots
x_r^{j_r}$ for $J=(j_1,\dots,j_r)$, each $g\in Q$ has a unique expression
$g=\sum_{J\in\BN^r}a_J\bsx^J$ with $a_J\in k$.  For each positive integer
$q$ set
 \[
g^{[q]}=\sum_{J\in\BN^r}a_J^q\bsx^J
 \]

Choose, by Cohen's Structure theorem, a surjective homomorphism
$\psi\colon Q\to\wh R$, such that $\psi(x_j)=v_j$ for $j=1,\dots,r$.
Let $\bsf=\{f_1,\dots,f_c\}$ be a minimal generating set of $\Ker\psi$.
The choices made so far imply $\bsf\subseteq(\bsx)^2$.  Set
 \[
\bsf^{q}=\{f^{q}_1,\dots,f^{q}_c\}
 \qquad\text{and}\qquad
\bsf^{[q]}=\{f^{[q]}_1,\dots,f^{[q]}_c\}
 \]

Fix an integer $i\ge1$, set $q=p^i$ and $\phi=\phi_{\wh R}^i$.  Let
$\bsy=\{y_1,\dots,y_r\}$ denote a second family of formal indeterminates.
With the data above, form the diagram
\begin{equation*}
\xymatrixrowsep{3.3pc}
\xymatrixcolsep{1.2pc}
\xymatrix{
k[[\bsx]]
\ar@{=}[r]
&Q
\ar@{->}[rrrrr]^-{\phi^i_Q}
\ar@{->}[dd]_-{\psi}
&&&&&
Q
\ar@{=}[r]
\ar@{->}[d]^-{\psi'}
& k[[\bsx]]
\\
&&
{\displaystyle\frac{k[[\bsx]]}{(\bsf^{[q]})}}[[\bsy]]
\ar@{=}[r]
&R'
\ar@{->}[drrr]^-{\phi'}
\ar@{->}[rrr]^-{\rho}
&&&R''
\ar@{=}[r]
\ar@{->}[d]^-{\phi''}
& {\displaystyle\frac{k[[\bsx]]}{(\bsf^{q})}}
\\
{\displaystyle\frac{k[[\bsx]]}{(\bsf)}}
\ar@{=}[r]
&\wh R
\ar@{->}[rrrrr]^-{\phi}
\ar@{->}[urr]^-{\dot\phi}
&&&&& \wh R
\ar@{=}[r]
&{\displaystyle\frac{k[[\bsx]]}{(\bsf)}}
}
\end{equation*}
of local homomorphisms, where the new objects are defined as follows.

The rings $R'$ and $R''$ are described by the respective equalities.

The maps $\psi'$ and $\phi''$ are the canonical surjections.

The maps $\phi^i_Q$ and $\dot\phi$ are given by the formulas
 \[
\phi^i_Q(g)=g^{q}
\qquad\text{and}\qquad
\dot\phi\big(g+(\bsf)\big)=g^{[q]}+(\bsf^{[q]})
 \]

The map $\rho$ is the unique homomorphism of complete $k$-algebras
satisfying
 \[
\rho\big(x_i+(\bsf^{[q]}))=x_i^q+(\bsf^{q})
\qquad\text{and}\qquad
\rho(y_i)=x_i+(\bsf^{q})
 \]
for $i=1,\dots,r$; note that $\rho$ is surjective.

The map $\phi'$ is the composition $\phi''\circ\rho$.

The definitions above show that the diagram commutes.
 \end{subconstruction}

\begin{subproposition}
\label{ultimate}
The maps in Construction \eqref{frob factorization} have the
properties below.
\begin{enumerate}[\quad\rm(1)]
\item
$\phi'\dot\phi$ is a minimal Cohen factorization of $\phi$.
\item
There is an equality of formal Laurent series
 \[
\rP{\phi}N=\rP{\phi''}{\wh N}\cdot(1+t)^r
 \]
\end{enumerate}
\end{subproposition}

\noindent\emph{Proof}.  (1) Since $\dot\phi(x_i)=x_i$ for $i=1,\dots,r$,
the ring $R'/\fm R'$ is isomorphic to the regular ring $k[[\bsy]]$.
To prove that $\dot\phi$ is flat, consider the composition
 \[
Q=k[[\bsx]]\xra{\phi^i_Q}k[[\bsx]]\xra{\iota}k[[\bsx,\bsy]]=Q'
 \]
where $\iota$ is the natural inclusion.  Thus, $\iota\circ\phi^i_Q$
maps the $Q$-regular set $\bsx$ to the $Q'$-regular set $\bsx^q$.
Computing $\Tor^Q(k,Q')$ from the resolution $\koszul{\bsx}Q$ of $k$
over $Q$, one gets $\betti nQ{Q'}=0$ for $n>0$.  It follows that $Q'$
is flat over $Q$, see \eqref{hdim:classical}.  As $\dot\phi$ is obtained
from $\iota\circ\phi^i_Q$ by base change along $\psi$, we conclude that
$\dot\phi$ is flat, as desired.

As $\phi'$ is surjective, $\phi=\phi'\dot\phi$ is a Cohen factorization.

(2)  The kernel of $\rho$ is generated by the set
$\{x_1-y_1^q,\dots,x_r-y_r^q\}$, which is regular and superficial.
Remark \eqref{series:cohen} and Theorem \eqref{superficial:prop1}
yield
 \begin{xxalignat}{3}
&\phantom{square}
 &\rP{\phi}N=\rP{\phi'}{\wh N}&=\rP{\phi''}{\wh N}\cdot(1+t)^r
&&\square
\end{xxalignat}

\begin{proof}[Proof of Theorem \emph{\ref{frob:ci}}]
We use Construction \eqref{frob factorization}.  As $R$ is complete
intersection, the set $\bsf$ is regular, see \eqref{ci rings}.
It follows that so is $\bsf^q$; note that $\bsf^q$ is contained
in $\fl(\bsf Q)$, where $\fl$ is the maximal ideal of $Q$.  Lemma
(\ref{ultimate}.2), Proposition \eqref{superficial:prop2}, and Example
(\ref{separated:examples}) provide the equalities below
 \begin{align*}
\rP{\phi}N&=\rP{\phi''}{\wh N}\cdot(1+t)^r\\
&=\rP{\phi''\circ\psi'}{\wh N}\cdot\frac{1}{(1-t^2)^c}\cdot(1+t)^r\\
&=\rP Qk\cdot\frac{K^R_N(t)}{(1+t)^r}\cdot\frac{(1+t)^d}{(1-t)^{c}}\\
&=K^R_N(t)\cdot\frac{(1+t)^d}{(1-t)^{c}}
 \end{align*}

Theorem \eqref{duality} yields the desired expression for $\rI{\vf}N$.
 \end{proof}

\section{Local homomorphisms}
\label{Local homomorphisms}

In this section we use the techniques and results developed earlier in
the paper to study relations between the ring theoretical properties of
$R$ and $S$ and the homological properties of the $R$-module $S$.
First we look at descent problems.

\begin{theorem}
Let $\vf\colon R\to S$ be a local homomorphism and let $N$ be a
homologically finite complex of $S$-modules with $\HH(N)\ne0$.
\begin{enumerate}[\rm\quad(1)]
\item
If $\fd_RN<\infty$ and $S$ is regular, then so is $R$.
\item
If $\curv_{\vf}S\leq 1$ and $S$ is complete intersection,
then so is $R$.
\end{enumerate}
\end{theorem}

\begin{Remark}
Part (1) of the theorem is due to Apassov \cite[Theorem R]{Ap1}.  Part (2)
significantly generalizes \cite[(5.10)]{Av:lci}, where it is proved that
maps of finite flat dimension descend the complete intersection property.
Dwyer, Greenlees, and Iyengar \cite{DGI} show that this property descends
even under the weaker hypothesis $\curv_{\vf}N\leq 1$.
 \end{Remark}

\begin{proof}
(1) One has $\betti n{\vf}N=0$ for $n\gg0$, see \eqref{koszul:sup}.
As $N$ is separated by Corollary \eqref{separated:regular}, the equality
in \eqref{separated complexes} implies $\betti nRk=0$ for all $n\gg0$.

(2)  Applying Theorem (\ref{ass composition}.1) to the composition
$R\to S\xra{=}S$ and the $S$-module $l$, we obtain the first inequality below:
 \[
\curv_{\vf}l \leq \max\{\curv_{\vf}S,\curv_{\idmap^S}l\}
             \le \max\{\curv_{\vf}S,1\}
             \leq 1
 \]
The second one comes from Theorem \eqref{generalcichar}, the third
from our hypothesis.  We conclude that $R$ is complete intersection by
referring once more to \eqref{generalcichar}.
 \end{proof}

Next we extend to arbitrary local homomorphisms a characterization of
Gorenstein rings, due to Peskine and Szpiro in the case of surjective
maps.

\begin{theorem}
\label{peskineszpiro}
For a local homomorphism $\vf\colon R\to S$ the condition 
$\id_RS<\infty$ holds if and only if 
$\fd_RS<\infty$ and the ring $R$ is Gorenstein.
\end{theorem}

\begin{proof}
Over a Gorenstein ring the flat dimension of a module is finite if and
only if its injective dimension is, see \cite[(2.2)]{LV}, so we have to
prove that $\id_RS<\infty$ implies $R$ is Gorenstein.  Let $R\to R'\to \wh
S$ be a minimal Cohen factorization of $\grave\vf$.  Corollary \eqref{dim
properties} gives $\id_{R'}\wh S<\infty$.  As $R'\to\wh S$ is surjective,
$R'$ is Gorenstein by Peskine and Szpiro \cite[(II.5.5)]{PS}.  By flat
descent, see \cite[(23.4)]{Ma}, so is $R$.
 \end{proof}

Finally, we turn to properties of a local ring $R$ equipped with
a contracting endomorphism $\phi$, for instance, a ring of prime
characteristic with its Frobenius endomorphism.  Some of our theorems
are stated in terms of homological properties of the $R$-module
$\up{\phi^i}R$, that is, $R$ viewed as a module over itself through
$\phi^i$.

The prototype of such results is a famous theorem of Kunz, \cite[(2.1)]{Ku}:
A ring $R$ of prime characteristic is regular if $\up{\phi^i}R$ is
flat for some $i\ge1$, only if $\up{\phi^i}R$ is flat for all $i$.
Later, Rodicio \cite[Theorem 2]{Ro} showed that the flatness hypothesis
on $\up{\phi^i}R$ can be replaced by one of finite flat dimension.
Our first criterion extends these results to all contracting endomorphisms
and provides tests for regularity by finite injective dimension, which
are new even for the Frobenius endomorphism

\begin{theorem}
\label{endo regchar}
For a contraction $\phi\colon R\to R$ the following
are equivalent.
\begin{enumerate}[{\quad\rm(i)}]
\item
$R$ is regular.
\item
$\fd_R\up{\phi^i}R=\dim R/(\phi^i(\fm)R)$ for all integers $i\ge1$.
\item
$\fd_R\up{\phi^i}R<\infty$ for some integer $i\geq 1$.
\item
$\id_R\up{\phi^i}R=\dim R$ for all integers $i\ge1$.
\item
$\id_R\up{\phi^i}R<\infty$ for some integer $i\geq 1$.
\end{enumerate}
 \end{theorem}

\begin{proof}
The implication \text{(i) $\implies$ (iv)} holds by Theorem
\eqref{bassroberts}, the implication \text{(iv) $\implies$ (v})
is clear, while \text{(v) $\implies$ (iii)} comes from Theorem
\eqref{peskineszpiro}.

\text{(iii) $\implies$ (i)}.  Since $\phi^i$ is a contraction, from Theorem
\eqref{endomorphism} one gets $\curv_Rk=\curv_{\phi^i}R=0$.  This
means $\pd_Rk$ is finite, that is, $R$ is regular.

\text{(i) $\implies$ (ii)}.  Let $\bsv$ be a minimal set of generators
of $\fm$.
The Koszul complex $\koszul{\bsv}R$ is a free resolution of $k$
over $R$, hence we get
\[
\Tor^R_n(k,\up{\phi^i}R)=\HH_n(\koszul{\bsv}{\up{\phi^i}R})
=\HH_n(\koszul{\phi^i(\bsv)}R)
\]
The largest $n$, such that $\Tor^R_n(k,\up{\phi^i}R)\ne0$, is equal to
$\fd_R(\up{\phi^i}R)$, see \eqref{hdim:classical}.  The largest $n$,
such that $\HH_n(\koszul{\phi^i(\bsv)}R)\ne0$, is equal to $\dim R-g$,
where $g$ is the maximal length of an $R$-regular sequence in the ideal
${\phi^i}(\fm)R$, see \cite[(16.8)]{Ma}.  The ring $R$, being regular, is
Cohen-Macaulay, so referring to \cite[(17.4.i)]{Ma} we conclude
\[
\dim R-g=\dim R-\height(\phi^i(\fm)R)=\dim R/(\phi^i(\fm)R)
\]

The implication \text{(ii) $\implies$ (iii}) is clear.
 \end{proof}

Next we show that the complete intersection property of $R$ can
also be read off of conditions on a contracting endomorphism $\phi$.
They are encoded in the growth of the Betti numbers of the module $R$
over $\phi^i$.  Indeed, the following result is abstracted from Theorem
\eqref{generalcichar}, Theorem \eqref{endomorphism}, and Corollary
\eqref{extremal frob}.

\begin{theorem}
\label{endo cichar}
For a contraction $\phi\colon R\to R$ the following
are equivalent.
\begin{enumerate}[{\quad\rm(i)}]
\item
$R$ is complete intersection.
\item
$\cxy_{\phi^i}R\le\codim R$ for all integers $i\ge1$.
\item
$\cxy_{\phi^i}R<\infty$ some integer $i\ge1$.
\item
$\curv_{\phi^i}R \leq 1$ for some integer $i\geq 1$.
\end{enumerate}
When $R$ has prime characteristic and $\phi$ is its Frobenius map
they are equivalent to
\begin{enumerate}[{\quad\rm(i)}]
\item[\rm(ii)$'$]
$\cxy_{\phi^i}R=\codim R$ for all integers $i\ge1$.
\qed
\end{enumerate}
\end{theorem}

As $R$ is regular if and only if $\codim R=0$, and the flat dimension
of $\up{\phi^i}R$ over $R$ is finite if and only if $\cxy_{\phi^i}R=0$,
the equivalence of conditions (i), (ii)$'$, and (iii) above constitutes
another broad generalization of the theorems of Kunz and Rodicio.
The theorem  also contains a characterization of complete intersections
in terms of Frobenius endomorphisms due to Blanco and Majadas
\cite[Proposition 1]{BM}:

\begin{remark}
Let $\phi$ be a contraction, such that for some $i$ and some Cohen
factorization $\wh R\to R'\to\wh R$ of $\wh\phi^i$ the $R'$-module $\wh R$
has finite CI-dimension in the sense of \cite{AGP}.  By \cite[(5.3)]{AGP}
one then has $\cxy_{R'}\wh R<\infty$, so $\cxy_{\wh\phi^i}\wh R$ is finite
by \eqref{cx cohen}.  Theorem \eqref{endo cichar} now shows that $\wh R$ is
complete intersection, and hence so is $R$.
 \end{remark}

\section*{Acknowledgements}

Part of the work on this paper was done while one of us was collaborating
with Sean Sather-Wagstaff on \cite{IS}, where the focus is on finite
Gorenstein dimension.  We thank Sean for many useful discussions, which
have been beneficial to this project, and for a very thorough reading
of its outcome.  Thanks are also due to Anders Frankild for raising the question
of characterizing regularity in terms of the injective dimension of the
Frobenius endomorphism, settled in Theorem \eqref{endo regchar},
and to the referee for useful suggestions on the exposition.


\end{document}